\documentclass[3p,times,11pt]{elsarticlearxiv}
\usepackage{graphicx}        
\usepackage{multicol}  
\usepackage{easyReview}
\usepackage{mathtools}
\usepackage[mathscr]{euscript}
\usepackage[bottom]{footmisc}
\usepackage{mdframed}
\usepackage[english]{babel}
\usepackage{amssymb,amstext,amsmath,amsthm,mathabx}
\usepackage{hyperref}
\usepackage{array}   
\newcolumntype{L}{>{$}l<{$}}
\theoremstyle{definition}
\newtheorem{definition}{Definition} 
\newtheorem{remark}{Remark} 
 
\newtheorem{theorem}{Theorem}

\newtheorem{example}{Example}

\usepackage{amssymb,amstext,amsmath,mathabx}

\newcommand{\no}[1]{\widebar{#1}}

\def\F{\mathcal F}
\def\P{\mathcal P}
\def\M{\mathcal M}
\def\I{\mathcal I}
\def\H{\mathcal H}
\def\pr{\mathbb{P}}
\def\prev{\mathbb{P}}
\def\K{\mathcal{K}}

\def\G{\mathcal{G}}

\def\V{\mathcal{V}}
\def\C{\mathscr{C}}
\def\C{\mathscr{C}}

\usepackage[utf8]{inputenc}
\usepackage{hyperref}
\usepackage{verbatim}
\usepackage{mathtools}
\usepackage{color}
\usepackage[normalem]{ulem} 
\usepackage{amssymb,amstext,amsmath,wasysym}

\journal{}
\begin{document}
\begin{frontmatter}

\title{
A probabilistic analysis of selected notions of iterated conditioning \\
under coherence
}

\author[lc]{Lydia Castronovo}
\address[lc]{Department of Mathematics and Computer Science, Via Archirafi 34, 90123 Palermo, Italy}
\ead{lydia.castronovo@unipa.it}
\cortext[cor1]{Corresponding author}
\author[lc]{Giuseppe Sanfilippo\corref{cor1}
}
\ead{giuseppe.sanfilippo@unipa.it}

\begin{abstract}
It is well know that basic conditionals satisfy some desirable basic logical and probabilistic properties, such as the compound probability theorem, but checking the validity of these becomes trickier when we switch to compound and iterated conditionals. 
We consider de Finetti’s notion of conditional as a three-valued object and as a conditional random quantity in  the betting framework. We recall the notions of conjunction and disjunction among conditionals in selected trivalent logics. First, in the framework of specific three-valued logics we analyze the notions of iterated conditioning introduced by Cooper-Calabrese, de Finetti and Farrell, respectively. We show that the compound probability theorem  and other  basic properties are not preserved by these objects,  by also computing some probability propagation rules. Then, for each trivalent logic we introduce an iterated conditional as a suitable random quantity which satisfies the compound prevision theorem and some of the desirable properties.  We also check the validity of two   generalized versions of Bayes' Rule for iterated conditionals. 
We study the p-validity of generalized versions of Modus Ponens and two-premise centering for  iterated conditionals.
Finally, we observe that all the basic properties are satisfied only by the  iterated conditional mainly developed in recent papers by Gilio and Sanfilippo in the setting of conditional random quantities. 
\end{abstract}

\end{frontmatter}

\section{Introduction}
The  study of conditionals (typically expressed by ``if–then'' sentences),  compound conditionals (which are obtained by combining conditionals 
with the logical operators such as ``and'', ``or'',  ``not''), and iterated conditionals (where both antecedent and consequent are conditionals) 
is a relevant research topic  in many fields, such as philosophy of science, psychology of uncertain reasoning, probability theory,  and conditional logics.  See for example, 
\cite{adams75,benferhat97,Cant22,CoSV15,crocco1995, Cruz20,Douven19,DuPr94,gilio12ijar,GoNW91,Kauf09,KrLM90,Miln97,NgWa94,OHEHS07,PeVa20,Stal70,vanFraassen76}. Compound and iterated conditionals  are largely used  in natural language to make decisions and inferences based on incomplete or uncertain information. Thus, how  to interpret and combine conditionals  so to  represent 
human rational reasoning realistically 
is a key  question  in the AI community.  Recently there has been a growing interest  in using  conditionals in  the field of AI and of knowledge representation. 
See \cite{KR2022-12,DeDP2020,FGGS22KR,FlGH20,GiSa19,Kern-Isb2023,Kern2023b,mundici21}.

Among the several tools for managing conditional uncertainty, subjective probability allows to  evaluate uncertainty on a conditional  by means of a conditional probability assessment, interpreted as a degree of belief.
In the subjective theory, the probability $P(E)$, that an individual in a given state of uncertain knowledge  attributes to  an event $E$, is a measure of  its degree of belief in the  occurrence  of $E$.
In order to operationally assess consistent probabilities, de Finetti proposed a coherence principle based on a suitable betting scheme
(\cite{definetti31})
.
In this framework, probability is assessed only to the involved  events and  coherence amounts to the avoidance of ``Dutch Book''. All the classical properties of a finitely additive probability follow from the coherence principle. 
This approach is therefore more flexible because it is not 
 necessary to give  probability values to  each event of a given Boolean  algebra.
 The subjective theory has been extended to conditionals  by means of conditional probability assessments.

In this framework,  given two events $A$ and $B$, the  conditional \emph{if $A$ then $B$} is represented not by the material conditional event \emph{not-$A$ or $B$}, but rather by a conditional event designated as $B|A$.  This is a three-valued object which can be assessed as ``true'',  ``false'', or ``void'' (\cite{defi36}). In this formulation, the probability of the conditional \emph{if $A$ then $B$} is interpreted by the 
conditional probability $P(B|A)$.  Subjective conditional probability is primitive, it does not  require an unconditional probability assessment  and hence it can be properly defined even if the conditioning event has probability zero (see Remark \ref{REM:CPT}).  
In particular an (unconditional) event $B$ coincides with the conditional event $B|\Omega$ and hence $P(B)=P(B|\Omega)$.

Trivalent logics have usually been used to combine  conditionals, as found in \cite{adams75,Cala1990,Cooper1968,defi36,EgRS20,Farrell1979}. 
Conjoined and disjoined  conditionals are interpreted in these logics as three-valued objects. However, such interpretations  lead to the invalidity of some basic logical properties.  Then, some basic probabilistic properties are not preserved when logical operations among conditionals belong to a  trivalent logic (\cite{GiSa22,SUM2018S}). In this paper we will study different definitions of iterated conditioning by   showing that  similar problems arise  both when   iterated conditionals are represented as three-valued objects and when they  depend on conjoined conditionals defined in  trivalent logics.
As a consequence, human-like reasoning about conditionals under uncertainty cannot properly be formalized in such terms. 

A different approach to compound and iterated  conditionals  has been followed in \cite{Kauf09} and in \cite{McGe89}, relevant to the definition of conjunction.  
A related study has been  developed in the setting of coherence in the recent papers \cite{GiSa14,GiSa19,GiSa21IJAR,GiSa21A}. In this approach, conjoined, disjoined, and iterated conditionals are defined not as three-valued objects, but  rather 
as suitable conditional random quantities having possibly more than three values in  the unit interval $[0,1]$.
A  betting interpretation governs these objects in a  coherent-setting.
An attractive advantage of this approach is that all the basic algebraic and  probabilistic properties  are preserved.  These include De Morgan’s Laws  and  Fr\'echet-Hoeffding bounds.  For a synthesis see \cite{GiSa21A}. 

A way to build a Boolean algebra of  conditionals satisfying  suitable properties has been introduced in \cite{FlGH20}.
In this context, 
the atomic structure
defines the compound conditionals  at a formal algebraic level.
A general theory of compound conditionals in the framework of  conditional random quantities can be framed in such a structure \cite{FGGS22KR}. 
Then based on a particular extension of a full conditional probability to this  Boolean structure (\cite{FlGH20}), 
some probabilistic properties of compound conditionals
coincide  with results obtained under coherence  in the framework of conditional random quantities (\cite{FGGS23}).


Conjunctions and   disjunctions among conditionals  have been introduced  and  studied quite commonly in  three-valued logics  (\cite{adams75,BaOP13,Cooper1968,Cantwell2008,defi36,EgRS20,Farrell1979, Cala1990,CiDu12,CiDu13,GoNW91}). In particular, de Finetti in 1935 (\cite{defi36}) also proposed  a three-valued logic  for conditional events by  introducing suitably defined notions of conjunction and disjunction.  These coincide with features of Kleen-Lukasiewicz-Heyting logic \cite{CiDu13} (see also \cite{Kleene38}).
Further, Calabrese (\cite{Cala1990}) and Cooper (\cite{Cooper1968}) introduced an algebra of conditionals by using the notions of quasi conjunction and quasi disjunction, similarly studied by Adams  (\cite{adams75}).
In his trivalent logic de Finetti introduced an operation of iterated conditioning called ``subordination'' and denoted here by $|_{dF}$.  This respects the requirement that, among other properties,  the Import-Export principle (\cite{McGe89}) is satisfied \ref{EQ:IMPEXPDF}).
Farrell  also 
introduced an operation  of iterated conditioning (denoted here by $|_F$)  in his trivalent logic (\cite{Farrell1979}), which uses the 
same notions of  conjunction and disjunction as de Finetti.  This also satisfies the Import-Export principle.
Cooper and Calabrese also introduced an operation  of iterated conditioning (denoted here by $|_C$) in their trivalent logic,  which  satisfies   the Import-Export principle as well.  

We should recall that  the validity of such a  principle, jointly with the requirement of preserving  classical probabilistic properties,
leads
to the well-known Lewis' triviality results  \cite{lewis76}.
However, the notion of iterated conditioning studied in the framework of conditional random quantities under coherence avoids the  Lewis’ triviality results because the Import-Export Principle is not satisfied, even satisfying basic probabilistic properties (see \cite{GiSa14,SGOP20,SPOG18}). 
This operation of iterated conditioning, denoted here  by $|_{gs}$, is based on the notion of the conjunction of two conditionals defined as a conditional random quantity ($\wedge_{gs}$) by means  of the following structure (\cite{SGOP20}) $
\Box|\bigcirc = \Box\wedge \bigcirc +\prev(\Box|\bigcirc)\no{\bigcirc}$, where $\prev$ is the symbol of prevision,   $\Box$ and $\bigcirc$ are (indicators of)   conditional events, and $\Box\wedge \bigcirc$ is a conditional random quantity with value in $[0,1]$. 
In the framework of subjective probability, the prevision  $\mu=\prev(\Box|\bigcirc)$  
	represents the amount that you agree to pay, knowing that you will receive the random quantity  $\Box\wedge \bigcirc +\prev(\Box|\bigcirc)\no{\bigcirc}$.
Moreover, 
by the linearity property of a coherent prevision and by observing  that $\no{\bigcirc}=1-\bigcirc$,  from the previous structure it follows that $\prev(\Box|\bigcirc)=\prev(\Box\wedge \bigcirc)+\prev(\Box|\bigcirc)[1-\prev(\bigcirc)]$, that is $\prev(\Box\wedge \bigcirc)=\prev(\Box|\bigcirc)\prev(\bigcirc)$. This last equation 
 is a generalized version of the compound probability theorem, whose standard version  is $P(A\wedge B)=P(B|A)P(A)$. Indeed, 
when $\Box$ and $\bigcirc$ are two events $A$ and $B$,
in a conditional bet on $B|A$,
the  structure above allows  
to numerically interpret  the indicator of the conditional event $B|A$ as the  random win $A\wedge B+P(B|A)\no{A}$, which takes value 1, or 0, or  $P(B|A)$, according to whether $AB$ is true, or $A\no{B}$ is true, or $\no{A}$ is true (\cite{gilio90,Lad95} see also \cite{GiSa14}).
Then, when $\Box$ denotes $B|K$,    $\bigcirc$ denotes $A|H$ (and hence $\no{\bigcirc}$ denotes $\no{A|H}=\no{A}|H$), and  $\Box\wedge \bigcirc$  is the conjunction $(B|K)\wedge_{gs} (A|H)$,  the iterated conditional $(B|K)|_{gs}(A|H)$ is defined as $(B|K)\wedge_{gs} (A|H)+\mu\,(\no{A}|H)$, where $\mu=\prev[(B|K)|_{gs}(A|H)]$. In addition, $\Box$ and $\bigcirc$ can be also replaced by conjoined conditionals and hence the previous structure allows to introduce a more general notion of iterated conditional $\Box|\bigcirc$ as done in \cite{GiSa21ECSQARU}.
The purpose of this paper is to investigate some of the basic properties valid for events and  conditional events with a view to different operations of iterated conditioning. Indeed, 
things get more problematic when we replace events with conditional events and we move to the properties of iterated conditionals. 
We recall four selected notions of conjunction in trivalent logics: Kleene-Lukasiewicz-Heyting-de Finetti ($\wedge_K$), Lukasiewicz ($\wedge_L$), Bochvar-Kleene ($\wedge_B$), and Sobocinski or quasi conjunction ($\wedge_S$). After recalling some logical and probabilistic results in the trivalent logics, we study  basic properties for the  notions of iterated conditioning introduced by Calabrese ($|_C$), by  de Finetti ($|_{dF}$), and by Farrell $(|_F)$, respectively. For each of them we also compute some  sets of coherent assessments on families of conditional events  involving  iterated conditionals. 
This study is based on a geometrical  approach for coherence checking of conditional probability assessments, which allows zero probability for conditioning events (\cite{gilio90}).
Then, we observe that none of this operations of iterated conditioning preserves the compound probability theorem. By exploiting the structure $
\Box|\bigcirc = \Box\wedge \bigcirc +\prev(\Box|\bigcirc)\no{\bigcirc}$, for each 
conjunction among the four selected
trivalent logics we introduce, in the framework of conditional random quantity, a suitable notion of iterated conditioning ($|_K$,$|_L$, $|_B$ and $|_S$).
We observe that all of them  satisfy the compound prevision theorem, we check the validity of some other basic properties, the Import-Export principle and two generalized versions of the    Bayes' Rule for iterated conditioning. 
Finally, we remark that, among the selected iterated conditionals,  $|_{gs}$ is the only one which 
 satisfies all the basic properties.\\

This paper is a revised and expanded version of the conference paper
\cite{CaSa22Sum}. We reorganized the structure of the conference paper and we added proofs and new results.
In particular, we expanded Section \ref{SEC:2} by also adding two new examples.
We included  the proofs 
in Section \ref{SEC:Calabrese}  and  in Section  \ref{SEC:deFinetti}.
We added several new results in Section \ref{SEC:GENITER} (theorems \ref{THM:LUK}--\ref{THM:LUS}).
We also added  Section \ref{SEC:Farrell} and 
Section \ref{SEC:MP}.

The paper is organized as follows. In Section \ref{SEC:2} 
we recall some basic notions and results which concern coherence, conditional random quantities, trivalent logics   and logical operations among conditional events.  We also give two examples on the geometrical interpretation of  coherence. We end the section by  recalling some 
logical and probabilistic properties satisfied by events and conditional events and we generalize them to compound and iterated conditionals.
Then, we check the validity of those selected properties
for the iterated conditioning defined by Cooper-Calabrese ($|_C$, Section \ref{SEC:Calabrese}), by de Finetti 
($|_{dF}$, Section \ref{SEC:deFinetti}) and by Farrell ($|_F$, Section \ref{SEC:Farrell}). In order to check some probabilistic properties, we also compute sets of coherent assessments on suitable families of conditional events.
In Section \ref{SEC:GENITER}, we recall some results on $|_{gs}$ and we introduce and study
the iterated conditioning  $|_K,|_L,|_B$, and $|_S$  in the framework of conditional random quantities under coherence.  We  also  consider the validity of two generalized versions of Bayes' Rule. 
In Section \ref{SEC:MP}, for selected operations of iterated conditioning  we check the p-validity of two generalized inference rules: Modus Ponens and two-premise centering. We also consider two examples of natural language.
Finally, in Section \ref{SEC:CONCLUSIONS}  we illustrate a brief summary on  the validity of the basic properties,   we  give  some conclusions and an outlook to future work.
 To improve the readability of the paper we put some proofs in \ref{SEC:APP}.
\section{Preliminary notions, results, and basic properties}
\label{SEC:2}
In this section we first recall some preliminary notions and results on coherence, conditional events and conditional random quantities.  We also  give some examples  and we   deepen some aspects of conditional random quantities. Then, we recall the logical operations among conditional events in  selected trivalent logics 
 and in the framework of conditional random quantities. Finally, we recall 
 some basic logical and probabilistic properties satisfied by events and conditional events and  we rewrite them for compound and iterated conditionals.
\subsection{Events, conditional events, conditional random quantities and coherence.}
\label{SEC:2.1}
An event $A$ is
a two-valued logical entity which is either  \emph{true} (T), or \emph{false} (F).
We use
 the same symbol to refer to an  event and its indicator, which can take value 1, or 0, according to whether, the event is true, or false, respectively.
We  denote by
$\Omega$ the sure event and by $\emptyset$ the impossible one.
We denote by $A\land B$ (resp., $A\vee B$), or simply by $AB$, the  conjunction (resp., disjunction) of $A$ and $B$. By $\no{A}$ we denote the negation of $A$. 
We simply write $A \subseteq B$ to denote
that $A$ logically implies $B$, i.e., $AB=A$.
Given two events $A$ and $H$, with $H \neq \emptyset$, the conditional event $A|H$  is a three-valued logical entity which is \emph{true}, or
	\emph{false}, or \emph{void} (V), according to whether $AH$ is true, or
	$\no{A}H$ is true, or $\no{H}$ is true, respectively.
Notice that conditional events are three-valued logical entities and hence  they are not in general events. However, as already observed, a conditional event $A|H$ reduces to the (unconditional) event $A$, when $H$ is the sure event $\Omega$, i.e. $A|\Omega=A$.
We recall that, given any conditional event $A|H$, it holds that  $AH|H=A|H$. Moreover,  the negation $\no{A|H}$ is defined as $\no{A|H}=\no{A}|H$ 
Given two conditional events  $A|H$ and $B|K$, we say that $A|H$ logically implies $B|K$, denoted by $A|H \subseteq B|K$, if and only if $AH$ logically implies $BK$  and $\no{B}K$ logically implies $\no{A}H$, that is  (\cite{GoNg88}), 
\begin{equation}\label{EQ:GN}
	A|H \subseteq B| K\;\; \Longleftrightarrow \;\;
	AH\subseteq BK \text{ and }  \no{B}K\subseteq \no{A}H.
\end{equation}		
In the betting framework of subjective probability, to assess $P(A|H)=x$ amounts to say that, for every real number $s$,  you are willing to pay 
an amount $s\,x$ and to receive $s$, or 0, or $s\, x$, according
to whether $AH$ is true, or $\no{A}H$ is true, or $\no{H}$
is true (the bet is called off), respectively. Hence, for the random gain $G=sH(A-x)$, the possible values are $s(1-x)$, or $-s\,x$, or $0$, according
to whether $AH$ is true, or $\no{A}H$ is true, or $\no{H}$
is true, respectively. 

We denote by $X$ a \emph{random quantity}, that is  an 
uncertain real quantity,  which has a well-determined but unknown value. 
In this paper we assume that  $X$ has a finite set of possible values. Given any event $H\neq \emptyset$, 
agreeing to the betting metaphor, if you  assess that the prevision of $``X$ {\em conditional on} $H$'' (or short:  $``X$ {\em given} $H$''), $\pr(X|H)$, is equal to $\mu$, this means that for any given  real number $s$ you are willing to pay an amount $s\mu$ and to receive  $sX$, or $s\mu$, according  to whether $H$ is true, or  false (bet  called off), respectively.
The random gain is 
\begin{equation}\label{EQ:RG}
G=s(XH+\mu \widebar{H})-s\mu=
sH(X-\mu).
\end{equation} In particular, when $X$ is (the indicator of) an event $A$, then 
\begin{equation}\label{EQ:PREVANDP}
\prev(X|H)=P(A|H).
\end{equation}

Given a prevision function $\prev$ defined on an arbitrary family $\mathcal{K}$ of finite
conditional random quantities, consider a finite subfamily $\F = \{X_1|H_1, \ldots,X_n|H_n\} \subseteq \mathcal{K}$ and the vector
$\M=(\mu_1,\ldots, \mu_n)$, where $\mu_i = \prev(X_i|H_i)$ is the
assessed prevision for the conditional random quantity $X_i|H_i$, $i\in \{1,\ldots,n\}$.
With the pair $(\F,\M)$ we associate the random gain $G =
\sum_{i=1}^ns_iH_i(X_i - \mu_i)$. We  denote by $\G_{\mathcal{H}_n}$ the set of possible values of $G$ restricted to $\H_n= H_1 \vee \cdots \vee H_n$. 
Then,  the notion of coherence is defined as below.
\begin{definition}\label{COER-RQ}{\rm
		The function $\prev$ defined on $\mathcal{K}$ is coherent if
		and only if $\forall n
		\geq 1$, $\forall \, s_1, \ldots,
		s_n$, $\forall \, \F=\{X_1|H_1, \ldots,X_n|H_n\} \subseteq \mathcal{K}$,    it holds that: $min \; \G_{\mathcal{H}_n} \; \leq 0 \leq max \;
		\G_{\mathcal{H}_n}$. }
\end{definition}
In other words, $\prev$ on $\mathcal{K}$ is incoherent, if
		and only if there exists a finite combination of $n$ bets such that, 
after discarding the case where  all the bets are called off, the values of the random gain are   all positive or all negative. In the particular case where $\mathcal{K}$ is a family of conditional events, by recalling (\ref{EQ:PREVANDP}),  then Definition \ref{COER-RQ} becomes  the well-known definition of coherence for a conditional probability function, denoted  as $P$. In this case, for a finite subfamily of conditional events $\F=\{E_1|H_1,\ldots,E_n|H_n\}$, we  denote by $\P=(p_1,\ldots,p_n)$, with $p_i=P(E_i|H_i)$, $i=1,\ldots,n$, the restriction  of $P$ to $\F$.  
We observe that, for the checking of coherence of the probability  assessment $p$ on a conditional event $A|H$, only the cases in which the bet is not called off are considered. Then,  we do not consider objects  $A|H$ with $H=\emptyset$, that is conditional events with an impossible conditioning event, because in this case   a bet on $A|\emptyset$ can only be  called off.
Given a conditional event $A|H$, with $H\neq \emptyset$, 
if $\emptyset\neq AH\neq H$, then any value $P(A|H)\in[0,1]$ is a coherent assessment for $A|H$. Coherence requires that $P(AH|H)=P(H|H)=1$ (resp., $P(A|H)=P(\emptyset|H)=0$) when  $AH=H$  (resp., $AH=\emptyset$).

\subsection{Geometrical interpretation of coherence.}
\label{SEC:2.2}
Given a family $\F = \{X_1|H_1,\ldots,X_n|H_n\}$ of $n$ conditional random quantities, for each $i \in \{1,\ldots,n\}$ we denote by $\{x_{i1}, \ldots,x_{ir_i}\}$ the set of possible values  of $X_i$ when  $H_i$ is true; then,  we set $A_{ij} = (X_i = x_{ij})$,   $i=1,\ldots,n$, $j = 1, \ldots, r_i$. 
We observe that, 
 for each $i$,  the family 
 $\{A_{i1}H_i,\ldots,A_{ir_{i}}H_i,\no{H}_i\}$, or equivalently $\{(X_{i}=x_{i1}),\ldots,(X_{i}=x_{ir_i}), \no{H}_{i}\}$, is a partition of the sure event $\Omega$, with  $A_{ij}H_i=A_{ij}$ and         $\bigvee_{j=1}^{r_i}A_{ij}=H_i$. 
Then,
\begin{equation}\label{EQ:OMEGARQ}
\Omega=(A_{11}\vee \cdots \vee A_{1r_1} \vee \no{H}_{1})\wedge \cdots\wedge (A_{n1}\vee \cdots \vee A_{nr_n} \vee \no{H}_{n}).
\end{equation}
By applying the distributive property, the
expression in (\ref{EQ:OMEGARQ}) becomes a  disjunction of ${(r_1+1)\cdots (r_n+1)}$ conjunctions.  
The non-impossible conjunctions 
are the constituents, or possible elementary outcomes,  generated by the family $\F$. 
We denote by $C_0$  the constituent $ \widebar{H}_1 \cdots \widebar{H}_n$ (if non-impossible) and  we denote by $C_1, \ldots, C_m$ the (remaining) constituents
which logically imply $\H_n=H_1\vee \cdots \vee H_n$, that is $C_h\subseteq  \H_n$, $h=1,\ldots,m$. Of course, $C_0,C_1,\ldots, C_m$
 form a partition of the sure event $\Omega$.

With each $C_h,\, h \in \{1,\ldots,m\}$, we associate a vector
$Q_h=(q_{h1},\ldots,q_{hn})$, where $q_{hi}=x_{ij}$ if $C_h \subseteq
A_{ij},\, j=1,\ldots,r_i$, while $q_{hi}=\mu_i$ if $C_h \subseteq \widebar{H}_i$;
with $C_0$ we associate  $Q_0=\M = (\mu_1,\ldots,\mu_n)$.
Denoting by $\I$ the convex hull of $Q_1, \ldots, Q_m$, the condition  $\M\in \I$ amounts to the existence of a vector $(\lambda_1,\ldots,\lambda_m)$ such that:
$ \sum_{h=1}^m \lambda_h Q_h = \M \,,\; \sum_{h=1}^m \lambda_h
= 1 \,,\; \lambda_h \geq 0 \,,\; \forall \, h$; in other words, $\M\in \I$ is equivalent to the solvability of the system $(\Sigma)$, associated with  $(\F,\M)$,
\begin{equation}\label{SYST-SIGMA}
\begin{small}
(\Sigma) \quad
\begin{array}{ll}
\sum_{h=1}^m \lambda_h q_{hi} =
\mu_i \,,\; i \in\{1,\ldots,n\} \,, 
\sum_{h=1}^m \lambda_h = 1,\;\;\lambda_h \geq 0 \,,\;  \,h \in\{1,\ldots,m\}\,.
\end{array}
\end{small}
\end{equation}
Given the assessment $\M =(\mu_1,\ldots,\mu_n)$ on  $\F =
\{X_1|H_1,\ldots,X_n|H_n\}$, let $S$ be the set of solutions $\Lambda = (\lambda_1, \ldots,\lambda_m)$ of system $(\Sigma)$.   
We point out that the solvability of  system $(\Sigma)$  is a necessary (but not sufficient) condition for coherence of $\M$ on $\F$. When $(\Sigma)$ is solvable, that is  $S \neq \emptyset$, we define:
\begin{equation}\label{EQ:I0}
\begin{small}
\begin{array}{l}
\Phi_{i}(\Lambda ) = \Phi_{i}(\lambda_{1}, \ldots , \lambda_{m}) =
\sum_{r :
	C_{r} \subseteq H_{i}} \lambda_{r}, \, i\in \{1,\ldots,n\}, \; \Lambda \in S;\\
M_{i} = \max_{\Lambda \in S },\Phi_{i}(\Lambda ),\;
i\in \{1,\ldots,n\}\,;\\
I_{0} = \{ i \, : \, M_{i}=0 \} \,,\F_0 =
\bigcup_{i\in I_0} \{X_i|H_i\}.  
\end{array}
\end{small}
\end{equation}
We also denote by $\M_0$ the sub-assessment of $\M$ on the sub-family $\F_0$.
For what concerns the probabilistic meaning of $I_0$, it holds that  $i\in I_0$ if and only if the (unique) coherent extension of $\M$ to $H_i|\H_n$ is zero.
Then, the following theorem can be proved  (
 see e.g., \cite{GiSa21IJAR}
)
\begin{theorem}
	\label{CNES-PREV-I_0-INT}{\rm 
		A conditional prevision assessment ${\M} = (\mu_1,\ldots,\mu_n)$ on
		the family $\F = \{X_1|H_1,\ldots,X_n|H_n\}$ is coherent if
		and only if  the following conditions are satisfied: 
		(i) the system $(\Sigma)$ defined in (\ref{SYST-SIGMA}) is solvable; (ii) if $I_0 \neq \emptyset$, then $\M_0$ is coherent. }
\end{theorem}
Let  $\mathcal{S}'$ be a nonempty subset of the set of solutions
$\mathcal{S}$ of system $(\Sigma) $. We denote by  $I_0'$ the set $I_0$ defined as in (\ref{EQ:I0}), where $\mathcal{S}$ is replaced by $\mathcal{S}'$, that is
\begin{equation}\label{EQ:I0'}
I_{0}' = \{ i:  \, M'_{i}=0 \}, \mbox{ where } M'_{i}=\max_{\Lambda \in \mathcal{S}' } \; \Phi_{i}(\Lambda ) \; ,
\; \; \;
i\in \{1,\ldots,n\}.
\end{equation}
 Moreover, we denote by $(\F_0',\M_0')$ the pair associated with $I_0'$. Then, we obtain 
\begin{theorem}
	\label{CNES-PREV-I_0-INT'}The assessment $\mathcal{M}$ on $\mathcal{F}$ is coherent if
		and only if the following conditions are satisfied: (i) 	the system $(\Sigma) $ associated with the pair $(\mathcal{F},
	\mathcal{M})$ is solvable; (ii) if $I_{0}' \neq \emptyset $, then
	$\mathcal{M}_{0}'$ is coherent.
\end{theorem}

Of course, the previous results  can  be used in the case of a probability assessment, which will be denoted by  $\P$, on a family of  conditional events $\F$.   
More precisely, given a family  
 $\mathcal{F} = \{E_{1}|H
_{1}, \ldots, E_{n}|H_{n}\}$ of $n$ conditional events, we observe that, for each $i$, the family $\{E_iH_i, \no{E}_iH_i, \no{H}_i\}$ is a partition of $\Omega$.
Then, 
$\Omega =\bigwedge_{i=1}^n(E_iH_i \vee \widebar{E}_iH_i \vee \widebar{H}_i)$. 
By  applying the  distributive property  it follows that 
$\Omega$
can also  be written as the disjunction of $3^n$ logical
conjunctions, some of which may be impossible. 
The remaining ones, denoted by $C_0,C_1,\ldots,C_m$,  where $C_0=\no{H}_1 \cdots \no{H}_n$ (if non-impossible),
are the constituents generated by  $\mathcal{F}$. 
Of course, $m+1\leq 3^n$, with  ${m+1=3^n}$ when, for example, the events in $\F$ are logically independent. 
Let $\P=(p_1,\ldots,p_n)$ be a probability assessment on $\F$.  In this case, for  each constituent $C_h$, $h=1,\ldots,m$, 
it holds that   
$Q_h = (q_{h1}, \ldots, q_{hn})$, where $q_{hi} = 1$, or 0, or $p_i$, according to whether $C_h \subseteq E_iH_i$, or $C_h \subseteq \no{E}_iH_i$, or $C_h \subseteq \no{H}_i$. 
The point $Q_0=\P$
is associated with $C_0$. 
In the following example, which will be useful in Section \ref{SEC:Calabrese} to prove Theorem \ref{ThCalabrese2}, we illustrate   how the constituents and the associated points are generated in order to check the coherence of a probability assessment.

\begin{example}\label{EX:EX1}
Let $\F=\{E_1|H_1,E_2|H_2\}=\{B|K, B|(K\wedge (\no{H}\vee A))\}$, where 
 $A$, $B$, $H$, $K$ are  four logically independent events, and 
let $\P=(p_1,p_2)$ be a probability assessment on  $\F$. We verify that  $\P$ is coherent for every $(p_1,p_2)\in[0,1]^2$. 
It holds that
\[
\begin{array}{ll}
\Omega=(E_1H_1\vee \no{E}_1H_1\vee \no{H}_1)\wedge (E_2H_2\vee \no{E}_2H_2\vee \no{H}_2)=\\
=
[BK\vee\no{B}K\vee \no{K}]\wedge [(BK\wedge(\no{H}\vee A)) \vee (\no{B}K\wedge(\no{H}\vee A)) \vee (\no{K}\vee \no{A}H) ]=
\\
=(BK\wedge(\no{H}\vee A))\vee 
\emptyset \vee \no{A}HBK \vee\emptyset  \vee  (\no{B}K\wedge(\no{H}\vee A))\vee \no{A}H\no{B}K\vee \emptyset \vee \emptyset \vee (\no{K} \vee\no{A}H\no{K}) =\\
=C_1\vee C_2\vee C_3 \vee C_4 \vee C_0,
\end{array}
\]
where the constituents are
\[
\begin{array}{l}
C_1=BK\wedge(\no{H}\vee A)=AHBK \vee \no{H}BK,\ 
C_2=\no{A}HBK, 
C_3=\no{B}K\wedge(\no{H}\vee A)=AH\no{B}K \vee \no{H}\no{B}K,\ 
\\ C_4=\no{A}H\no{B}K,\ C_0=\no{K} \vee\no{A}H\no{K}=\no{K}.
\end{array}
\]
The points $Q_h$'s associated with the pair $(\F,\P)$ are
\[
\begin{array}{l}
Q_1=(1,1),\ Q_2=(1,p_2),\ Q_3=(0,0),\ Q_4=(0,p_2),\ Q_0=\P=(p_1,p_2).
\end{array}
\] 
We denote by $\I$ the convex hull of points $Q_1,Q_2,Q_3,Q_4$ (see Figure \ref{FIG:EX1}).)
\begin{figure}[tbph]
	\centering
	\includegraphics[width=0.5\linewidth]{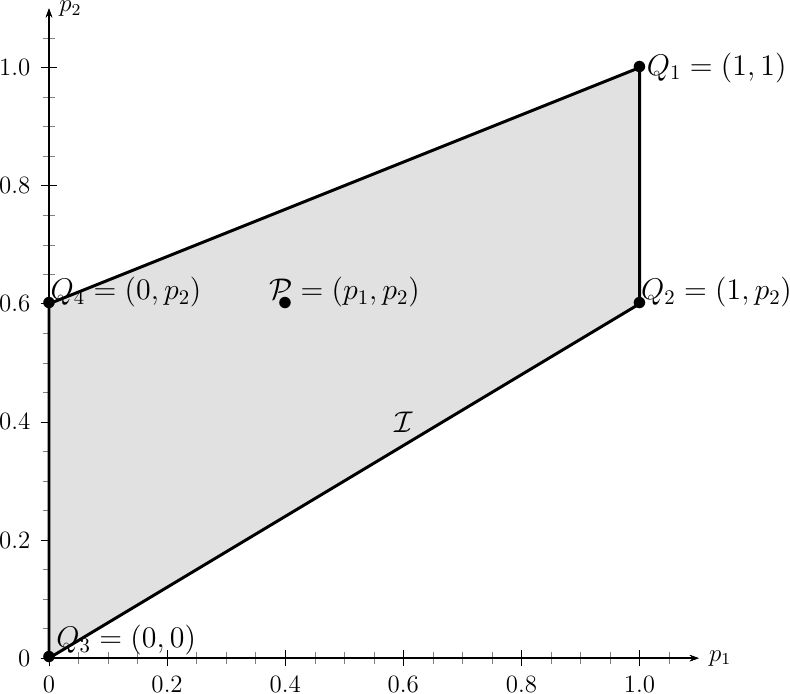}
	\caption{Convex hull of  the points $Q_1, Q_2,Q_3,Q_4$ associated with the pair $(\F,\P)$, where  $\F=\{B|K, B|(K\wedge (\no{H}\vee A))\}$ and $\P=(p_1,p_2)$.  
		In the figure the numerical  values are: $p_1=0.4$ and  $p_2=0.6$.}
	\label{FIG:EX1}
\end{figure}
The system $(\Sigma)$ in (\ref{SYST-SIGMA}) associated with the pair $(\F,\P)$  is\\
\begin{equation}
\left\{
    \begin{array}{llll}
     \lambda_1+\lambda_2=p_1, \\
    \lambda_1+p_2\lambda_2+p_2\lambda_4=p_2, \\
    \lambda_1+\lambda_2+\lambda_3+\lambda_4=1,\\
    \lambda_i\geq0, \ \forall i=1, \dots, 4.\\
    \end{array}
    \label{sisEs}
\right.    
\end{equation}
We observe that $\P=(p_1,p_2)=p_1Q_2+(1-p_1)Q_4$ and hence the system (\ref{sisEs}) is solvable and a solution is $\Lambda=(0,p_1,0,1-p_1)$, with $p_1\in[0,1]$. \\
By considering the  function  $\phi$ as defined in (\ref{EQ:I0}),  it holds that 
\[
\phi_{1}(\Lambda)=\sum_{h:C_h\subseteq H}\lambda_h=\lambda_1+\lambda_2+\lambda_3+\lambda_4=1,
\]
\[
\phi_{2}(\Lambda)=\sum_{h:C_h\subseteq (K \wedge(A\vee \no{H}))}\lambda_h=\lambda_1+\lambda_3=0.
\]
We set $\mathcal{S'}=\{\Lambda\}$ and we get
 $\I_0'=\{2\}$. We observe that  the sub-assessment $p_2$ on $B|(K\wedge (\no{H}\vee A))$ is coherent for every $p_2\in[0,1]$. Thus, by Theorem \ref{CNES-PREV-I_0-INT'}, the assessment $\P=(p_1,p_2)$  on $\F$ is coherent $\forall (p_1,p_2)\in[0,1]^2$.
\end{example}
The next example, which  is related to the compound prevision theorem listed in Section \ref{SEC:3}, illustrates 
 how to check coherence for (conditional) prevision assessments.

\begin{example}\footnote{This example was inspired by  Angelo Gilio's talk ``On Coherence and Conditionals'' presented  at the workshop   ``Reasoning and uncertainty: probabilistic, logical, and psychological perspectives'', Regensburg, August 9-10, 2022.}[Compound prevision theorem]\label{EX:CPT}
Let $\F=\{XH,X|H,H\}$, where $H$ is an event, with
 $H\neq \emptyset$, $H\neq  \Omega$, and $X$ is a  finite random quantity. We denote by  $\{x_1,\ldots,x_n\}$ the possible values of $X$ when $H$ is true. We also  set $x'=\min\{x_1,\ldots,x_n\}$  and $x''=\max\{x_1,\ldots,x_n\}$.
 Let $\M=(z,\mu, p)$ be a prevision assessment on  $\F$. We verify that  $\M$ is coherent if and only if  $\mu\in [x',x'']$, $p\in[0,1]$, and $z=\mu p$.  Of course,  coherence requires that $p\in[0,1]$.
It holds that
\[
\begin{array}{ll}
\Omega= C_1\vee \cdots \vee C_n\vee  C_{n+1},
\end{array}
\]
where the constituents are
\[
\begin{array}{l}
C_1=(X=x_1),\cdots, C_n=(X=x_{n}),\; C_{n+1}=\no{H}.
\end{array}
\]
The points $Q_h$'s associated with the pair $(\F,\M)$ are
\[
\begin{array}{l}
Q_1=(q_{11},q_{12},q_{13})=(x_1,x_1,1),\ldots,\ Q_n=(q_{n1},q_{n2},q_{n3})=(x_n,x_n,1),\\ Q_{n+1}=(q_{(n+1)1},q_{(n+1)2},q_{(n+1)3})=(0,\mu,0).
\end{array}
\] 
We denote by $\I$ the convex hull of points $Q_1,\ldots,Q_{n+1}$.
The system $(\Sigma)$ in (\ref{SYST-SIGMA}) associated with the pair $(\F,\M)$  is\\
\begin{equation}
\label{EQ:SIGMACPT}
\left\{
    \begin{array}{llll}
     \lambda_1x_1+\cdots+\lambda_nx_n=z, \\
     \lambda_1x_1+\cdots+\lambda_nx_n+\lambda_{n+1}\mu=\mu,
     \\
\lambda_1+\cdots+\lambda_n=p, \\
\lambda_1+\cdots+\lambda_{n+1}=1, \\

    \lambda_i\geq0, \ i=1, \dots, n+1.\\
    \end{array}
\right.    
\Longleftrightarrow \left\{
    \begin{array}{llll}
     \lambda_1x_1+\cdots+\lambda_nx_n=z, \\
     z =\mu p,
     \\
\lambda_1+\cdots+\lambda_n=p, \\
\lambda_{n+1}=1-p, \\

    \lambda_i\geq0, \  i=1, \dots, n+1.\\
    \end{array}
\right.
\end{equation}
We distinguish two cases: $(i)$ $p=0$; and $p\in (0,1]$.
\\
Case $(i)$. System (\ref{EQ:SIGMACPT}) is solvable only if $z=0$, with the unique solution given by $\Lambda=(\lambda_1,\ldots,\lambda_{n},\lambda_{n+1})=(0,\ldots,0,1)$. 
By considering the  function  $\phi$ as defined in (\ref{EQ:I0}),  it holds that 
\[
\phi_{1}(\Lambda)=\phi_{3}(\Lambda)=\sum_{h:C_h\subseteq \Omega}\lambda_h=\lambda_1+\cdots+\lambda_{n+1}=1>0
\]
and 
\[
\phi_{2}(\Lambda)=\sum_{h:C_h\subseteq H}\lambda_h=\lambda_1+\cdots+\lambda_{n}=p=0.
\]
Then $I_0=\{2\}$. We consider the  sub-assessment $\M_0=(\mu)$ on $\F_0=\{X|H\}$.
It can be easily proved that  $\mu=\prev(X|H)$ is coherent if and only if $\mu\in[x',x'']$. Then, by Theorem \ref{CNES-PREV-I_0-INT}, the assessment $\M=(0,\mu,0)$ on $\F$ is coherent if and only if $\mu\in[x',x'']$.\\
Case $(ii)$.  In this case System (\ref{EQ:SIGMACPT})   becomes
\begin{equation}
	\label{EQ:SIGMACPTBIS}
	\left\{
	\begin{array}{llll}
		z=\mu p, \\
		\mu=\frac{z}{p} =\frac{\lambda_1x_1+\cdots+\lambda_nx_n}{\lambda_1+\cdots+\lambda_n},
		\\
		\lambda_1+\cdots+\lambda_n=p, \\
		\lambda_{n+1}=1-p, \\
		
		\lambda_i\geq0, \  i=1, \dots, n+1,\\
	\end{array}
	\right.
\end{equation}
which is solvable when $z=\mu p$  and  $\mu\in[x',x'']$, because $\mu$  is a convex linear combination of $\{x_1,\ldots,x_n\}$ with weights $\frac{\lambda_1}{\lambda_1+\cdots+\lambda_n},\cdots, \frac{\lambda_n}{\lambda_1+\cdots+\lambda_n}$ . 
For each solution $\Lambda$ of System (\ref{EQ:SIGMACPTBIS})
it holds that $\phi_{1}(\Lambda)=\phi_{3}(\Lambda)=1>0$ and $\phi_{2}(\Lambda)=p>0$. Then, as
 $I_0=\emptyset$, by Theorem \ref{CNES-PREV-I_0-INT} the assessment $(z,\mu,p)$, with $p\in(0,1]$ is coherent if and only $\mu\in[x',x'']$ and $z=\mu p$. \\ Therefore
$\M$ is coherent if and only if  $\mu\in [x',x'']$, $p\in[0,1]$, and $z=\mu p$.
  We now show that, in agreement to Definition \ref{COER-RQ}, an incoherent assessment leads to   the  existence of a combination of bets where the values of the random gain are all positive or all negative (Dutch Book). 
We set $\M=(z,\mu,p)$, with $p\in[0,1]$, $\mu\in[x',x'']$ and $z\neq \mu p$.
We observe that the points $Q_1,\ldots,  Q_{n+1}$ belong to the plane $\pi: -\mathcal{X}+\mathcal{Y}+\mu \mathcal{Z}=\mu$, where $\mathcal{X}, \mathcal{Y}, \mathcal{Z}$ are the axes coordinates. Moreover, $\M$ does not belong to the convex hull $\I$ because, as $z\neq\mu p$,  System (\ref{EQ:SIGMACPT}) is not solvable.
	Now, by considering the function $f(\mathcal{X},\mathcal{Y},\mathcal{Z})= -\mathcal{X}+\mathcal{Y}+\mu \mathcal{Z}-\mu$, we observe that for each constant $k$ the equation $f(\mathcal{X},\mathcal{Y},\mathcal{Z})=k$ represents a plane which is
	parallel to $\pi$ and coincides with $\pi$ when $k=0$. Then,  $f(Q_1)=\cdots=f(Q_{n+1})=0$.
	Moreover,  $f(Q_h)-f(\M)=z-\mu p\neq 0$, $h=1,\ldots,n+1$. 
We recall that for $h=1,\ldots, n+1$,  the value $g_h$, of the random gain 
\[
G=s_1(XH-z)+s_2H(X-\mu)+s_3(H-p) 
\]
associated to the constituent $C_h$  is
\[
g_h=s_1(q_{h1}-z)+s_2H(q_{h2}-\mu)+s_3(q_{h3}-p).
\]
We observe that, by setting the stakes $s_1=-1$, $s_2=1$, $s_3=-\mu $, it holds that  $g_h=f(Q_h)-f(\M)=z-\mu p\neq 0$, $h=1,\dots,n+1$. Therefore,  $\min g_h  \cdot \max g_h>0$, when $s_1=-1$, $s_2=1$, $s_3=-\mu $, that is  a combination of bets where the value of the random gain are all positive or all negative.

\end{example}
\begin{remark}\label{REM:CPT}
By Example \ref{EX:CPT} coherence requires that (\emph{compound prevision theorem})
\begin{equation}\label{EQ:CPREVT}
\prev(XH)=\prev(X|H)P(H).
\end{equation}
In particular, when $X$ is (the indicator of) an event $E$, 
Equation (\ref{EQ:CPREVT}) becomes (\emph{compound probability theorem})
\begin{equation}\label{EQ:CPT}
    P(E\wedge H)=P(E|H)P(H).
\end{equation}   
Moreover, under logical independence of $E$ and $H$,  Example \ref{EX:CPT} shows that 
$P(E|H)$  coincides with  the ratio $\frac{P(E\wedge  H)}{P(H)}$  when $P(H)>0$, and $P(E|H)$ can be any value in $[0,1]$ when $P(H)=0$.  
Then, 
differently from the ``standard'' approach where $P(E|H)$ is defined only when when $P(H)>0$ by the ratio $\frac{P(E\wedge H)}{P(H)}$, 
in the  coherence-based approach the conditional probability $P(E|H)$ is a primitive notion which is properly  defined even if ${P(H)=0}$. 
\end{remark}
\subsection{Numerical interpretation of a conditional random quantity}
Given a conditional event $A|H$, 
with  $P(A|H) = x$,
the indicator of $A|H$,  denoted by the same symbol, is the following random quantity (see, e.g., \cite{gilio90,Lad95})
\begin{equation}\label{EQ:AgH}
A|H=
AH+x \no{H}=AH+x (1-H)=\left\{\begin{array}{ll}
1, &\mbox{if $AH$ is true,}\\
0, &\mbox{if $\no{A}H$ is true,}\\
x, &\mbox{if $\no{H}$ is true.}\\
\end{array}
\right.
\end{equation} 
Notice that, by the linearity property of a coherent  prevision, it holds that
\[\prev(AH + x\no{H})=xP(H)+xP(\no{H})=x=P(A|H).\]
Then, in a conditional bet on $A|H$, the indicator in  (\ref{EQ:AgH}) represents the {\em random win} that you receive when you pay  the amount $P(A|H)=x\in [0,1]$. Indeed, you  receive 1 (you win), or 0 (you lose), or $x$ (the  bet is called off), according to whether $AH$ is true, or $\no{A}H$ is true, or $\no{H}$ is true, respectively.

We also observe that, the   value $x$ of the random quantity  $A|H$  (subjectively) depends on the assessed probability  $P(A|H)=x$. 
 When $H\subseteq A$ (i.e., $AH=H$), it holds that $P(A|H)=1$; then,
for the indicator $A|H$, when $H\subseteq A$, it holds that 
$A|H=AH+x\no{H}=H+\no{H}=1. 
$
Similarly, if $AH=\emptyset$, as  $P(A|H)=0$, it follows that 
$A|H=0+0\no{H}=0$.
For the indicator of the negation  of $A|H$, as $P(\no{A}|H)=1-P(A|H)$, it holds that  $\no{A}|H=1-A|H$.
Given two conditional events $A|H$ and $B|K$, for every coherent assessment $(x,y)$ on $\{A|H,B|K\}$, 
it holds that (\cite[formula (15)]{GiSa21A})
\begin{equation*}
AH+x\no{H}\leq BK+y\no{K} 
\;\Longleftrightarrow 
\; \mbox{either } A|H \subseteq B|K, \mbox{ or } AH=\emptyset, \mbox{ or } K\subseteq B,
\end{equation*}
that is,   between the numerical values of $A|H$ and $B|K$, under coherence it holds that
\begin{equation}
\label{EQ:LEQ}
A|H\leq B|K 
\;\Longleftrightarrow \;
\mbox{either } A|H \subseteq B|K, \mbox{ or } AH=\emptyset, \mbox{ or } K\subseteq B.
\end{equation}
Of course, the relation $A|H \leq  B|K$ requires that $P(A|H)=x\leq y=P(B|K)$. Then, 
\begin{equation}
\label{EQ:LEQP}
P(A|H)\leq P(B|K)\; \forall \text{coherent } P  
\;\Longleftrightarrow \;
\mbox{either } A|H \subseteq B|K, \mbox{ or } AH=\emptyset, \mbox{ or } K\subseteq B,
\end{equation}
which, when $H=K=\Omega$, reduces to 
\[
P(A)\leq P(B)\, \forall \text{coherent } P  \Longleftrightarrow A\subseteq B \Longleftrightarrow A\leq  B.
\]
By following the approach given in \cite{CoSc99a,GiSa14,lad96}, once  a coherent assessment $\mu=\prev(X|H)$ is specified,  the conditional random quantity $X|H$ (is not looked at as the restriction to $H$, but)
is defined as  $X$, or $\mu$,
according
to whether $H$ is true, or $\widebar{H}$
is true; that is,  
\begin{equation}\label{EQ:XgH}
X|H=XH+\mu \widebar{H}.
\end{equation}
As shown in (\ref{EQ:XgH}), given any random quantity $X$ and any event $H\neq \emptyset$, in the framework of subjective probability, in order to define $X|H$ we just need to specify the value $\mu$ of the conditional prevision $\prev(X|H)$.  Indeed, once  the value $\mu$ is specified, the object $X|H$ is (subjectively) determined.
We observe that (\ref{EQ:XgH}) is consistent because
\begin{equation*}\label{EQ:PREVXgH}
\prev(XH+\mu\no{H})=\prev(XH)+\mu P(\no{H})=\prev(X|H)P(H)+\mu P(\no{H})=\mu P(H)+\mu P(\no{H})=\mu.
\end{equation*}
By (\ref{EQ:XgH}), the random gain associated with a bet on $X|H$ can be represented as 
$G=s(X|H-\mu)$,
that is $G$ is the difference between what you 
receive, $s X|H$, and what you pay, $s \mu$.
In what follows,
for any given conditional random quantity $X|H$, we assume that, when $H$ is true, the set of possible values of $X$ is finite.
 In this case we say 
that $X|H$ is a finite conditional random quantity.  
Denoting by $\{x_{1}, \ldots,x_{r}\}$ the set of possible values  of $X$ restricted to $H$ and by setting $A_{j} = (X = x_{j})$, $j=1,\ldots, r$, it holds that 
$\bigvee_{j=1}^r A_j=H$ and
$
X|H=XH+\mu \widebar{H}=x_1A_1+\cdots +x_rA_r+\mu\widebar{H}.
$

  The  result below (\cite[Theorem 4]{GiSa14}) shows that if  
  two conditional random quantities $X|H$, $Y|K$ coincide when  $H\vee K$ is true, then $X|H$ and $ Y|K$ also coincide when  $H\vee K$ is false, and hence $X|H$ coincides with $Y|K$ in all cases.
 \begin{theorem}\label{THM:EQ-CRQ}{\rm Given any events $H\neq \emptyset$ and  $K\neq \emptyset$, and any random quantities $X$ and $Y$, let $\Pi$ be the set of the coherent prevision assessments $\prev(X|H)=\mu$ and $\prev(Y|K)=\nu$. \\
 		$(i)$ Assume that, for every $(\mu,\nu)\in \Pi$,  $X|H=Y|K$  when  $H\vee K$ is true; then   $\mu=\nu$ for every $(\mu,\nu)\in \Pi$. \\
		$(ii)$ For every $(\mu,\nu)\
 		\in \Pi$,   $X|H=Y|K$  when  $H\vee K$ is true  if and only if $X|H=Y|K$.
 }\end{theorem}	
 \begin{remark}\label{REM:INEQ-CRQ}
 Theorem \ref{THM:EQ-CRQ} has been generalized in \cite[Theorem 6]{GiSa19} by replacing the symbol ``$=$'' by ``$\leq$'' in statements  $(i)$ and $(ii)$. In other words,
 if  $X|H\leq Y|K$ when  $H\vee K$ is true, then $\prev(X|H)\leq  \prev(Y|K)$ and hence  $X|H\leq  Y|K$ in all cases.	
 \end{remark}

\subsection{Trivalent logics, logical operations of conditionals and conditional random quantities}
\label{SUBSEC:CONJUNCTIONS}
We recall some notions of conjunction among conditional events in some trivalent logics: Kleene-Lukasiewicz-Heyting conjunction ($\wedge_K$), or de Finetti conjunction (\cite{defi36});  Lukasiewicz conjunction  ($\wedge_L$);  Bochvar internal conjunction, or Kleene weak conjunction ($\wedge_B$);  Sobocinski conjunction, or quasi conjunction ($\wedge_S$).  In all these definitions the result of the conjunction   is still a conditional event with set of truth values $\{\text{true}, \text{false},\text{void}\}$ (see, e.g., \cite{CiDu12,CiDu13}). 
We also recall the notions of conjunction among conditional events,  $\wedge_{gs}$,  introduced as a suitable conditional random quantity in a betting-scheme context(\cite{GiSa14,GiSa19}, see also \cite{Kauf09,McGe89}).
We list below in an explicit way the five conjunctions  and the associated disjunctions obtained by De Morgan's law (\cite{GiSa22}):
\begin{enumerate}
	\item   $(A|H) \wedge_{K} (B|K)=AHBK|(HK\vee \overline{A}H\vee \overline{B}K)$,\\
		$(A|H) \vee_{K} (B|K) = (AH\vee BK)|(\no{A}H\no{B}K\vee AH\vee BK)$\,;
   \item $(A|H)\wedge_L (B|K)=AHBK|( AHBK \vee \overline{A} \, H\vee \, \overline{B} K \vee \overline{H} \, \overline{K})$,\\
   	$(A|H) \vee_{L} (B|K) = (AH \vee BK)|(\no{A}H\no{B}K \vee AH \vee BK\vee   \overline{H} \, \overline{K})$  \,;
	\item   $(A|H)\wedge_{B} (B|K)=AHBK|HK$,\\
		$(A|H)\vee_{B} (B|K) = (A \vee B)|HK$;
	\item 	$(A|H) \wedge_{S} (B|K)=((AH\vee\overline{H})\wedge (BK\vee \overline{K}))|(H\vee K)$,\\
	$(A|H) \vee_{S} (B|K)=(AH\vee BK)|(H\vee K)$;
	\item 	$(A|H) \wedge_{gs} (B|K)=(AHBK+P(A|H)\no{H}BK+P(B|K)AH\no{K})|(H\vee K)$,\\
	$(A|H) \vee_{gs} (B|K)=(AH\vee BK+P(A|H)\no{H}\no{B}K+P(B|K)\no{A}H\no{K})|(H\vee K)$.
\end{enumerate}
 The operations above are  all commutative and associative.
By setting   $P(A|H)=x$, $P(B|K)=y$,   $P[(A|H)\wedge_{i}(B|K)]=z_i$, $i\in\{K,L,B,S\}$, and $\prev[(A|H)\wedge_{gs}(B|K)]=z_{gs}$, based on  (\ref{EQ:AgH}) and on (\ref{EQ:XgH})
 the conjunctions  $(A|H)\wedge_{i}(B|K)$, $i\in\{K,L,B,S,gs\}$ can be also looked at as  random quantities  with set of possible value  illustrated in  Table \ref{TAB:TABLE1L}. A similar interpretation can also be given for  the associated disjunctions. 
 In Table \ref{tab:LUB} we list the lower and upper bounds for the coherent extension $z_i=P[(A|H)\wedge_i (B|K)]$, $i\in\{K,L,B,S\}$ (\cite{SUM2018S}) and $z_{gs}=\prev[(A|H)\wedge_{gs} (B|K)]$ (\cite{GiSa14})  of  the given assessment $P(A|H)=x$ and $P(B|K)=y$ ).
\begin{table}[]
	\def\True{$1$}
	\def\False{$0$}
	\def\Void{$z$}
	\def\Voidx{$x$}
	\def\Voidy{$y$}
	\centering
	\begin{tabular}{l|l|c|c|c|c|c|c}
		&  $A|H$   &  $B|K$   & $\wedge_{K}$ & $\wedge_{L}$ & $\wedge_{B}$ & $\wedge_{S}$ &$\wedge_{gs}$\\
		\hline		                 
	$AHBK  $           & \True  & \True  &   \True    & \True  &   \True    &    \True    & \True\\
			$	{AH}{\overline{B}K}$    & \True  & \False &   \False   & \False &   \False   &   \False    & \False\\
			$	AH\overline{K}  $       & \True  & \Voidy &   \Void$_K$   & \Void$_L$ &   \Void$_B$   &    \True    & \Voidy\\
		$	{\overline{A}H}{BK}$    & \False & \True  &   \False   & \False &   \False   &   \False    &\False    \\
		$	{\overline{A}H\overline{B}K} $& \False & \False &   \False   & \False &   \False   &   \False    &\False    \\
			$	{\overline{A}H} \, \overline{K}$  & \False & \Voidy &   \False   & \False &   \Void$_B$   &   \False    &\False    \\
		$	\overline{H}{BK}$       & \Voidx & \True  &   \Void$_K$   & \Void$_L$ &   \Void$_B$   &    \True  &  \Voidx      \\
		$	\overline{H}\,{\overline{B}K}$  & \Voidx & \False &   \False   & \False &   \Void$_B$   &   \False   & \False  \\
			$	\no{H}\,\no{K} $    & \Voidx & \Voidy &   \Void$_K$    &   \False   & \Void$_B$ &   \Void $_S$ & $z_{gs}$ 
	\end{tabular}
	\caption{Numerical values   (of the indicator)  of the conjunctions $(A|H)\wedge_{i}(B|K)$, $i\in\{K,L,B,S,gs\}$. The triplet  $(x,y,z_i)$ denotes a  coherent assessment on $\{A|H,B|K,(A|H)\wedge_{i}(B|K)\}$.}
	\label{TAB:TABLE1L}
\end{table} 
Notice that, differently from conditional events which are three-valued objects, the conjunction $(A|H) \wedge_{gs} (B|K)$ (and the associated disjunction)
is no longer a three-valued object, but  a five-valued object with values in $[0,1]$.  
In betting terms, the prevision $z_{gs}=\prev[(A|H) \wedge_{gs} (B|K)]$ represents the amount you agree to pay, with the proviso that you will receive the random quantity
$AHBK+x\no{H}BK+yAH\no{K}$, if $H\vee K$ is true, 
  $z_{gs}$ if  $\no{H}\,\no{K}$ is true.
In other words 
by paying $z_{gs}$
you receive:
1, if both conditional events are true; 0, if at least one of the conditional events is false; 
the probability of the conditional event that is void if one conditional event is void and the other one is true;   the amount $z_{gs}$ you paid if both conditional events are void. 
 The notion of conjunction $\wedge_{gs}$ (and disjunction $\vee_{gs}$) among conditional events has been generalized to the case of $n$ conditional events in \cite{GiSa19}.
For some applications see, .e.g., \cite{SGOP20,SPOG18}. Developments of this approach to general compound conditionals has been given in \cite{FGGS22KR}. 
Differently from the other notions of conjunctions, $\wedge_{gs}$ preserves the classical logical and probabilistic properties valid for unconditional events (see, e.g.,\cite{GiSa21A}).
In particular,  the Fréchet-Hoeffding bounds, 
i.e., the lower and upper bounds $z'=\max\{x+y-1,0\},z''=\min\{x,y\}$, 
obtained under logical independence in  the unconditional case for the coherent extensions $z=P(AB)$ of $P(A)=x$ and $P(B)=y$, 
 when  $A$ and $B$ are replaced by $A|H$ and $B|K$, are only satisfied by $z_{gs}$  (see Table~\ref{tab:LUB}).\ \\

\begin{table}[]
    \centering
    \begin{tabular}{c|c|c}
        Conjunction & Lower bound  & Upper bound \\
        \hline
        $(A|H)\wedge_K(B|K)$ & $z_K'=0$ 
        & $z_K''=\min\{x,y\}$
        $\begin{array}{l} ~\\ \\ \end{array}$
        \\
        $(A|H)\wedge_L(B|K)$ & $z_L'=0$ & $z_L''=\min\{x,y\}$ $\begin{array}{l} ~\\ \\ \end{array}$ \\
        $(A|H)\wedge_L(B|K)$ & $z_B'=0$ & $z_B''=1$ $\begin{array}{l} ~\\ \\ \end{array}$\\
        $(A|H)\wedge_S(B|K)$ & $z_S'=\max\{x+y-1,0\}$ & $z_S''=\begin{cases}
             \frac{x+y-2xy}{1-xy}, &  \text{if } (x,y)\neq (1,1)\\
             1,  & \text{if } (x,y)=(1,1)\\
         \end{cases}$

    \end{tabular}
    \caption{Lower and upper bounds for the coherent extension $z_i=P[(A|H)\wedge_i (B|K)]$, $i\in\{K,L,B,S\}$ and $z_{gs}=\prev[(A|H)\wedge_{gs} (B|K)]$  of  the given assessment $P(A|H)=x$ and $P(B|K)=y$.}
    \label{tab:LUB}
\end{table} 

\subsection{Some  basic properties and Import-Export principle}\label{SEC:3}
In  this section, after recalling some basic logical and probabilistic properties satisfied by events and conditional events,  we rewrite them for compound and iterated conditionals, by replacing  events with   conditional events.  We also recall the Import-Export principle and its connection with the Lewis' triviality results.\\ \ \\ 
Given two events $A$ and $B$, with $A\neq \emptyset$, 
it is well-known the validity of following properties 
\begin{enumerate}
\item[1.] $B|A=AB|A$;
\item[2.]  $AB\leq   B|A$ and hence $P(AB)\leq P(B|A)$;
\item[3.] $P(AB)=P(B|A)P(A)$  (\emph{compound probability theorem}, see Remark \ref{REM:CPT});
\item[4.] if $A$ and  $B$ are logically independent, by setting  $P(A)=x$  and $P(B)=y$,  the  extension  $\mu =P(B|A)$ is coherent if and only if  $\mu \in[\mu',\mu'' ]$, where (see, e.g. \cite[Theorem 6]{SPOG18})
\begin{equation}\label{EQ:LU}
\mu'=\left\{
\begin{array}{ll}
\frac{\max\{x+y-1,0\}}{x}, &\text{ if } x\neq 0,\\
0, &\text{ if } x= 0,\\
\end{array}
\right.,\;\;
\mu''=\left\{
\begin{array}{ll}
\frac{\min\{x,y\}}{x}, &\text{ if } x\neq 0,\\
1, &\text{ if } x= 0.\\
\end{array}
\right.
\end{equation}
\end{enumerate}

\noindent By replacing events $A,B$ by conditional events $A|H, B|K$, and for the compound conditionals the symbol of probability $P$ by the symbol of prevision $\prev$, the properties 1, 2, 3, and  4 become:\\
\begin{itemize}
\item[P1.] $(B|K)|(A|H)=[(A|H)\wedge (B|K)]|(A|H)$;
\item[P2.]  $(A|H)\wedge (B|K)\leq    (B|K)|(A|H)$ and hence $\prev[(A|H)\wedge (B|K)]\leq \prev[(B|K)|(A|H)]$;
\item [P3.] $\prev[(A|H)\wedge (B|K)]=\prev[(B|K)|(A|H)]P(A|H)$  (\emph{compound formula for iterated conditional});
\item [P4.] if $A,B,H,K$ are logically independent events, denoting $P(A|H)=x$  and $P(B|H)=y$,  the  extension  $\mu$ on $(B|K)|(A|H)$ is coherent if and only if  $\mu \in[\mu',\mu'' ]$, where $\mu'$ and $\mu''$ are given in  formula 
(\ref{EQ:LU}).
\end{itemize}
Another classical property that can be checked for the iterated conditional is the  Import-Export principle. Given three events $B,K,A$, with $AK\neq \emptyset$, we say that   the Import-Export principle (\cite{McGe89}) is satisfied for the iterated conditional  $(B|K)|A$, if 
\begin{equation}\label{EQ:IE}
(B|K)|A = B|AK.
\end{equation} 
We also recall that the validity of the Import-Export principle (together with the validity of the total probability theorem) could lead to the counter-intuitive consequences related to Lewis’ triviality results (\cite{lewis76}, see also \cite{SGOP20}). 
Indeed, assuming that the total probability theorem holds for iterated conditionals, that is
\[
P(C|A)=P((C|A)\wedge C)+P((C|A)\wedge \no{C})=P((C|A)|C)P(C)+P((C|A)|\no{C})P(\no{C});
\]
if the Import-Export principle is valid, by applying (\ref{EQ:IE}) and by observing that $P(C|AC)=1$ and $P(C|A\no{C})=0$, it  follows that 
\begin{equation}\label{EQ:LEWIS}
P(C|A)=
P(C|AC)P(C)+P(C|A\no{C})P(\no{C})=P(C),
\end{equation}
 which of course is not valid in general for conditional events.
Then, the non validity of the Import-Export principle may avoid  Lewis’ triviality results. 

In sections \ref{SEC:Calabrese}, \ref{SEC:deFinetti}, \ref{SEC:Farrell},
we will check  the validity of  the previous properties  for notions of compound and  iterated conditional introduced in different trivalent logics  as suitable conditional events (in these cases, in properties P2 and P4, the symbol of prevision $\prev$ is replaced by the symbol of probability $P$, because the involved objects are conditional events). 
Then, in Section \ref{SEC:GENITER} we will check the basic properties   for the iterated conditionals, defined as conditional random quantities, built using the structure $
\Box|\bigcirc = \Box\wedge \bigcirc +\prev(\Box|\bigcirc)\no{\bigcirc}$ and the different notions of conjunction in trivalent logic  recalled in  Section \ref{SUBSEC:CONJUNCTIONS}.
\section{The Iterated conditional in the trivalent logic of  Cooper-Calabrese}\label{SEC:Calabrese}
In this section, in the framework of a trivalent logic, we consider the validity of the Import-Export principle and of the properties P1-P4 for the  notion of iterated conditional, here denoted by $(B|K)|_C (A|H)$,  studied by Cooper (\cite{Cooper1968}) and by  
Calabrese  (\cite{Cala1990}, see also \cite{Cala17,Cala2021}).\\
 We recall that   the notions of  conjunction and disjunction of conditionals used by Cooper and Calabrese  coincide with  $\wedge_S$ and $\vee_S$, respectively
 \footnote{ 
Note that Cantwell  (\cite{Cantwell2008})
defines the iterated conditional as Cooper and Calabrese, but in his trivalent logic conjunction and disjunction are defined by $\wedge_K$ and $\vee_K$, respectively.}
\begin{definition}\label{DEF:ITERCALA}
Given any pair of conditional events $A|H$ and $B|K$, the iterated conditional $(B|K)|_C (A|H)$ is defined as the following conditional event 
\begin{equation}\label{EQ:ITCAL}
(B|K)|_C (A|H)= B|(K \wedge (\no{H}\vee A)).
\end{equation}
\end{definition}
 We observe that in (\ref{EQ:ITCAL}) the conditioning event is the conjunction of the conditioning event $K$ of the consequent $B|K$ and the material conditional $\no{H}\vee A$ associated with  the antecedent $A|H$. 
\begin{remark}[Import-Export principle for $|_C$]
 By applying Definition \ref{DEF:ITERCALA}  with $H=\Omega$, it holds that \begin{equation}
 	\label{EQ:IMPEXPCAL}
 	(B|K)|_C A=ABK|AK=B|AK,
\end{equation}  
 which shows that the Import-Export principle is satisfied by  $|_C$.
 \end{remark}

\subsection{Property P1}
\noindent We observe that 
\begin{equation}\label{EQ:ITP1CAL}
\begin{array}{ll}
[(A|H)\wedge_S(B|K)]|_C(A|H)=
[(AHBK\vee AH\no{K}\vee \no{H}BK)|(H\vee K)]|_C(A|H)=\\
=(AHBK\vee AH\no{K}\vee \no{H}BK)|(AK\vee \no{H}K\vee AH\no{K}).
\end{array}
\end{equation}
From (\ref{EQ:ITCAL}) and (\ref{EQ:ITP1CAL}) it follows that
$(B|K)|_{C}(A|H)\neq  ((A|H)\wedge_S(B|K))|_C(A|H).$
Indeed, as illustrated by Table \ref{tab:nf_c}, 
when the constituent $AH\no{K}$ is true, it holds that $(B|K)|_{C}(A|H)$ is void, while  $[(A|H)\wedge_S(B|K)]|_C(A|H)$ is true.
Then, property P1 is not satisfied by the pair  $(\wedge_S,|_C)$.
\begin{table}[h]
	\centering
	\begin{tabular}{|c|c| c| c|}
\hline
 {$C_h$} & {$(A|H)\wedge_{S}(B|K)$} & {$(B|K)|_{C}(A|H)$} & {$[(A|H)\wedge_S(B|K)]|_C(A|H)$}\\
		\hline
\rule{0pt}{1em}
 {$AHBK\vee \no{H}BK$} & {True} & {True} & {True}\\
 {$AH\no{B}K\vee \no{H}\no{B}K$} & {False} & {False} & {False} \\
 {$AH\no{K}$} & {True} & {Void} & {True}\\
 {$\no{A}H$} & {False} & {Void} & {Void}\\
{$\no{H}\,\no{K}$} & {Void} & {Void} & {Void} \\
		\hline
	\end{tabular}
	\caption{Truth values of  $(A|H)\wedge_{S}(B|K)$, $(B|K)|_{C}(A|H)$,  and $[(A|H)\wedge_S(B|K)]|_C(A|H)$.}
	\label{tab:nf_c}
\end{table}
\subsection{Property P2}
From Table \ref{tab:nf_c} we also 
obtain that  property P2 is not satisfied by  $(\wedge_S, |_C)$. Indeed, when $AH\no{K}$ is true, it holds  that 
$(A|H)\wedge_S (B|K)$ is true, while  $(B|K)|_{C}(A|H)$ is void and hence $((A|H)\wedge_S (B|K))\nsubseteq (B|K)|_{C}(A|H).$\\

\subsection{Property P3}
Now we focus our attention on the following result regarding the coherence of a probability assessment on  $\{A|H, (B|K)|_{C}(A|H), (A|H)\wedge_{S}(B|K)\}$ (Theorem \ref{ThCalabrese}), then we use this result in order to check the validity of property P3 for the pair $(\wedge_S, |_C)$. 

\begin{theorem} Let $A$, $B$, $H$, $K$ be any logically independent events.
A probability assessment $\P=(x, y, z)$ on the family of conditional events $\F=\{A|H, (B|K)|_{C}(A|H), (A|H)\wedge_{S}(B|K)\}$ is coherent if and only if
$(x, y)\in[0,1]^2$ and $z \in [z', z'']$, where $z'=xy$ and $z''=max(x,y)$.
\label{ThCalabrese}
\end{theorem}
\begin{proof}
	The constituents $C_h$'s and the points $Q_h$'s associated with the assessment $\P=(x,y,z)$ on $\F$ are (see also Table \ref{tab:th7})
	\[
	\begin{array}{l}
		C_1=AHBK, C_2=\no{A}HBK\vee \no{A}H\no{B}K\vee \no{A}H\,\no{K}=\no{A}H, C_3=\no{H}BK, C_4=AH\no{B}K,\\ C_5=\no{H}\no{B}K, C_6=AH\no{K}, C_0=\no{H}\,\no{K},
	\end{array}
	\]
	and
	\[
	\begin{array}{l}
		Q_1=(1,1,1),\ Q_2=(0,y,0),\ Q_3=(x,1,1),Q_4=(1,0,0),\\ Q_5=(x,0,0),\ Q_6=(1,y,1),  \P=Q_0=(x,y,z).
	\end{array}
	\]
		\begin{table}[h]
		\centering
		\begin{tabular}{|l|c|c c c| c|}
			\hline
			{} & {$C_h$} & {$A|H$} & {$(B|K)|_{C}(A|H)$} & {$(A|H)\wedge_{S}(B|K)$} & {$Q_h$}\\
			\hline
			{$C_1$} & {$AHBK$} & {1} & {1} & {1} & {$Q_1$} \\
			{$C_2$} & {$\no{A}H$} & {0} & {$y$} & {0} & {$Q_2$} \\
			{$C_3$} & {$\no{H}BK$} & {$x$} & {1} & {1} & {$Q_3$} \\
			{$C_4$} & {$AH\no{B}K$} & {1} & {0} & {0} & {$Q_4$} \\
			{$C_5$} & {$\no{H}\no{B}K$} & {$x$} & {$0$} & {0} & {$Q_5$} \\
			{$C_6$} & {$AH\no{K}$} & {$1$} & {$y$} & {$1$} & {$Q_6$} \\
			{$C_0$} & {$\no{H}\,\no{K}$} & {$x$} & {$y$} & {$z$} & {$Q_0$} \\
			\hline
		\end{tabular}
		\caption{Constituents $C_h$ and points $Q_h$ associated with $\F=\{A|H, (B|K)|_{C}(A|H), (A|H)\wedge_{S}(B|K)\}$, the probability assessment $\P=(x, y, z)$.}    
		\label{tab:th7}
	\end{table}
	The system $(\Sigma)$ in (\ref{SYST-SIGMA}) associated with the pair $(\F,\P)$  becomes\\
	\begin{equation}
		\left\{
		\begin{array}{llll}
			\lambda_1+x\lambda_3+\lambda_4+x\lambda_5+\lambda_6=x, \\
			\lambda_1+y\lambda_2+\lambda_3+y\lambda_6=y,\\
			\lambda_1+\lambda_3+\lambda_6=z,\\
			\lambda_1+\lambda_2+\lambda_3+\lambda_4+\lambda_5+\lambda_6=1, \\ \lambda_i\geq0 \ \forall i=1, \dots, 6.\\
		\end{array}
		\label{sis}
		\right.    
	\end{equation}
	\paragraph{Lower bound}
	We  first prove that the assessment $(x, y, xy)$ 
	is coherent for every 
	$(x,y)\in[0,1]^2$. Then, in order to prove that $z'=xy$ is the lower bound for $z=P((A|H)\wedge_{S}(B|K))$,
	we verify that the assessment  $(x,y,z)$ is not coherent when $z<z'=xy$.\\
	We observe that $\P=(x,y,xy)=xyQ_1+(1-x)Q_2+x(1-y)Q_4$, so a solution of (\ref{sis}) is given by $\Lambda=(xy, 1-x,0, x(1-y), 0,0).$\\
	Then, by considering the  function  $\phi$ as defined in (\ref{EQ:I0}),  it holds that 
	\[
	\phi_{1}(\Lambda)=\sum_{h:C_h\subseteq H}\lambda_h=\lambda_1+\lambda_2+\lambda_4+\lambda_6=xy+(1-x)+x(1-y)=1>0,
	\]
	\[
	\phi_{2}(\Lambda)=\sum_{h:C_h\subseteq (A\vee \no{H}) \lor K}\lambda_h=\lambda_1+\lambda_3+\lambda_4+\lambda_5=xy+x(1-y)=x,
	\]
	\[
	\phi_{3}(\Lambda)=\sum_{h:C_h\subseteq H\lor K}\lambda_h=\lambda_1+\cdots+\lambda_6=1>0.
	\]
	Let $\mathcal{S}'=\{(xy, 1-x,0, x(1-y), 0,0)\}$  denote a subset of the set $\mathcal{S}$ of all solutions of (\ref{sis}). We have that $M_1'=1$, $M_2'=x$, $M_3'=1$  (as defined in (\ref{EQ:I0'})). 
	We distinguish two cases: $(i)$ $x>0$, $(ii)$ $x=0$. In the case $(i)$ we get 
	$M_1'>0$, $M_2'>0$, $M_3'>0$ and then $I'_0=\emptyset$. By Theorem \ref{CNES-PREV-I_0-INT}, the assessment $(x,y,xy)$ is coherent $\forall (x,y)\in[0,1]^2$. 
	In the case $(ii)$ we have that $M_1'>0$, $M_2'=0$, $M_3'>0$, hence $I_0'=2$. We observe that  the sub-assessment $\P'_0=y$ on $\F'_0=\{(B|K)|_C (A|H)\}$ is coherent for every $y\in[0,1]$. Then, by Theorem \ref{CNES-PREV-I_0-INT}, the assessment $(x,y,xy)$  on $\F$ is coherent $\forall (x,y)\in[0,1]^2$.
	\\

	In order to verify that $z'=xy$ is the lower bound for $z$, we observe that the points $Q_1, Q_2, Q_4$ belong to the plane $\pi: yX+Y-Z=y$, where $X, Y, Z$ are the axes coordinates.\\
	Now, by considering the function $f(X,Y,Z)= yX+Y - Z$, we observe that for each constant $k$ the equation $f(X,Y,Z)=k$ represents a plane which is
	parallel to $\pi$ and coincides with $\pi$ when $k=y$. We also observe that $f(Q_1)=f(Q_2)=f(Q_4)=y$, $f(Q_3)= f(Q_5)=xy \leq y$ and $f(Q_6)= f(1,y,1)= y+y-1=2y-1 \leq y$.\\
	Then, for every $\P=\sum_{h=1}^{6} \lambda_h Q_h$, with  $\sum_{h=1}^{6} \lambda_h=1$ and $\lambda_h\geq 0$, that is $\P \in \I$, it holds that $f(\P)=f(\sum_{h=1}^{6} \lambda_h Q_h)=
	\sum_{h=1}^{6}\lambda_h f(Q_h)\leq y$. On the other side, given $a>0$, by considering $\P=(x,y,xy-a)$ it holds that $f(\P)=y+a> y$ and hence  $\P=(x,y,xy-a) \notin \I$.
	Therefore, for any given $a>0$ the assessment $(x,y,xy-a)$ is not coherent because  $(x,y,xy-a) \notin \I$. Then the lower bound of $z=P((A|H)\wedge_S(B|K))$ is $z'=xy$.\\
	\paragraph{Upper bound}
	We  verify that the assessment  $(x,y,\max(x,y))$  on $\F$
	is coherent  for every  $(x,y)\in[0,1]^2$. 
	Moreover, we show that  $z''=\max(x,y)$ is the upper bound for $z=P((A|H)\wedge_S(B|K))$ by showing that any assessment $(x,y,z)$ on $\F$ with $(x,y)\in[0,1]^2$ and  $z>\max\{x,y\}$ is not coherent.
	We distinguish two cases: $(i)$ $x\geq y$, $(ii)$ $x<y$.\\
	$(i)$ We have that $\max(x,y)=x$ and hence
	$$\P=(x,y,x)=(1-x)Q_2+xQ_6.$$
	Then, the vector $\Lambda=(0,1-x,0,0,0,x)$ is a solution of (\ref{sis}). Moreover,  it holds that 
	\[
	\phi_{1}(\Lambda)=\sum_{h:C_h\subseteq H}\lambda_h=\lambda_1+\lambda_2+\lambda_4+\lambda_6=1-x+x=1>0,
	\]
	\[
	\phi_{2}(\Lambda)=\sum_{h:C_h\subseteq (K\wedge(A\lor \bar{H}))}\lambda_h=\lambda_1+\lambda_3+\lambda_4+\lambda_5=0,
	\]
	\[
	\phi_{3}(\Lambda)=\sum_{h:C_h\subseteq (H\lor K)}\lambda_h=\lambda_1+\cdots+\lambda_6=1>0.
	\]
	Let $\mathcal{S}'=\{(0,1-x,0,0,0,x)\}$  denote a subset of the set $\mathcal{S}$ of all solutions of (\ref{sis}). We have that $M_1'=1$, $M_2'=0$, $M_3'=1$. It follows that $\I'_0=2$. As the sub-assessment $\P'_0=y$ on $\F'_0=\{(B|K)|_C (A|H)\}$ is coherent $\forall y \in [0,1]$, { by Theorem \ref{CNES-PREV-I_0-INT}}, it follows that the assessment $(x,y,\max(x,y))$ is coherent.\\
	In order to verify that $z''=\max(x,y)=x$ is the upper bound for $z$, 
	we observe that  if $\max\{x,y\}=1$, then $(1,y,z)$ with $z>1$ is incoherent. Let us assume that $y\leq x<1$.
	We observe that the points $Q_2, Q_3, Q_6$ belong to the plane $\pi:X+ \frac{1-x}{1-y}Y-Z
	=\frac{y(1-x)}{1-y}
	$, where $X, Y, Z$ are the axes coordinates.\\
	Now, by considering the function $f(X,Y,Z)=
	X+ \frac{1-x}{1-y}Y-Z
	-\frac{y(1-x)}{1-y}
	$, we observe that 
	$f(Q_1)=1-x>0$,  $f(Q_2)=f(Q_3)=f(Q_6)=0$, $f(Q_4)=1-y\frac{1-x}{1-y}\geq 0$, $f(Q_5)=\frac{x-y}{1-y}\geq 0$.
	Then, for every $\P=\sum_{h=1}^{6} \lambda_h Q_h$, with  $\sum_{h=1}^{6} \lambda_h=1$ and $\lambda_h\geq 0$, that is $\P \in \I$, it holds that $f(\P)=f(\sum_{h=1}^{6} \lambda_h Q_h)=
	\sum_{h=1}^{6}\lambda_h f(Q_h)\geq 0$. On the other side, given $z>x$, by considering $\P=(x,y,z)$ it holds that $f(\P)=x-z<0$ and hence  $\P=(x,y,z) \notin \I$.
	Therefore, for any given $z>x$ the assessment $(x,y,z)$ is not coherent because  $(x,y,z) \notin \I$. Then the upper bound on $z$ is $z''=z=\max\{x,y\}$.\\
	
	$(ii)$   In this case $\max(x,y)=y$. We prove that the assessment $(x,y,\max(x,y))$ is coherent. We observe that 
	$$(x,y,y)=yQ_3+(1-y)Q_5.$$
	Then, the vector $\Lambda=(0,0,y,0,1-y,0)$ is a solution of (\ref{sis}).
	We have
	\[
	\phi_{1}(\Lambda)=\sum_{h:C_h\subseteq H}\lambda_h=0,
	\]
	\[
	\phi_{2}(\Lambda)=\sum_{h:C_h\subseteq (K\wedge(A\lor \bar{H}))}\lambda_h=1>0,
	\]
	\[
	\phi_{3}(\Lambda)=\sum_{h:C_h\subseteq (H\lor K)}\lambda_h=1>0.
	\]
	Let $\mathcal{S}'=\{(0,0,y,0,1-y,0)\}$  denote a subset of the set $\mathcal{S}$ of all solutions of (\ref{sis}). We have that $M_1'=0$, $M_2'=1$, $M_3'=1$. It follows that $\I'_0=1$.
The sub-assessment $x$ on $\{A|H\}$ is coherent $\forall x \in [0,1]$. Then, by Theorem \ref{CNES-PREV-I_0-INT}, the assessment $(x,y,\max(x,y))$  on $\F$ is coherent $\forall (x,y)\in[0,1]^2$.
	
	In order to verify that $z''=\max(x,y)=y$ is the upper bound for $z$, 
	we observe that  if $\max\{x,y\}=1$, then $(x,1,z)$ with $z>1$ is incoherent. Let us assume that $x\leq y<1$.
	We observe that the points $Q_3, Q_5, Q_6$ belong to the plane $\pi: \frac{1-y}{1-x}X-\frac{x(1-y)}{1-x}+Y-Z=0$, where $X, Y, Z$ are the axes coordinates.\\
	Now, by considering the function $f(X,Y,Z)= \frac{1-y}{1-x}X-\frac{x(1-y)}{1-x}+Y-Z$, we observe that 
	$f(Q_1)=f(Q_4)=1-y>0$, $f(Q_2)=\frac{y-x}{1-x}> 0$ $f(Q_3)=f(Q_5)=f(Q_6)=0$.
	Then, for every $\P=\sum_{h=1}^{6} \lambda_h Q_h$, with  $\sum_{h=1}^{6} \lambda_h=1$ and $\lambda_h\geq 0$, that is $\P \in \I$, it holds that $f(\P)=f(\sum_{h=1}^{6} \lambda_h Q_h)=
	\sum_{h=1}^{6}\lambda_h f(Q_h)\geq 0$. On the other side, given $z>y$, by considering $\P=(x,y,z)$ it holds that $f(\P)=y-z<0$ and hence  $\P=(x,y,z) \notin \I$.
	Therefore, for any given $z>y$ the assessment $(x,y,z)$ is not coherent because  $(x,y,z) \notin \I$. Then the upper bound on $z$ is $z''=y=\max\{x,y\}$.\\
	Finally, for each given $(x,y)\in[0,1]$,  as $(x,y,z')$ and $(x,y,z'')$ are coherent, 
	by the Fundamental theorem of probability, 
 any assessment $(x,y,z)$ with $z\in[z',z'']$ is coherent too. Then, the assessment $(x,y,z)$  on $\F$ is coherent for every $(x,y)\in[0,1]^2$ and $z\in[z',z'']$.

\end{proof}
\begin{remark}[Property P3]
From {Theorem \ref{ThCalabrese}} any probability assessment $(x,y,z)$ on $\F=\{A|H, (B|K)|_{C}(A|H),(A|H)\wedge_{S}(B|K)\}$, with $(x,y)\in[0,1]^2$ and  $xy\leq z\leq max(x,y)$ is coherent.
Thus, as $z=xy$ is not the unique coherent extension of the conjunction  $(A|H)\wedge_S (B|K)$, in general the quantity $
P[(A|H)\wedge_S (B|K)]$ do not coincide  with the product $P[(B|K)|_{C}(A|H)]P(A|H)$. 
For example,  it could be that  $
P[(B|K)|_{C}(A|H)]P(A|H)=0<
P[(A|H)\wedge_S (B|K)]=1$, because  the assessment $(1,0,1)$ on $\F$ is coherent (while it is not coherent on $\{A,B|A,AB\}$).
Then, property P3 is not satisfied by the pair $(\wedge_S,|_C)$.
\end{remark}

\subsection{Property P4}
We check the validity of property P4 for the iterated conditioning $|_C$ by studying  the set of all coherent probability assessments on the family $\{A|H, B|K, (B|K)|_{C}(A|H)\}$ (Theorem \ref{ThCalabrese2}).

\begin{theorem}
Let $A$, $B$, $H$, $K$ be any logically independent events.
The probability assessments  $\P=(x, y, z)$ on the family of conditional events $\F=\{A|H, B|K, (B|K)|_{C}(A|H)\}$
is  coherent for every  $(x,y,z)\in[0,1]^3$.
\label{ThCalabrese2}
\end{theorem}
 \begin{proof}
We recall that $(B|K)|_{C}(A|H)=B|(K\wedge (\no{H}\vee A))$. Then, 
the constituents $C_h$'s and the points $Q_h$'s associated with the assessment $\P=(x,y,z)$ on $\F=\{A|H, B|K, B|(K\wedge (\no{H}\vee A))\}$ are (see also Table \ref{tab:th8})
\[ 
\begin{array}{l}
C_1=AHBK,\  C_2=AH\no{B}K, \ C_3=AH\no{K},\ C_4=\no{A}HBK,\ C_5=\no{A}H\no{B}K,\\ C_6=\no{A}H\no{K},\ C_7=\no{H}BK,\ C_8=\no{H}\no{B}K,\ C_0=\no{H}\,\no{K},
\end{array}
\]
and
\[
\begin{array}{l}
Q_1=(1,1,1), Q_2=(1,0,0), Q_3=(1,y,z), Q_4=(0,1,z), Q_5=(0,0,z), \\
Q_6=(0,y,z), Q_7=(x,1,1), Q_8=(x,0,0), \P=Q_0=(x,y,z).
\end{array}
\]
\begin{table}[h]
    \centering
    \begin{tabular}{|l|c|c c c| c|}
    \hline
        {} & {$C_h$} & {$A|H$} & {$B|K$} & {$(B|K)|_{C}(A|H)$} & {$Q_h$}\\
        \hline
        {$C_1$} & {$AHBK$} & {1} & {1} & {1} & {$Q_1$} \\
        {$C_2$} & {$H\no{B}K$} & {1} & {0} & {0} & {$Q_2$} \\
        {$C_3$} & {$AH\no{K}$} & {1} & {$y$} & {$z$} & {$Q_3$} \\
        {$C_4$} & {$\no{A}HBK$} & {0} & {1} & {$z$} & {$Q_4$} \\
        {$C_5$} & {$\no{A}H\no{B}K$} & {0} & {0} & {$z$} & {$Q_5$} \\
        {$C_6$} & {$\no{A}H\no{K}$} & {0} & {$y$} & {$z$} & {$Q_6$} \\
        {$C_7$} & {$\no{H}BK$} & {$x$} & {$1$} & {1} & {$Q_7$} \\
        {$C_8$} & {$\no{H}\no{B}K$} & {$x$} & {$0$} & {0} & {$Q_8$} \\
        {$C_0$} & {$\no{H}\,\no{K}$} & {$x$} & {$y$} & {$z$} & {$Q_0$} \\
    \hline
    \end{tabular}
\caption{Constituents and points $Q_h$' associated with $\F=\{A|H, B|K, (B|K)|_{C}(A|H)\}$, the probability assessment $\P=(x, y, z)$.}    
\label{tab:th8}
\end{table}
The system $(\Sigma)$ in (\ref{SYST-SIGMA}) associated with the pair $(\F,\P)$  becomes\\
\begin{equation}
\left\{
    \begin{array}{llll}
     \lambda_1+\lambda_2+\lambda_3+x\lambda_7+x\lambda_8=x, \\
     \lambda_1+y\lambda_3+\lambda_4+y\lambda_6+\lambda_7=y,\\
     \lambda_1+z\lambda_3+z\lambda_4+z\lambda_5+z\lambda_6+\lambda_7=z,\\
    \lambda_1+\lambda_2+\lambda_3+\lambda_4+\lambda_5+\lambda_6+\lambda_7+\lambda_8=1, \\ \lambda_i\geq0 \ \forall i=1, \dots, 8.\\
    \end{array}
    \label{sc}
\right.    
\end{equation}

\begin{figure}[tbph]
\centering
\includegraphics[width=0.75\linewidth]{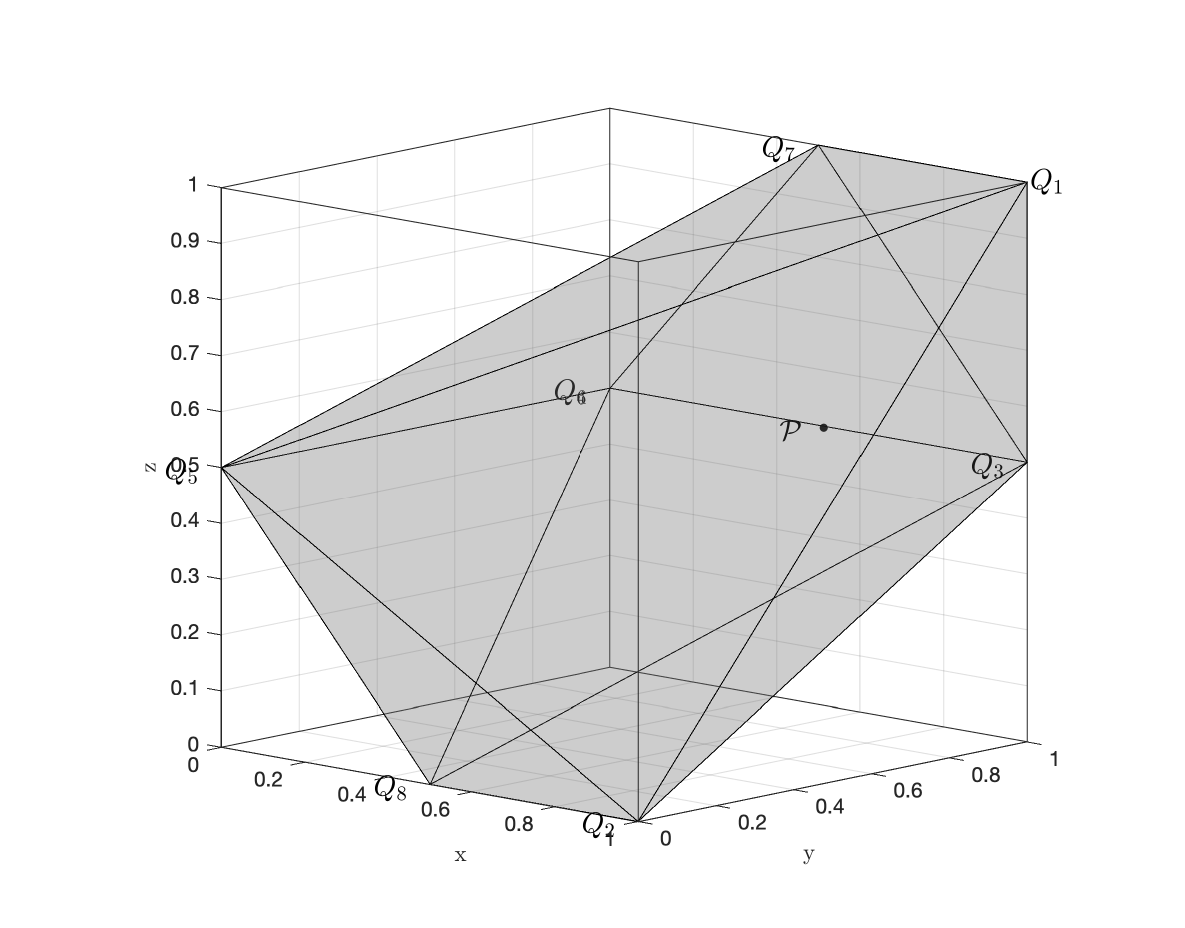}
\caption{Convex hull of  the points $Q_i$ for $i=1, \dots,8$ associated with the pair $(\F,\P)$, where $\P=(x,y,z)$ and $\F=\{A|H, B|K, (B|K)|_{C}(A|H)\}$.  
In the figure the numerical  values are: $x=0.5$, $y=1$, $z=0.5$.}
\label{FIG:fig_it_c}
\end{figure}

We observe that $\P$ belongs to the segment with end points $Q_3$, $Q_6$; indeed $(x,y,z)=xQ_3+(1-x)Q_6=x(1,y,z)+(1-x)(0,y,z)$.
The vector $\Lambda=(0,0,x,0,0,1-x,0,0)$ is a solution of (\ref{sc}), with
\[
\phi_{1}(\Lambda)=\sum_{h:C_h\subseteq H}\lambda_h=\lambda_1+\lambda_2+\lambda_3+\lambda_4+\lambda_5+\lambda_6=1>0,
\]
\[
\phi_{2}(\Lambda)=\sum_{h:C_h\subseteq K}\lambda_h=\lambda_1+\lambda_2+\lambda_4+\lambda_5+\lambda_7+\lambda_8=0,
\]
\[
\phi_{3}(\Lambda)=\sum_{h:C_h\subseteq ((A \vee \no{H})\wedge K)}\lambda_h=\lambda_1+\lambda_2+\lambda_7+\lambda_8=0.
\]

We set $\mathcal{S'}=\{\Lambda\}$ and we 
 get $\I_0'=\{2,3\}$. So we obtain $\F_0'=\{B|K, B|(K\wedge (\no{H}\vee A))\}$ and $\P_0'=(y, z)$. By recalling Example \ref{EX:EX1} the sub-assessment $(y, z)$ on $\{B|K, B|(K\wedge (\no{H}\vee A))\}$ is coherent for every $(y, z)\in[0,1]^2$. Thus, by  Theorem \ref{CNES-PREV-I_0-INT'}, the assessment $(x,y,z)$  on $\F$ is coherent $\forall (x,y,z)\in[0,1]^3$.
\end{proof}

\begin{remark}[Property P4]
We observe that  the probability propagation rule  valid for unconditional events (Property P4)  is no longer valid for Calabrese's iterated conditional. Indeed,  from {Theorem \ref{ThCalabrese2}}, any probability assessment $(x,y,z)$ on $\F=\{ A|H,B|K, (B|K)|_{C}(A|H)\}$, with $(x,y,z)\in[0,1]^3$ is coherent. For instance, the assessment $(1,1,0)$ on $\F$ is coherent, while it is not coherent on $\{A, B, B|A\}$.
\end{remark}

\section{The Iterated conditional in the trivalent logic of  {d}e Finetti.}\label{SEC:deFinetti}
In this section we analyze the notion of iterated conditional introduced by de Finetti in \cite{defi36}. 
After recalling that the notion of  conjunction and disjunction of conditionals introduced  by de Finetti in \cite{defi36} coincide with  $\wedge_K$ and $\vee_K$ (see Section \ref{SUBSEC:CONJUNCTIONS}),  we check the validity of the Import-Export principle and of the properties P1-P4 for this iterated conditional in the corresponding trivalent logic. \\

\begin{definition}
\label{DEF:ITERDEF}
Given any pair of conditional events $A|H$ and $B|K$, de Finetti iterated conditional, denoted by $(B|K)|_{dF} (A|H)$,  is defined as
\begin{equation}\label{EQ:ITdf}
(B|K)|_{dF} (A|H)= B|(AHK).
\end{equation}
\end{definition}
\begin{remark}[Import-Export principle for $|_{dF}$]
By applying Definition \ref{DEF:ITERDEF}  with $H=\Omega$, it holds that 
\begin{equation}
\label{EQ:IMPEXPDF}
(B|K)|_{dF}A=B|AK,
\end{equation}
which shows that the Import-Export principle \cite{McGe89} is satisfied by $|_{dF}$. 
Then,  from (\ref{EQ:IMPEXPCAL}), it follows that 
\[
(B|K)|_{C}A=(B|K)|_{dF}A=B|AK.
\]
\end{remark}
\subsection{Property P1}
\noindent To check property P1 we observe that from (\ref{EQ:ITdf}) it holds that
\begin{equation}\label{EQ:ITP1dF}
\begin{array}{ll}
    [(A|H)\wedge_K(B|K)]|_{dF}(A|H)=
    [AHBK|(HK\lor \no{A}H \lor \no{B}K)]|_{dF}(A|H)=\\
  AHBK|(AHK\lor AH\no{B}K)=
    AHBK|AHK=(B|K)|_{dF}(A|H).
\end{array}
\end{equation}
Then, property P1 is satisfied by the pair $(\wedge_K, |_{dF})$ (see  also Table \ref{tab:rel2}).
\begin{table}[]
    \centering
    \begin{tabular}{|c|c|c|c|}
    \hline
        {$C_h$} & {$(A|H)\wedge_K (B|K)$} & {$(B|K)|_{dF}(A|H)$} & {$[(A|H)\wedge_K(B|K)]|_{dF}(A|H)$} \\
        \hline
        {$AHBK$} & {True} & {True} & {True} \\
        {$AH\no{B}K$} & {False} & {False} & {False} \\
        {$AH\no{K}\vee \no{H}BK \vee \no{H}\,\no{K}$} & {Void} & {Void} & {Void} \\
        {$\no{A}H \vee \no{H}\no{B}K$} & {False} & {Void} & {Void} \\
    \hline
    \end{tabular}
\caption{Truth table of $(A|H)\wedge_K (B|K),\ (B|K)|_{dF}(A|H)$, and $[(A|H)\wedge_K(B|K)]|_{dF}(A|H)$.}    \label{tab:rel2}
\end{table}

\subsection{Property P2}
From Table \ref{tab:rel2} we also observe that relation P2 is satisfied by $(\wedge_K, |_{dF})$.
Indeed, according to (\ref{EQ:GN}), if $(A|H)\wedge_K (B|K)$ is true, then  $(B|K)|_{dF}(A|H)$
is  true;  
if  $(B|K)|_{dF}(A|H)$ is false, then  $(A|H)\wedge_K (B|K)$ is  false.

\subsection{Property P3}
We consider now the following results regarding the coherence of a probability assessment on  $\{A|H, (B|K)|_{dF}(A|H), (A|H)\wedge_{K}(B|K)\}$ (Theorem \ref{ThdeFinetti}) in order to check the validity of property P3 for the pair $(\wedge_K, |_{dF})$.

\begin{theorem} Let $A$, $B$, $H$, $K$ be any logically independent events.
A probability assessment $\P=(x, y, z)$ on the family of conditional events $\F=\{A|H, (B|K)|_{dF}(A|H), (A|H)\wedge_{K}(B|K)\}$ is coherent if and only if 
$(x, y)\in[0,1]^2$ and $z \in [z', z'']$, where $z'=0$ and $z''=xy$.
\label{ThdeFinetti}
\end{theorem}
 \begin{proof}
See  \ref{SEC:APPdf1}.
\end{proof}

\begin{remark}[Property P3]
From {Theorem \ref{ThdeFinetti}} any probability assessment $(x,y,z)$ on $\F=\{A|H, (B|K)|_{dF}(A|H),(A|H)\wedge_{K}(B|K)\}$, with $(x,y)\in[0,1]^2$ and  $0\leq z\leq xy$,   is coherent.
Thus, as $z=xy$ is not the unique coherent extension of the conjunction  $(A|H)\wedge_K (B|K)$,  the quantity $
P[(A|H)\wedge_K (B|K)]$ could not coincide with the product $P[(B|K)|_{dF}(A|H)]P(A|H)$. For example, if we choose the probability assessment $\P=(1,1,0)$, we observe that $\P$ is  coherent on $\F$ because $0\in[z',z'']=[0,1]$. However, we observe that $\P$
 on $\{A, B|A, AB\}$ is not coherent 
because $P(AB)=0\neq P(B|A)P(A)$.\\
Then, property P3 is not satisfied by the pair $(\wedge_K,|_{dF})$.
\end{remark}
\subsection{Property P4}
To check the validity of property P4 for the iterated conditioning $|_{dF}$ we study the set of all  coherent probability assessments on the family $\{A|H, B|K, (B|K)|_{dF}(A|H)\}$ (Theorem \ref{ThdeFinetti2}).
\begin{theorem} 
Let $A$, $B$, $H$, $K$ be any logically independent events.
The probability assessments  $\P=(x, y, z)$ on the family of conditional events $\F=\{A|H, B|K, (B|K)|_{dF}(A|H)\}$
is  coherent for every  $(x,y,z)\in[0,1]^3$.
\label{ThdeFinetti2}
\end{theorem}
\begin{proof}
The proof can be found in  \ref{SEC:APPdf2}.
\end{proof}

\begin{remark}[Property P4]
We observe that  the probability propagation rule valid for unconditional events P4 is no longer valid for de Finetti's iterated conditional.
Indeed,  from {Theorem \ref{ThdeFinetti2}}, any probability assessment $(x,y,z)$ on $\F=\{A|H, B|K, (B|K)|_{dF}(A|H)\}$, with $(x,y,z)\in[0,1]^3$ is coherent. For instance, the assessment $(1,1,0)$ is coherent on $\F$ but it is not coherent on $\{A, B, B|A\}$.
\end{remark}

\section{The Iterated conditional in the trivalent logic of  Farrell.}\label{SEC:Farrell}
In this section we consider another definition of iterated conditional as a suitable conditional event which was proposed by Farrell in \cite[p. 385]{Farrell1979} (see also \cite{EgRS20}). 
In his trivalent logic, the author uses "$\wedge_K$" and "$\vee_K$"  as conjunction and disjunction of conditional events (Section \ref{SUBSEC:CONJUNCTIONS}), respectively.
We first recall the definition of the iterated conditional, then we check   the validity of the Import-Export principle and  of the properties P1-P4. 
\begin{definition}
\label{DEF:ITERF}
Given any pair of conditional events $A|H$ and $B|K$, Farrell iterated conditional, here denoted by $(B|K)|_{F}(A|H)$, is defined as the conditional event
\begin{equation}\label{EQ:ITERF}
    (B|K)|_{F}(A|H)=AHBK|(AHBK \vee AH\no{B}K \vee \no{H}\no{B}K).
\end{equation}
\end{definition}
\begin{remark}[Import-Export principle for $|_F$]
By applying (\ref{EQ:ITERF}) with $H=\Omega$, 
it holds that 
\[
(B|K)|_F A=ABK|(ABK \vee A\no{B}K)=B|AK.
\]
Then,  as $(B|K)|_F A=B|AK$, the Import-Export principle  is satisfied by $|_F$. Moreover, by recalling  (\ref{EQ:IMPEXPCAL}) and (\ref{EQ:IMPEXPDF}),  it follows that 
$$(B|K)|_FA=(B|K)|_{C}A=(B|K)|_{dF}A=B|AK.$$
\end{remark}

\subsection{Property P1}
By recalling that $(A|H)\wedge_K(B|K)=AHBK|(AHBK\lor \no{A}H \lor \no{B}K)$, from  (\ref{EQ:ITERF}) it follows that 
\begin{equation}\label{EQ:ITP1F}
\begin{array}{ll}
[(A|H)\wedge_K(B|K)]|_{F}(A|H)=
    [AHBK|(AHBK\lor \no{A}H \lor \no{B}K)]|_{F}(A|H)=\\
  AHBK|(AHBK \vee AH\no{B}K \vee \no{H}\no{B}K)=(B|K)|_{F}(A|H).
\end{array}
\end{equation}
Then, as $   [(A|H)\wedge_K(B|K)]|_{F}(A|H)=(B|K)|_{F}(A|H)$, 
 Property P1 is satisfied by the pair $(\wedge_K, |_{F})$.  This relation can also be obtained by observing that  the truth values of $(B|K)|_{F}(A|H)$ and of $[(A|H)\wedge_K(B|K)]|_{F}(A|H)$ in Table \ref{tab:relF1} coincide. 

\begin{table}[]
    \centering
    \begin{tabular}{|c|c|c|c|}
    \hline
        {$C_h$} & {$(A|H)\wedge_K (B|K)$} & {$(B|K)|_{F}(A|H)$} & {$[(A|H)\wedge_K(B|K)]|_{F}(A|H)$} \\
        \hline
        {$AHBK$} & {True} & {True} & {True} \\
        {$AH\no{B}K \vee \no{H}\no{B}K$} & {False} & {False} & {False} \\
        {$AH\no{K}\vee \no{H}BK \vee \no{H}\,\no{K}$} & {Void} & {Void} & {Void} \\
        {$\no{A}H$} & {False} & {Void} & {Void} \\
    \hline
    \end{tabular}
\caption{Truth table of $(A|H)\wedge_K (B|K),\ (B|K)|_{F}(A|H)$, and $[(A|H)\wedge_K(B|K)]|_{F}(A|H)$.}    \label{tab:relF1}
\end{table}

\subsection{Property P2}
We observe that (see  Table \ref{tab:relF1}) the iterated conditional  
$(B|K)|_F(A|H)$ is true when the cojunction  $(A|H)\wedge_K (B|K)$ is true; moreover, the condjunction  $(A|H)\wedge_K (B|K)$ is false when the iterated conditional $(B|K)|_F(A|H)$ is false.
Then, $(A|H)\wedge_K (B|K)\subseteq (B|K)|_F(A|H)$. Thus, it follows from  (\ref{EQ:LEQ}) that $(A|H)\wedge_K (B|K)\leq (B|K)|_F(A|H)$, which means that 
P2 is satisfied by $(\wedge_K, |_F)$.

\subsection{Property P3}
To check the validity of  property P3 
for $(B|K)|_{F}(A|H)$  we study the set of coherent probability assessment on the family 
$\F=\{A|H, (B|K)|_{F}(A|H), (A|H)\wedge_{K} (B|K)\}.$
\begin{theorem}
    Let $A$, $B$, $H$, $K$, be any logically independent events. A probability assessment $\P=(x, y,z)$ on the family of conditional events $\F=\{A|H, (B|K)|_{F}(A|H), (A|H)\wedge_{K} (B|K)\}$ is coherent if and only if $(x,y)\in[0,1]^2$ and $z\in[z', z'']$, where $z'=0$ and  $z''=T^{H}_0(x,y)$,
 where 
 \[
T^{H}_0(x,y)=\begin{cases}
                0, & \text{if } x=0 \text{ or } y=0\\
             \frac{xy}{x+y-xy}, &  \text{if } x\neq 0 \text{ and } y\neq 0,\\
         \end{cases}
 \]   
 is the Hamacher t-norm with parameter $\lambda=0$.
    \label{ThFarrell}
\end{theorem}
\begin{proof}
See \ref{SEC:APPF1}.
\end{proof}

\begin{remark}[Property P3]
From {Theorem \ref{ThFarrell}} any probability assessment $(x,y,z)$ on $\F=\{A|H, (B|K)|_{F}(A|H),(A|H)\wedge_{K}(B|K)\}$, with $(x,y)\in[0,1]^2$ and $z\in[z', z'']$, where $z'=0$ and  $z''=T^{H}_0(x,y)$ is coherent. We observe that $z=xy \in [z', z'']$, then $z=xy$ is not the unique coherent extension of the conjunction  $(A|H)\wedge_K (B|K)$,  the quantity $
P[(A|H)\wedge_K (B|K)]$ could not coincide with the product $P[(B|K)|_{F}(A|H)]P(A|H)$. For example, if we choose the probability assessment $\P=(1,1,0)$, we observe that  $\P$  is coherent on $\F$, because $0\in[z',z'']=[0,1]$. However, we observe that the assessment  $\P=(1,1,0)$ on $\{A, B|A, AB\}$ is not coherent, because $P(AB)=0\neq P(B|A)P(A)=1$.
Then, property P3 is not satisfied by the pair $(\wedge_K,|_{F})$.
\end{remark}

\subsection{Property P4}
To check the validity of property P4 for the iterated conditioning $|_{F}$ we study the set of all coherent probability assessments on the family $\{A|H, B|K, (B|K)|_{F}(A|H)\}$. 
\begin{theorem} 
Let $A$, $B$, $H$, $K$ be any logically independent events.
The probability assessments  $\P=(x, y, z)$ on the family of conditional events $\F=\{A|H, B|K, (B|K)|_{F}(A|H)\}$
is  coherent for every  $(x,y,z)\in[0,1]^3$.
\label{ThFarrell2}
\end{theorem}
\begin{proof}
See \ref{SEC:APPF2}.
\end{proof}

\begin{remark}[Property P4]
We observe that  the probability propagation rule valid for unconditional events P4 is no longer valid for Farrell's iterated conditional.
Indeed,  from {Theorem \ref{ThFarrell2}}, any probability assessment $(x,y,z)$ on $\F=\{A|H, B|K, (B|K)|_{F}(A|H)\}$, with $(x,y,z)\in[0,1]^3$ is coherent. For instance, the assessment $(1,1,0)$ is coherent on $\F$ but it is not coherent on $\{A, B, B|A\}$.
\end{remark}
\section{Iterated conditionals and compound prevision theorem}\label{SEC:GENITER}
We observe that none of the iterated conditioning operations, $|_C$, $|_{dF}$, and $|_F$, studied in the previous sections  satisfy the compound probability theorem  P3. 
In this section,  we consider  iterated conditionals  which, among other things, satisfy  property 
P3. After recalling a structure used for a notion of iterated conditional defined in the framework of conditional random quantities, we introduce four notions of iterated conditional which are  based on the same structure and on the four conjunctions of trivalent logics recalled in Section \ref{SEC:2}. Then,  for all these new objects we check the validity of the basic properties.  

We recall that in \cite{SGOP20} (see also \cite{GiSa14}), 
by using the structure 
\begin{equation}\label{EQ:STRUCT}
\Box|\bigcirc = \Box\wedge \bigcirc +\prev(\Box|\bigcirc)\no{\bigcirc},
\end{equation} 
with $\Box=B|K$, $\bigcirc=A|H\neq \emptyset$, and $\Box\wedge \bigcirc=(B|K)\wedge_{gs}(A|H)$, the iterated conditional $	(B|K)|_{gs}(A|H)$ has been defined as the following conditional random quantity
\begin{equation}\label{EQ:itergs}
	(B|K)|_{gs}(A|H)=(A|H)\wedge_{gs} (B|K)+\mu_{gs}\, (\no{A}|H),
\end{equation}
where $\mu_{gs}=\prev[(B|K)|_{gs}(A|H)]$. We underline that when $\Box=A, \bigcirc =H$ formula (\ref{EQ:STRUCT}) reduces to formula (\ref{EQ:AgH}). 
We also  recall that (\cite{GiSa14})
\begin{equation}\label{EQ:CPTgs}
    \prev[(A|H) \wedge_{gs} (B|K)]
=\prev[(B|K)|_{gs}(A|H)]P(A|H).
\end{equation} 
Moreover, 
in \cite[Definition 7]{GiSa21ECSQARU}, by exploiting the structure (\ref{EQ:STRUCT}),  
the iterated conditioning $|_{gs}$   has been extended to the case of conjoined conditionals. More precisely, denoting by $\C(\F_1)$ and $\C(\F_2)$  the $\wedge_{gs}$-conjunctions  of the conditional events in two finite families $\F_1$ and $\F_2$ (\cite{GiSa19}), the iterated conditional
$\C(\F_2)|_{gs}\C(\F_1)$ has been defined as
\begin{equation}
\label{EQ:GENITERCOND}
\C(\F_2)|_{gs}\C(\F_1) =	\C(\F_2)\wedge_{gs}\C(\F_1) + \mu (1-\C(\F_1)) 
		=\C(\F_1 \cup \F_2)+\mu (1-\C(\F_1)),
\end{equation}
where $\mu = \prev[\C(\F_2)|_{gs}\C(\F_1)]$. In addition, it holds \cite[Equation (8)]{GiSa21ECSQARU}
that 
\begin{equation}\label{EQ:PREVCOMP}
\prev[\C(\F_2)\wedge_{gs} \C(\F_1)]=\prev[\C(\F_2)|_{gs}\C(\F_1)]\prev[\C(\F_1)].
\end{equation} Formulas (\ref{EQ:CPTgs} and \ref{EQ:PREVCOMP}) generalize the  compound probability theorem recalled in (\ref{EQ:CPT}).

We now introduce the  different notions of iterated conditioning, beyond $|_{gs}$,  obtained with the structure  (\ref{EQ:STRUCT}), by using the    trivalent logic conjunctions  $\wedge_K, \wedge _L, \wedge_B$, and  $\wedge_S$, recalled in Section \ref{SEC:2}.
Among other things, we will show that a formula like (\ref{EQ:CPTgs}) still holds for each of these new objects.

\begin{definition}\label{DEF:ITECMPTHM}
Given two conditional events $A|H$, $B|K$, with $AH\neq \emptyset$, for each $i\in\{K,L,B,S,gs\}$, we define the iterated conditional $(B|K)|_i(A|H)$ as
\begin{equation}\label{EQ:ITECMPTHM}
(B|K)|_i(A|H)=(A|H)\wedge_i (B|K)+\mu_i \,(\no{A}|H),
\end{equation}
where $\mu_i=\prev[(B|K)|_i(A|H)]$.
\end{definition}

Now,  we  check  the validity of properties P1--P4, introduced in Section \ref{SEC:3}, for each iterated conditioning $|_i$, $i\in\{K,L,B,S,gs\}$.
\subsection{Property P1}  In the next result we show  that   each iterated conditioning  $|_{i}$, $i\in\{K,L,B,S\}$, satisfies property P1.
\begin{theorem}
Given two conditional events $A|H$, $B|K$, with $AH\neq \emptyset$, it holds that
\begin{equation}\label{EQ:CANITER}
   ((A|H)\wedge_i(B|K))|_i(A|H)=(B|K)|_i(A|H),\;\;  i\in\{K,L,B,S\}.
\end{equation}
\end{theorem}
\begin{proof}
Let be given  $i\in\{K,L,B,S\}$. 
We set $\prev[(B|K)|_i(A|H)]=\mu_i$
 and $\prev[((A|H)\wedge(B|K))|_i(A|H)]=\nu_i$.
By Definition \ref{DEF:ITECMPTHM}, 
as $(A|H)\wedge_i(A|H)\wedge_i (B|K)=(A|H)\wedge_i (B|K)$, 
it holds that 
\begin{equation}\label{EQ:ITECMPTHM2}
((A|H)\wedge_i(B|K))|_i(A|H)=(A|H)\wedge_i (B|K)+\nu_i (\no{A}|H).
\end{equation}
Then, as $(B|K)|_i(A|H)=(A|H)\wedge_i (B|K)+\mu_i (\no{A}|H)$, in order to prove (\ref{EQ:CANITER}) it is enough to verify that $\nu_i=\mu_i$.
We observe that 
$
((A|H)\wedge_i(B|K))|_i(A|H)-
(B|K)|_i(A|H)=
(\nu_i-\mu_i)(\no{A}|H),
$
where $\nu_i-\mu_i=
\prev[((A|H)\wedge_i(B|K))|_i(A|H)-
(B|K)|_i(A|H)].$
By setting $P(A|H)=x$, it holds that 
\[
(\nu_i-\mu_i)(\no{A}|H)=
\left\{
\begin{array}{ll}
0, \mbox{ if } A|H=1,\\
\nu_i-\mu_i, \mbox{ if } A|H=0,\\
(\nu_i-\mu_i)(1-x), \mbox{ if } A|H=x,\; 0<x<1.
\end{array}
\right.
\]
Notice that, in the betting scheme, 
 $\nu_i-\mu_i$ is the amount to be paid in order to receive the random amount $(\nu_i-\mu_i)(\no{A}|H)$. Then, by coherence, $\nu_i-\mu_i$ must be a linear convex combination of the possible values of $(\nu_i-\mu_i)(\no{A}|H)$, by discarding the cases  where
the bet is called off, that is the cases where you  receive back the paid amount 
$\nu_i-\mu_i$, whatever $\nu_i-\mu_i$ be. In other words, coherence requires that $\nu_i-\mu_i$ must belong to the convex hull of the set $\{0,(\nu_i-\mu_i)(1-x)\}$, that is $\nu_i-\mu_i=\alpha \cdot 0+(1-\alpha)(\nu_i-\mu_i) (1-x)$,  for some $\alpha\in[0,1]$. Then, as $0<x<1$, we observe that the previous equality  holds if and only if $\nu_i-\mu_i=0$, that is $\nu_i=\mu_i$. Therefore, equality (\ref{EQ:CANITER}) holds.
\end{proof}
We recall that  
$|_{gs}$ satisfies property P1. Indeed, for the generalized iterated conditional (\ref{EQ:GENITERCOND}) 
it holds that  (\cite[Theorem 5]{GiSa21ECSQARU}) 
\begin{equation}\label{EQ:P1gs}
\mathscr{C}(\F_2)|_{gs}\mathscr{C}(\F_1)=\mathscr{C}(\F_1\cup \F_2)|_{gs}\mathscr{C}(\F_1).
\end{equation}
Then, by applying (\ref{EQ:P1gs}) with $\F_1=\{A|H\}$, $\F_2=\{B|K\}$, it follows that  
\[(B|K)|_{gs}(A|H)=((A|H)\wedge_{gs}(B|K))|_{gs}(A|H).\]
Thus, property P1 is satisfied by each  iterated conditioning $|_i\in\{|_K,|_L,|_B,|_S,|_{gs}\}$.

\subsection{Property P2}
  We show  that   each iterated conditional  $|_{i}$, $i\in\{K,L,B,S,gs\}$, satisfies property P2. 
  We observe that,  
  coherence requires $\mu_i\geq 0$, $i\in\{K,L,B,S,gs\}$ (see Remark \ref{REM:mupositive} in Section \ref{SEC:GENITER:P3}). 
  Then, for each $i\in\{K,L,B,S,gs\}$, 
  as $(A|H)\wedge_i (B|K)\leq 
(A|H)\wedge_i (B|K)+\mu_i(\no{A}|H)$, 
  from  
 (\ref{EQ:ITECMPTHM}) it follows that  
\begin{equation}
(A|H)\wedge_i (B|K)\leq 
(B|K)|_i(A|H),
\end{equation}
and hence 
$
\prev[(A|H)\wedge_i (B|K)]\leq \prev[(B|K)|_i(A|H)].
$
\subsection{Property P3}
\label{SEC:GENITER:P3}
We already recalled in equation (\ref{EQ:CPTgs}) that the iterated conditional $|_{gs}$  satisfies P3. 
Then, by setting 
$z_{gs}=\prev[(A|H) \wedge_{gs} (B|K)]$, $\mu_{gs}=\prev[(B|K)|_{gs}(A|H)]$, and $x=P(A|H)$ it holds that $z_{gs}=x \mu_{gs}$.
By exploiting the  structure (\ref{EQ:STRUCT}), 
we show below that P3 is also valid
for $|_K,|_L,|_B,|_S$.
 Indeed, by the  linearity property of a coherent prevision,  from (\ref{EQ:ITECMPTHM}) we obtain that
\begin{equation}
\begin{array}{lll}
\mu_i &=& \prev[(B|K)|_i(A|H)] = 
P[(B|K) \wedge_i (A|H)]+\mu_iP(\no{A}|H)
=\\&=&P[(B|K) \wedge_i (A|H)] + \mu_i\, P(\no{A}|H) = z_i + \mu_i(1-x) \,,
\end{array}
\label{prevision}
\end{equation}
where $x=P(A|H)$ and $z_i=P[(A|H) \wedge_i (B|K)]$, $i\in\{K,L,B,S\}$. As $\mu_i=z_i + \mu_i(1-x)$, it follows that 
\[ z_i=\mu_i x, \;\; i\in\{K,L,B,S\}.
\] Therefore,  coherence requires that 
\begin{equation}\label{EQ:COMPTHM}
\prev[(A|H) \wedge_i (B|K)]
=\prev[(B|K)|_i(A|H)]P(A|H),\; i\in\{K,L,B,S,gs\},
\end{equation}
which states that the compound prevision formula for iterated conditional (property P3) is valid for  each  iterated conditioning $|_i,\  i\in\{K,L,B,S,gs\}$. 
\begin{remark}\label{REM:mupositive}
    Notice that, as $P(A|H)\geq 0$ and $\prev[(A|H) \wedge_i (B|K)]\geq 0$ for $i\in\{K,L,B,S,gs\}$, from 
    (\ref{EQ:COMPTHM}) it follows that $\prev[(B|K)|_i(A|H)]\geq 0$, $i\in\{K,L,B,S,gs\}$. 
\end{remark}

We set $P(A|H)=x$, $P(B|K)=y$, 
$P[(A|H) \wedge_{i} (B|K)]=z_i$, 
$i\in\{K,L,B,S\}$, $\prev[(A|H) \wedge_{i} (B|K)]=z_{gs}$,  
$P[(B|K) |_{i} (A|H)]=\mu_i$,
$i\in\{K,L,B,S, gs\}$. Then,  based on  Definition \ref{DEF:ITECMPTHM} and on the compound prevision theorem, for each $|_i$, $i\in\{K,L,B,S, gs\}$,  we obtain the  random quantities 
 $(B|K)|_{i}(A|H)$ illustrated below.
\begin{itemize}
    \item  $(|_K)$.  We recall that    $(A|H) \wedge_{K} (B|K)=AHBK|(HK\vee \overline{A}H\vee \overline{B}K)$, then we have
\begin{equation} \label{EQ:iteK}
\begin{array}{l}
(B|K)|_{K}(A|H) = AHBK|(HK\vee \overline{A}H\vee \overline{B}K) + \mu_K (\bar{A}|H)
   =\\= \begin{cases}
    1,               & {\ \text{if} \ AHBK \ \text{is true}}, \\
    0,               & {\ \text{if} \ AH\no{B}K \ \text{is true}},\\
    z_K,               & {\ \text{if} \ AH\no{K} \ \text{is true}},\\
    z_K +\mu_K (1-x),    & {\ \text{if} \ \no{H}BK\vee \no{H}\no{K} \ \text{is true}},\\
    \mu_K (1-x),       & {\ \text{if} \ \no{H}\no{B}K \ \text{is true}},\\
    \mu_K,             & {\ \text{if} \ \no{A}H \ \text{is \ true}}.
    \end{cases}
\end{array}
\end{equation}
  From (\ref{EQ:COMPTHM}),  coherence requires that $z_K=x\mu_K$ and $z_K+\mu_K(1-x)=\mu_K$. Then, we obtain 
\begin{equation} \label{EQ:iteK'}
\begin{array}{l}
(B|K)|_{K}(A|H) 
   = \begin{cases}
    1,               & {\ \text{if} \ AHBK \ \text{is true}}, \\
    0,               & {\ \text{if} \ AH\no{B}K \ \text{is true}},\\
    x\mu_K,               & {\ \text{if} \ AH\no{K} \ \text{is true}},\\
    \mu_K (1-x),       & {\ \text{if} \ \no{H}\no{B}K \ \text{is true}},\\
    \mu_K,             & {\ \text{if} \ \no{A}H \vee \no{H}BK\vee \no{H}\no{K} \ \text{is \ true}}.
    \end{cases}
\end{array}
\end{equation}
\item $(|_L).$ We recall that    $(A|H) \wedge_{L}(B|K)=AHBK|( AHBK \vee \overline{A} \, H\vee \, \overline{B} K \vee \overline{H} \, \overline{K})$, then we have
\begin{equation}\label{EQ:iteL}
\begin{array}{l}
(B|K)|_{L}(A|H) = AHBK|( AHBK \vee \overline{A} \, H\vee \, \overline{B} K \vee \overline{H} \, \overline{K}) + \mu_L (\bar{A}|H)
=\\= \begin{cases}
    1,               & {\ \text{if} \ AHBK \ \text{is true}},\\
    0,               & {\ \text{if} \ AH\no{B}K \ \text{is true}},\\
    z_L,              & {\ \text{if} \ AH\no{K} \ \text{is true}},\\
    z_L +\mu_L (1-x),    & {\ \text{if} \ \no{H}BK \ \text{is true}},\\
    \mu_L (1-x),       & {\ \text{if} \ \no{H}\no{B}K\lor \no{H}\no{K}     \ \text{is true}},\\
    \mu_L,             & {\ \text{if} \ \no{A}H \ \text{is true}}.
    \end{cases}
\end{array}
\end{equation}
  From (\ref{EQ:COMPTHM}),  coherence requires that $z_L=x\mu_L$ and $z_L+\mu_L(1-x)=\mu_L$. Then, we obtain 
\begin{equation}\label{EQ:iteL'}
\begin{array}{l}
(B|K)|_{L}(A|H) 
= \begin{cases}
    1,               & {\ \text{if} \ AHBK \ \text{is true}},\\
    0,               & {\ \text{if} \ AH\no{B}K \ \text{is true}},\\
    x\mu_L,              & {\ \text{if} \ AH\no{K} \ \text{is true}},\\
    \mu_L (1-x),       & {\ \text{if} \ \no{H}\no{B}K\lor \no{H}\no{K}     \ \text{is true}},\\
    \mu_L,             & {\ \text{if} \ \no{A}H\vee \no{H}BK \ \text{is true}}.
    \end{cases}
\end{array}
\end{equation}
\item $(|_B).$ We recall that    $(A|H) \wedge_{B}(B|K)= AHBK|BK$, then we have
\begin{equation}\label{EQ:iteB}
\begin{array}{l}
(B|K)|_{B}(A|H) = AHBK|HK + \mu_B (\bar{A}|H) 
=\\=\begin{cases}
    1,               & {\ \text{if} \ AHBK \ \text{is true}},\\
    0,               & {\ \text{if} \ AH\no{B}K \ \text{is true}},\\
    z_B,               & {\ \text{if} \ AH\no{K} \ \text{is true}},\\
    z_B+\mu_B,           & {\ \text{if} \ \no{A}H\no{K} \ \text{is true}},\\
    z_B +\mu_B (1-x),    & {\ \text{if}\ \no{H}\ \text{is true}},\\
    \mu_B,             & {\ \text{if} \ \no{A}HBK \lor \no{A}H\no{B}K \ \text{is true}}.
    \end{cases}
\end{array}
\end{equation}
  From (\ref{EQ:COMPTHM}),  coherence requires that $z_B=x\mu_B$ and $z_B+\mu_B(1-x)=\mu_B$. Then, we obtain
\begin{equation}\label{EQ:iteB'}
\begin{array}{l}
(B|K)|_{B}(A|H) 
=\begin{cases}
    1,               & {\ \text{if} \ AHBK \ \text{is true}},\\
    0,               & {\ \text{if} \ AH\no{B}K \ \text{is true}},\\
    x\mu_B,               & {\ \text{if} \ AH\no{K} \ \text{is true}},\\
    \mu_B(1+x),           & {\ \text{if} \ \no{A}H\no{K} \ \text{is true}},\\
    \mu_B,             & {\ \text{if} \ \no{H}\vee \no{A}HBK \lor \no{A}H\no{B}K \ \text{is true}}.
    \end{cases}
\end{array}
\end{equation}
\item $(|_S).$ We recall that    $ (A|H)\wedge_S (B|K)=((AH\vee\overline{H})\wedge (BK\vee \overline{K}))$, then we have
\begin{equation}\label{EQ:iteS}
\begin{array}{l}
(B|K)|_S(A|H) = ((AH\vee\overline{H})\wedge (BK\vee \overline{K}))|(H\vee K) + \mu_S (\no{A}|H)
   = \\=\begin{cases}
    1,               & {\ \text{if} \ AHB \lor AH\no{K} \ \text{is true}},\\
    0,               & {\ \text{if} \ AH\no{B}K \ \text{is true}},\\
    1 +\mu_S (1-x),    & {\ \text{if} \ \no{H}BK \ \text{is true}},\\
    \mu_S (1-x),       & {\ \text{if} \ \no{H}\no{B}K \ \text{is true}},\\
    z_S +\mu_S (1-x),    & {\ \text{if}\ \no{H}\no{K}\ \text{is true}},\\
    \mu_S             & {\ \text{if} \ \no{A}H \ \text{is true}}.
    \end{cases}
\end{array}
\end{equation}
 From (\ref{EQ:COMPTHM}),  coherence requires that $z_S=x\mu_S$ and $z_S+\mu_S(1-x)=\mu_S$. Then, we obtain
 \begin{equation}\label{EQ:iteS'}
\begin{array}{l}
(B|K)|_S(A|H) 
   = \begin{cases}
    1,               & {\ \text{if} \ AHB \lor AH\no{K} \ \text{is true}},\\
    0,               & {\ \text{if} \ AH\no{B}K \ \text{is true}},\\
    1 +\mu_S (1-x),    & {\ \text{if} \ \no{H}BK \ \text{is true}},\\
    \mu_S (1-x),       & {\ \text{if} \ \no{H}\no{B}K \ \text{is true}},\\
    \mu_S             & {\ \text{if} \ \no{A}H\vee \no{H}\no{K} \ \text{is true}}.
    \end{cases}
\end{array}
\end{equation}
\item $(|_{gs}).$ We recall that    $(A|H) \wedge_{gs} (B|K)=(AH\vee BK+x\no{H}\no{B}K+y\no{A}H\no{K})|(H\vee K)$, then we have (\cite{GiSa14})
\begin{equation}\label{EQ:itegs}
\begin{array}{l}
(B|K)|_{gs}(A|H) = (AHBK+x\no{H}BK+yAH\no{K})|(H\vee K) + \mu_{gs} (\no{A}|H)
   = \\=\begin{cases}
    1,               & {\ \text{if} \ AHBK \ \text{is true}},\\
    0,               & {\ \text{if} \ AH\no{B}K \ \text{is true}},\\
    y,    & {\ \text{if} \ AH\no{K} \ \text{is true}},\\
    x +\mu_{gs} (1-x),    & {\ \text{if} \ \no{H}BK \ \text{is true}},\\
    \mu_{gs} (1-x),       & {\ \text{if} \ \no{H}\no{B}K \ \text{is true}},\\
    z_{gs} +\mu_{gs} (1-x),    & {\ \text{if}\ \no{H}\,\no{K}\ \text{is true}},\\
    \mu_{gs},             & {\ \text{if} \ \no{A}H \ \text{is true}}.
    \end{cases}
\end{array}
\end{equation}
From (\ref{EQ:COMPTHM}),  coherence requires that $z_{gs}=x\mu_{gs}$ and $z_{gs}+\mu_{gs}(1-x)=\mu_{gs}$. Then, we recall that 
\begin{equation}\label{EQ:itegs'}
\begin{array}{l}
(B|K)|_{gs}(A|H) 
   = \begin{cases}
    1,               & {\ \text{if} \ AHBK \ \text{is true}},\\
    0,               & {\ \text{if} \ AH\no{B}K \ \text{is true}},\\
    y,    & {\ \text{if} \ AH\no{K} \ \text{is true}},\\
    x +\mu_{gs} (1-x),    & {\ \text{if} \ \no{H}BK \ \text{is true}},\\
    \mu_{gs} (1-x),       & {\ \text{if} \ \no{H}\no{B}K \ \text{is true}},\\
    \mu_{gs},             & {\ \text{if} \ \no{A}H \vee \no{H}\,\no{K}\ \text{is true}}.
    \end{cases}
\end{array}
\end{equation}
\end{itemize} 
In Table \ref{TAB:ITECMPPRB} we summarize the possible values of  the   iterated conditionals $(B|K)|_i(A|H)$, $i\in\{K,L,B,S,gs\}$.
\begin{table}[h]
	\begin{small}
	\def\True{$1$}
	\def\False{$0$}
	\def\Void{$z$}
	\def\Voidx{$x$}
	\def\Voidy{$y$}
	\centering
	\begin{tabular}{l|c|c|c|c|c}
		& $(B|K)|_{K}(A|H)$ & $(B|K)|_{L}(A|H)$ & $(B|K)|_{B}(A|H)$ &  $(B|K)|_{S}(A|H)$  &   $(B|K)|_{gs}(A|H)$   \\ \hline
		$AHBK  $                           &       \True       &       \True       &       \True       &        \True        &         \True          \\
		$	{AH}{\overline{B}K}$             &      \False       &      \False       &      \False       &       \False        &         \False         \\
		$	AH\overline{K}  $                &     $x\mu_K$     &     $x\mu_L$     &     $x\mu_B$     &        \True        &         \Voidy         \\
		$	{\overline{A}H}BK$             &      $\mu_K$      &      $\mu_L$      &      $\mu_B$      &       $\mu_S$       &       $\mu_{gs}$       \\
  		$	{\overline{A}H}\overline{B}K$             &      $\mu_K$      &      $\mu_L$      &      $\mu_B$      &       $\mu_S$       &       $\mu_{gs}$       \\
		$	{\overline{A}H} \, \overline{K}$ &      $\mu_K$      &      $\mu_L$      &    $\mu_B(1+x)$    &       $\mu_S$       &       $\mu_{gs}$       \\
		$	\overline{H}{BK}$                & $\mu_K$  & $\mu_L$  & $\mu_B$  &   1+ $\mu_S(1-x)$   &   $x+\mu_{gs}(1-x)$    \\
		$	\overline{H}\,{\overline{B}K}$   &   $\mu_K(1-x)$    &   $\mu_L(1-x)$    & $\mu_B$  &    $\mu_S(1-x)$     &    $\mu_{gs}(1-x)$     \\
		$	\overline{H}\,\overline{K} $     & $\mu_K$  &   $\mu_L(1-x)$    & $\mu_B$  & $\mu_S$ & $\mu_{gs}$ 
	\end{tabular}
	\caption{Numerical values   of  $(B|K)|_{i}(A|H)$, $i\in\{K,L,B,S,gs\}$. 
 We denotes $P(A|H)=x$, $P(B|K)=y$, 
	and $\prev[(B|K)|_i(A|H)]=\mu_i$, $i\in\{K,L,B,S,gs\}$. 
	}
	\label{TAB:ITECMPPRB}
\end{small}
\end{table} 
\subsection{Property P4} 
In this section, for each $i\in\{K,L,B,S\}$,  we  find the set  of all the  coherent assessments on the family $\{A|H, B|K, (A|H)\wedge_{i}(B|K), (B|K)|_{i}(A|H)\}$. 
Moreover, for each $i\in\{K,L,B,S\}$, 
we determine  the interval of  coherent extensions   on
 $(B|K)|_{i}(A|H)$ and we check the validity of property P4. Then we recall that property P4 is satisfied by $|_{gs}$.

\paragraph{The iterated conditioning $|_K$}  
\begin{theorem}\label{THM:PIK}
Let $A$, $B$, $H$, $K$ be any logically independent events. The set $\Pi$ of all the coherent assessment $(x,y,z,\mu)$ on the family $\F=\{A|H, B|K, (A|H)\wedge_{K}(B|K), (B|K)|_{K}(A|H)\}$ is $\Pi=\Pi' \cup \Pi''$, where 
$\Pi'=\{(x,y,z,\mu): x\in(0,1], y\in [0,1], z\in[z', z''], \mu=\frac{z}{x}\}$
 with $ z'=0$,    $z''=min\{x,y\}$, and
$\Pi''=\{(0,y,0,\mu): (y,\mu)\in[0,1]^2\}.$
\end{theorem}
 \begin{proof}
 It is well-known that the assessment $(x,y)$ on $\{A|H, B|K\}$ is coherent for every $(x,y)\in [0,1]^2$. 
 By Table (\ref{tab:LUB}), the assessment $z=P((A|H)\wedge_{K}(B|K))$ is a coherent extension of $(x,y)$ if and only if $z \in [z',z'']$ where $z'=0$ and $z''=min\{x,y\}$. Assuming $x>0$, from (\ref{EQ:COMPTHM}) it holds that $\mu=\frac{z}{x}$. Then, every $(x,y,z,\mu)\in \Pi'$ is coherent, i.e., $\Pi'\subseteq \Pi$. Of course, if $x>0$ and $(x,y,z,\mu)\notin \Pi'$, then the assessment $(x,y,z,\mu)$ is not coherent, i.e. $(x,y,z,\mu)\notin \Pi$.\\
Let us consider now the case $x=0$, so that $z'=0$ and $z''=0$. We show that the assessment $(0,y,0,\mu)$ is coherent if and only if $(y,\mu)\in [0,1]^2$, that is $\Pi''\subseteq \Pi$. 
As $x=0$, from (\ref{EQ:iteK'}), it holds that 
\begin{equation} \label{EQ:iteK'x=0}
\begin{array}{l}
(B|K)|_{K}(A|H) 
   = \begin{cases}
    1,               & {\ \text{if} \ AHBK \ \text{is true}}, \\
    0,               & {\ \text{if} \ AH\no{B}K\vee AH\no{K} \ \text{is true}},\\
    \mu,             & {\ \text{if} \ \no{A}H \vee \no{H}BK\vee \no{H}\no{B}K \vee \no{H}\no{K} \ \text{is \ true}},
    \end{cases}
\end{array}
\end{equation}
that is $(B|K)|_{K}(A|H)=AHBK+\mu(\no{A}\vee \no{H})$. Moreover we have that $AHBK+\mu(\no{A}\vee \no{H})$ and the conditional event  $BK|AH=AHBK+\eta(\no{A}\vee \no{H})$, where $\eta=P(BK|AH)$, coincide when $AH$ is true. Then, from Theorem \ref{THM:EQ-CRQ}, it follows that 
$\mu=\eta$ and hence
$(B|K)|_{K}(A|H)$  and $BK|AH$ should also coincide when $AH$ is false. Thus,  $(B|K)|_{K}(A|H)$ and $BK|AH$ coincide in all cases. Therefore $\F=\{A|H, B|K, (A|H)\wedge_{K}(B|K), (B|K)|_{K}(A|H)\}=\{A|H, B|K, AHBK|(AHBK\vee \overline{A}H\vee \overline{B}K), BK|AH\}$.
The constituents $C_h$'s and the points $Q_h$'s associated with $(\F, \P)$, where $\P=(0,y,0,\mu)$ are the following:
\[
\begin{array}{l}
C_1=AHBK, C_2=AH\no{B}K, C_3=AH\no{K}, C_4=\no{A}HBK, C_5=\no{A}H\no{B}K,\\ C_6=\no{A}H\no{K}, C_7=\no{H}BK, C_8=\no{H}\no{B}K, C_0=\no{H}\,\no{K},
\end{array}
\]
and
\[
\begin{array}{l}
Q_1=(1,1,1,1), Q_2=(1,0,0,0), Q_3=(1,y,0,0), Q_4=(0,1,0,\mu), 
Q_5=(0,0,0,\mu),\\ Q_6=(0,y,0,\mu),Q_7=(0,1,0,\mu), Q_8=(0,0,0,\mu), 
Q_0=(0,y,0,\mu).
\end{array}
\]
The constituents contained in $\H_4=H\vee K$ are $C_1,\ldots,C_8$.
The system $(\Sigma)$ in (\ref{SYST-SIGMA}) associated with the pair $(\F,\P)$  is\\
\begin{equation}\label{EQ:SIGMAK}
\left\{
    \begin{array}{llll}
     \lambda_1+\lambda_2+\lambda_3=0, \\
     \lambda_1+y\lambda_3+\lambda_4+y\lambda_6+\lambda_7=y,\\
     \lambda_1=0,\\
     \lambda_1+\mu \lambda_4+\mu \lambda_5+\mu \lambda_6+\mu \lambda_7+\mu \lambda_8=\mu\\
    \lambda_1+ \lambda_2+\lambda_3+\lambda_4+\lambda_5+\lambda_6+\lambda_7+ \lambda_8=1, \\ \lambda_i\geq0 \ \forall i=1, \dots, 8.\\
    \end{array}
\right.    
\end{equation}
We observe that, as $y\in[0,1]$,  $\P$ belongs to the segment with end points $Q_4$, $Q_5$ because $\P=yQ_4+(1-y)Q_5$. Then,  $\Lambda=(0,0,0,y,1-y,0,0,0)$ is a solution of (\ref{EQ:SIGMAK}) with
\[
\phi_{1}(\Lambda)=\sum_{h:C_h\subseteq H} \lambda_h=\lambda_1+\lambda_2+\lambda_3+\lambda_4+\lambda_5+\lambda_6=1>0,
\]
\[
\phi_{2}(\Lambda)=\sum_{h:C_h\subseteq K}\lambda_h=\lambda_1+\lambda_2+\lambda_4+\lambda_5+\lambda_7+\lambda_8=1>0,
\]
\[
\phi_{3}(\Lambda)=\sum_{h:C_h\subseteq (AHBK\vee \overline{A}H\vee \overline{B}K)}\lambda_h=\lambda_1+\lambda_2+\lambda_4+\lambda_5+\lambda_6+\lambda_8=1>0,
\]
\[
\phi_{4}(\Lambda)=\sum_{h:C_h\subseteq AH}\lambda_h=\lambda_1+\lambda_2+\lambda_3=0.
\]
Let $\mathcal{S}'=\{(0,0,0,y,1-y,0,0,0)\}$  denote a subset of the set $\mathcal{S}$ of all solutions of (\ref{EQ:SIGMAK}). We have that $M_1'=M_2'=M_3'=1$ and $M_4'=0$   (as defined in (\ref{EQ:I0'})). 
It follows that $I'_0=\{4\}$. As the sub-assessment $\P'_0=\mu$ on $\F'_0=\{BK|AH\}$ is coherent $\forall \mu \in [0,1]$,  by Theorem \ref{CNES-PREV-I_0-INT}, it follows that 
the assessment  $(0,y,0,\mu)$ on $\F=\{A|H, B|K, (A|H)\wedge_{K}(B|K), (B|K)|_{K}(A|H)\}$ is coherent for every $(y,\mu)\in[0,1]^2$, that is 
$(0,y,0,\mu)\in \Pi''$. Thus $\Pi''\subseteq \Pi$.
Of course, if $(0,y,z,\mu)\notin\Pi''$ the assessment $(0,y,z,\mu)$ is not coherent and hence $(0,y,z,\mu)\notin\Pi$.
Therefore
$\Pi=\Pi'\cup \Pi'' $.
 \end{proof}
Based on Theorem~\ref{THM:PIK}, we obtain  
\begin{theorem}\label{THM:LUK}
	Let $A,B,H,K$ be any logically independent events. Given any assessment $(x,y)\in[0,1]^2$ on $\{A|H,B|K\}$, for the iterated conditional $(B|K)|_K(A|H)$ the  extension $\mu_K=\prev((B|K)|_K(A|H))$ is coherent if and only if 
	$\mu_K\in [\mu_K', \mu_K'']$, where	\begin{equation}\label{EQ:P4K}
	\begin{array}{ll}
	\mu_K'=0,\;\;
	\mu_K''=\left\{
	\begin{array}{ll}
	\min\left\{1,\tfrac yx\right\},& \text{ if } x>0;\\
	1,& \text{ if } x=0.\\
	\end{array}
	\right.
	\end{array}
	\end{equation} 
\end{theorem}
\begin{proof} Of course, $(x,y)$ on $\{A|H,B|K\}$ is coherent.
Assume that $x>0$. We simply write $\mu$ instead of $\mu_K$. From Theorem \ref{THM:PIK} it follows that the set of all coherent assessments $(x,y,z,\mu)$ on $\F=\{A|H, B|K, (A|H)\wedge_{K}(B|K), (B|K)|_{K}(A|H)\}$ is $\Pi'=\{(x,y,z,\mu): 0<x\leq 1, 0\leq y\leq 1, z' \leq z\leq z'', \mu=\frac{z}{x}\} $, where $z'=0$ and $z''=\min\{x,y\}$ (\ref{tab:LUB}). Then,  $\mu$ is a coherent extension of $(x,y)$ if and only if $\mu \in [\mu', \mu'']$, where $\mu'=0$ and $\mu''=\frac{z''}{x}=\min\left\{1,\tfrac yx\right\}$.
Assume that  $x=0$. From Theorem \ref{THM:PIK} it follows that the set of all coherent assessments $(x,y,z,\mu)$ on $\F=\{A|H,B|K$, $(A|H)\wedge_K (B|K), (B|K)|_K(A|H)\}$ is  $\Pi''=\{(0,y,0,\mu): (y,\mu)\in[0,1]^2\}$. Then, 
$\mu$ is a coherent extension of $(x,y)$ if and only if $\mu \in[\mu',\mu'']$, where $\mu'=0$ and $\mu''=1$.

\end{proof}
\begin{remark}\label{REM:ITKP4}
We notice that the lower and upper bounds $\mu'_K$ and  $\mu''_K$ for $(B|K)|_K(A|H)$, given in 
(\ref{EQ:P4K}, do not coincide with  the lower and upper bound $\mu'$ and $\mu''$ for $B|A$, given in (\ref{EQ:LU}).  In particular, the extension $\mu_K=0$ to  $(B|K)|_{K}(A|H)$ of $(1,1)$ on $\{A|H, B|K\}$ is coherent because, as $\mu'_K=0$, the assessment  $(1,1,0)$ on $\{A|H, B|K, (B|K)|_{K}(A|H)\}$ is coherent. However, $\mu=0$ is not a coherent extension to  $B|A$ of $(1,1)$ on $\{A, B\}$, because, as $\mu'=\frac{\max\{1+1-1,0\}}{1}=1$, the assessment $(1,1,0)$ on $\{A, B, B|A\}$ is not coherent. 
Therefore, property P4 is not satisfied by $|_K$.
\end{remark}
\paragraph{The iterated conditioning $|_L$}
We obtain results similar to the case $|_K$. Indeed we have
\begin{theorem}\label{THM:PIL}
Let $A$, $B$, $H$, $K$ be any logically independent events. The set $\Pi$ of all the coherent assessment $(x,y,z,\mu)$ on the family $\F=\{A|H, B|K, (A|H)\wedge_{L}(B|K), (B|K)|_{L}(A|H)\}$ is $\Pi=\Pi' \cup \Pi''$, where 
$\Pi'=\{(x,y,z,\mu): x\in(0,1], y\in [0,1], z\in[z', z''], \mu=\frac{z}{x}\}$
 with $ z'=0$,    $z''=min\{x,y\}$, and
$\Pi''=\{(0,y,0,\mu): (y,\mu)\in[0,1]^2\}.$
\end{theorem}
\begin{proof}
See \ref{SEC:APPL}.
\end{proof}

Based on Theorem~\ref{THM:PIL}, we obtain  
\begin{theorem}\label{THM:LU}
	Let $A,B,H,K$ be any logically independent events. Given a coherent assessment $(x,y)$ on $\{A|H,B|K\}$, for the iterated conditional $(B|K)|_L(A|H)$ the  extension $\mu_L=\prev((B|K)|_L(A|H))$ is coherent if and only if 
	$\mu_L \in [\mu_L', \mu_L'']$, where	\begin{equation}\label{EQ:P4L}
	\begin{array}{ll}
	\mu_L'=0,\;\;
	\mu_L''=\left\{
	\begin{array}{ll}
	\min\left\{1,\tfrac yx\right\},& \text{ if } x>0;\\
	1,& \text{ if } x=0.\\
	\end{array}
	\right.
	\end{array}
	\end{equation} 
\end{theorem}
\begin{proof}
The proof is the same as in Theorem \ref{THM:LUK} where $K$ is replaced by $L$ and Theorem \ref{THM:PIK} is replaced by Theorem \ref{THM:PIL}.
\end{proof}
\begin{remark}
Theorem \ref{THM:LU} shows that the interval of coherent extensions $[\mu'_L, \mu''_L]$ for $\mu_L=\prev((B|K)|_L(A|H))$ does not coincide with the interval $[\mu', \mu'']$ of coherent assessment for $\mu=P(B|A)$ where $\mu'$ and $\mu''$ are as in (\ref{EQ:LU}). Thus, property P4 is not satisfied by $|_L$.
\end{remark}
 \paragraph{The iterated conditioning  $|_B$}
 We continue by analyzing the set of coherent extensions for the iterated conditional $|_B$. 
We first  show that, when we  evaluate $\prev[(B|K)|_B(A|H)]$,  
if $P(A|H)=0$, then 
coherence requires that $P[(A|H)\wedge_B(B|K)]=0$. 

\begin{remark}\label{coherenceB}
We recall that, given  any $(x,y)\in[0,1]^2$, where $x=P(A|H)$ and $y=P(B|K)$,  the interval of coherent extensions on $z_B=P[(A|H)\wedge_{B}(B|K)]$ is $[z_B', z_B'']=[0,1]$ (Table \ref{tab:LUB}). In particular it is coherent  to assess $(0,y,z_B)$, with $z_B>0$, on $\{A|H,B|K,(A|H)\wedge_B(B|K)\}$. However, for the object 
$(B|K)|_B(A|H)$ coherence  also requires that $z_B=\mu_Bx$ (see \ref{prevision}). Then, 
coherence requires that $z_B=0$ when we consider the assessment $(0,y,z_B,\mu_B)$ on
$\{A|H, B|K,(A|H)\wedge_{B}(B|K),(B|K)|_{B}(A|H)\}$.
\end{remark} 
 
\begin{theorem}\label{THM:PIB}
Let $A$, $B$, $H$, $K$ be any logically independent events. The set $\Pi$ of all the coherent assessment $(x,y,z,\mu)$ on the family $\F=\{A|H, B|K, (A|H)\wedge_{B}(B|K), (B|K)|_{B}(A|H)\}$ is $\Pi=\Pi' \cup \Pi''$, where 
$\Pi'=\{(x,y,z,\mu): x\in(0,1], y\in [0,1], z\in[z', z''], \mu=\frac{z}{x}\}$
 with $ z'=0$,    $z''=1$, and
$\Pi''=\{(0,y,0,\mu): (y,\mu)\in[0,1]^2\}.$
\end{theorem}
\begin{proof}
See \ref{SEC:APPB}.
\end{proof}
\begin{theorem}\label{THM:LU_B}
	Let $A,B,H,K$ be any logically independent events. Given a coherent assessment $(x,y)$ on $\{A|H,B|K\}$, for the iterated conditional $(B|K)|_{B}(A|H)$ the  extension $\mu_B=\prev((B|K)|_{B}(A|H))$ is coherent if and only if 
	$\mu_B \in [\mu_B', \mu_B'']$, where	\begin{equation}\label{EQ:P4B}
	\begin{array}{ll}
	\mu_B'=0,\;\;
	\mu_B''=\left\{
	\begin{array}{ll}
	\frac{1}{x},& \text{ if } 0<x<1,\\
	1,& \text{ if } x=0 \vee x=1.\\
	\end{array}
	\right.
	\end{array}
	\end{equation} 
\end{theorem}
\begin{proof}
See \ref{SEC:APPB2}.
\end{proof}
\begin{remark}
We notice that the lower and upper bounds $\mu'_B$ and  $\mu''_B$ for $(B|K)|_K(A|H)$, given in 
(\ref{EQ:P4B}, do not coincide with  the lower and upper bound $\mu'$ and $\mu''$ for $B|A$, given in (\ref{EQ:LU}).  In particular, the extension $\mu_B=\frac{1}{x}$ to  $(B|K)|_{B}(A|H)$ of $(x,y)$ on $\{A|H, B|K\}$, with $0<x<1$, is coherent. Indeed, as $\mu''_B=\frac{1}{x}$, from (\ref{EQ:P4B}) the assessment  $(x,y,\frac{1}{x})$ on $\{A|H, B|K, (B|K)|_{B}(A|H)\}$ is coherent. However, we recall that $\mu=\frac{1}{x}$ is not a coherent extension to  $B|A$ of $(x,y)$ on $\{A, B\}$, with $0<x<1$. Indeed,  as $\frac{1}{x}>1\geq \mu''=min\{x,y\}$, from  (\ref{EQ:LU}) the assessment $(x,y,\frac{1}{x})$ on $\{A, B, B|A\}$ is not coherent. 
Therefore, property P4 is not satisfied by $|_B$.
\\
We also notice that, when $P(A|H)>0$, as $(A|H)\wedge_B(B|K)=AH|BK$, from (\ref{EQ:COMPTHM}) it holds that
\begin{equation}\label{EQ:HYPERP}
\prev[(B|K)|_{B}(A|H)]=\frac{
P[(A|H)\wedge_{B}(B|K)]}{P(A|H)}=\frac{P(AH|BK)}{P(A|H)}.
\end{equation}
We  point out that the prevision of $(B|K)|_{B}(A|H)$ given in  (\ref{EQ:HYPERP}) coincides with the probability of the \emph{super-conditional event}  $B_K|A_H$, denoted by $P^*(B_K|A_H)$, introduced in \cite[formula (20)]{BrGi85}.
As  $P^*(B_K|A_H)$  can assume values bigger than 1, in  \cite{BrGi85} it has been called  \emph{conditional hyper-probability}.
\end{remark}

\paragraph{The iterated conditioning $|_S$}
We first  show that, when we  evaluate $\prev[(B|K)|_S(A|H)]$,  
if $P(A|H)=0$, then 
coherence requires that $P[(A|H)\wedge_S(B|K)]=0$. 
\begin{remark}\label{coherenceS}
We recall that, given  any $(x,y)\in[0,1]^2$, where $x=P(A|H)$ and $y=P(B|K)$,  the interval of coherent extensions on $z_S=P[(A|H)\wedge_{S}(B|K)]$ is $[z', z'']$ where $ z'=\max\{x+y-1, 0\}$ and $
z''=\left\{
\begin{array}{ll}
	\frac{x+y-2xy}{1-xy}, &\text{ if } (x,y)\neq (1,1),\\
	1, &\text{ if } (x,y)= (1,1).\\
\end{array}
\right.
$ (Table \ref{tab:LUB}). In particular it is coherent  to assess  $z_S>0$ when $x=0$. However, for the object 
$(B|K)|_S(A|H)$ coherence  also requires that $z_S=\mu_S x$ (see \ref{prevision}). Then, 
coherence requires that $z_S=0$ when we consider the assessment $(0,y,z_S,\mu_S)$ on
$\{A|H,B|K, (A|H)\wedge_{S}(B|K),(B|K)|_{S}(A|H)\}$.

\end{remark}  
 
\begin{theorem}\label{THM:PIS}
Let $A$, $B$, $H$, $K$ be any logically independent events. The set $\Pi$ of all the coherent assessment $(x,y,z,\mu)$ on the family $\F=\{A|H, B|K, (A|H)\wedge_{S}(B|K), (B|K)|_{S}(A|H)\}$ is $\Pi=\Pi' \cup \Pi''$, where 
\begin{equation}\label{EQ:PI'}
\begin{array}{ll}
\Pi'=\{(x,y,z,\mu): x\in(0,1], y\in [0,1], z\in[z', z''], \mu=\frac{z}{x}\},\\
 \text{with } z'=\max\{x+y-1, 0\},  
z''=\left\{
\begin{array}{ll}
	\frac{x+y-2xy}{1-xy}, &\text{ if } (x,y)\neq (1,1),\\
	1, &\text{ if } (x,y)= (1,1),\\
\end{array}
\right.
\end{array}
\end{equation} 
 and
\begin{equation}\label{EQ:PI''}
\Pi''=\{(0,y,0,\mu): y\in[0,1], \mu \geq 0\}.
\end{equation}
\end{theorem}
 \begin{proof}
See \ref{SEC:APPS}
 \end{proof} 
\begin{remark}
From Theorem \ref{THM:PIS} we observe that when  $x=0$ the set of all coherent assessments $(x,y,z,\mu)$ on $\F$ is  $\Pi''$ given in \ref{EQ:PI''}. Then,  in this case, $\mu=\prev[(B|K)|_{S}(A|H)]=\prev[(AHBK+AH\no{K}+(1+\mu)\no{H}BK)|(AH\vee \no{H}BK)]$, is coherent for every value  $\mu\geq 0$.
\end{remark}
Based on Theorem~\ref{THM:PIS}, when $x>0$, we obtain the following result on  the  lower and upper bounds for $\prev[(B|K)|_{S}(A|H)]$.
\begin{theorem}\label{THM:LUS}
	Let $A,B,H,K$ be any logically independent events. Given a coherent assessment $(x,y)$ on $\{A|H,B|K\}$, with $x\neq0$, for the iterated conditional $(B|K)|_{S}(A|H)$ the  extension $\mu_S=\prev((B|K)|_{S}(A|H))$ is coherent if and only if 
	$\mu_S \in [\mu_S', \mu_S'']$, where
\begin{equation}\label{EQ:P4S}
	\mu_S'=\max\{\frac{x+y-1}{x}, 0\} \text{ and }\mu_S''=\left\{
	\begin{array}{ll}
	\frac{x+y-2xy}{x(1-xy)},& \text{ if } (x,y)\neq (1,1);\\
	1,& \text{ if } (x,y)=(1,1).\\
	\end{array}
	\right.
\end{equation} 
\end{theorem}
\begin{proof}
 See \ref{SEC:APPS2}
\end{proof}

\begin{remark}
We observe that the lower and upper bounds $\mu'_S$ and  $\mu''_S$ for $(B|K)|_S(A|H)$, given in 
(\ref{EQ:P4S}), do not coincide with  the lower and upper bound $\mu'$ and $\mu''$ for $B|A$, given in (\ref{EQ:LU}). Then, formula P4 is not satisfied by $|_S$.
In particular the extension $\mu_S=\frac{1}{x}$ to $(B|K)|_S(A|H)$ of $(x,1)$ on $\{A|H,B|K\}$,  with $0<x<1$,
is coherent because
$\mu_S''=\frac{1}{x}$. However, the assessment $\mu=\frac{1}{x}$ is not a coherent extension on $\mu=P(B|A)$ of $P(A)=x$, $P(B)=1$,  because by (\ref{EQ:LU}) it holds that $\mu''\leq 1< \frac{1}{x}$.
\end{remark}

\paragraph{The iterated conditioning $|_{gs}$}

Finally, differently from the other iterated conditionals, we recall that   $|_{gs}$ satisfies P4 (\cite[Theorem 4]{SPOG18}). Indeed,    given a coherent assessment $(x,y)$ on $\{A|H, B|K\}$, under logical independence, for the iterated conditional $(B|K)|_{gs}(A|H)$ the extension $\mu=\prev((B|K)|_{gs}(A|H))$ is coherent if and only if  $\mu \in [\mu_{gs}', \mu_{gs}'']$, where
$
\small
\mu_{gs}'=\left\{
\begin{array}{ll}
	\max\{\frac{x+y-1}{x},0\}, &\text{ if } x\neq 0,\\
	0, &\text{ if } x= 0,\\
\end{array}
\right.,\;\;
\mu_{gs}''=\left\{
\begin{array}{ll}
	\min\{1,\frac{y}{x}\}, &\text{ if } x\neq 0,\\
	1, &\text{ if } x= 0,\\
\end{array}
\right.
$ 
which coincide with $\mu'$ and $\mu''$ given in (\ref{EQ:LU}), respectively.
As we can see from Table \ref{TAB:LU-ALL}, only the iterated conditioning $|_{gs}$ satisfies property P4. 
	\begin{table}[!ht]
	\centering
	\begin{tabular}{c|c}
Iterated conditioning         & Interval of coherent extensions   \\ \hline
    \rule[-4mm]{0mm}{1cm}
		$ |_C$ &  $[0,1]$\\  \hline
		\rule[-4mm]{0mm}{1cm}
		$|_{dF}$ &  $[0,1]$    \\ \hline
		\rule[-4mm]{0mm}{1cm}
		$|_{F}$ &  $[0,1]$    \\ \hline	
		\rule[-6mm]{0mm}{1.5cm}
		$|_{K}$ &  $\left [0,\; \left\{
	\begin{array}{ll}
	\min\left\{1,\tfrac yx\right\},& \text{ if } x\neq 0;\\
	1,& \text{ if } x=0.\\
	\end{array}
	\right.\right]$\\ \hline
	\rule[-6mm]{0mm}{1.5cm}
		$|_{L}$ &  $\left [0,\; \left\{
	\begin{array}{ll}
	\min\left\{1,\tfrac yx\right\},& \text{ if } x\neq 0;\\
	1,& \text{ if } x=0.\\
	\end{array}
	\right.\right]$\\	\hline
	\rule[-6mm]{0mm}{1.5cm}
		$|_{B}$ &  $\left [	0,\; \left\{
	\begin{array}{ll}
	\frac{1}{x},& \text{ if } 0<x<1,\\
	1,& \text{ if } x=0 \vee x=1.\\
	\end{array}
	\right.\right]$    \\	\hline
		\rule[-9.5mm]{0mm}{2.5cm}
		$|_{S}$ &  $\left\{\begin{array}{ll} \left [\max\{\frac{x+y-1}{x}, 0\}, \; \left\{
	\begin{array}{ll}
	\frac{x+y-2xy}{x(1-xy)},& \text{ if } (x,y)\neq (1,1);\\
	1,& \text{ if } (x,y)=(1,1).\\
	\end{array}
	\right. \right] & \text{ if } x\neq 0;\\ \rule{0mm}{0.5cm}
	\mu \geq 0,  & \text{ if } x=0;
	\end{array}\right.$  \\	\hline	
		\rule[-6mm]{0mm}{1.5cm}
		$|_{gs}$ &  $\left[ \left\{
\begin{array}{ll}
	\max\{\frac{x+y-1}{x},0\}, &\text{ if } x\neq 0,\\
	0, &\text{ if } x= 0,\\
\end{array}, \;
\right.
\left\{
\begin{array}{ll}
	\min\{1,\frac{y}{x}\}, &\text{ if } x\neq 0,\\
	1, &\text{ if } x= 0.\\
\end{array}
\right.
		\right]$  \\	
  \hline				
	\end{tabular}
	\caption{ Interval $[\mu_i',\mu_i'']$  of coherent extensions of the assessment $(x,y)\in[0,1]^2$ on $\{A|H,B|K\}$ to the iterated conditional $(B|K)|_i(A|H)$, $i\in\{C, dF, F, K, L, B, S, gs\}$, under the assumption that  $A,H,B,K$ are logically independent.}	
	\label{TAB:LU-ALL}
\end{table}

\subsection{Import-Export principle}
 \begin{table}[h]
	\begin{small}
		\def\True{$1$}
		\def\False{$0$}
		\def\Void{$z$}
		\def\Voidx{$x$}
		\def\Voidy{$y$}
		\centering
		\begin{tabular}{|l|c|c|c|c|c|c|}
			\hline
			                                  & $(B|K)|_{K}A$ & $(B|K)|_{L}A$ & $(B|K)|_{B}A$ & $(B|K)|_{S}A$ & $(B|K)|_{gs}A$ & $B|AK$ \\ \hline
			$ABK  $                           &       1       &       1       &       1       &       1       &       1        &   1    \\
			$	{A}{\overline{B}K}$             &       0       &       0       &       0       &       0       &       0        &   0    \\
			$	A\overline{K}  $                &   $x\mu_K$    &   $x\mu_L$    &   $x\mu_B$    &       1       &      $y$       &      $z$  \\
			$	{\overline{A}}BK$               &    $\mu_K$    &    $\mu_L$    &    $\mu_B$    &    $\mu_S$    &   $\mu_{gs}$   & $z$       \\
			$	{\overline{A}}\;\overline{B}K$    &    $\mu_K$    &    $\mu_L$    &    $\mu_B$    &    $\mu_S$    &   $\mu_{gs}$   & $z$       \\
			$	{\overline{A}} \, \overline{K}$ &    $\mu_K$    &    $\mu_L$    & $\mu_B(1+x)$  &    $\mu_S$    &   $\mu_{gs}$   &     $z$   \\ \hline
		\end{tabular}
		\caption{Numerical values   of  $(B|K)|_{i}A$, $i\in\{K,L,B,S,gs\}$ and of the conditional event $B|AK$. We denotes $x=P(A)$, $y=P(B|K)$,   $z=P(B|AK)$ 
			and $\mu_i=\prev[(B|K)|_i A]$, $i\in\{K,L,B,S,gs\}$.  We observe that $(B|K)|_{i}A\neq B|AK$, $i\in\{K,L,B,S,gs\}$.
		}
		\label{TAB:IE}
	\end{small}
\end{table} 
In this section we show that none of the iterated conditional $\{K,L,B,S,gs\}$ satisfy the Import-Export principle.
We recall that the Import-Export principle is valid when $(B|K)|A=B|AK$ (see equation (\ref{EQ:IE}). 
We  remind that, in agreement with \cite{adams75,Kauf09} and differently from \cite{McGe89}, for the iterated conditional $(B|K)|_{gs}(A|H)$ the  Import-Export principle is not valid. As a consequence, 
as shown in \cite{GiSa14} (see also \cite{SGOP20,SPOG18}), 
Lewis’ triviality results (\cite{lewis76}) are avoided by $|_{gs}$. 
For what concerns the iterated conditionals $|_K, |_L, |_B, |_S$, it holds that 
$(B|K)|_{i}A\neq B|AK$, $i\in \{K,L,B,S\}$, 
because there are some constituents where the two objects may assume different values (see Table \ref{TAB:IE}). For instance, when $A\no{K}$ is true  and $z=P(B|AK)\neq 1$, it follows that $(B|K)|_{S}A=1\neq z=B|AK$. 
Then,
the
none of the  iterated conditional $|_K, |_L, |_B, |_S, |_{gs}$ satisfy
the Import-Export principle. However, we observe that Import-Export principle could be satisfied under some suitable logical relations among the events $A,B,K$. For instance, if $A\no{K}=\emptyset$, it can easily proved that $(B|K)|_iA=B|AK$, $i\in\{K,L,S,gs\}$. Thus,  when $A\no{K}=\emptyset$, the Import-Export principle is satisfied by $(B|K)|_KA,$ $(B|K)|_LA$, $(B|K)|_SA$, and $(B|K)|_{gs}A$.
\subsection{Generalized versions of  Bayes' Rule } \label{SEC:Bayes}
In this section, 
by exploiting property P4,
we  analyze  generalized versions of  Bayes' Rule for the iterated conditioning  $|_K, |_L, |_B, |_S$, and  $|_{gs}$. 

From (\ref{EQ:COMPTHM}) it follows that 
$\prev[(B|K)\wedge_{i}(A|H)]=\prev[(B|K)|_i(A|H)]P(A|H)=\prev[(A|H)|_i(B|K)]P(B|K)$, $i\in\{K,L,B,S,gs\}$. Then,  when $P(A|H)>0$ it holds that 
\begin{equation}\label{EQ:BAYES}
\prev[(B|K)|_{i}(A|H)]=\frac{\prev[(A|H)|_i(B|K)]P(B|K)}{P(A|H)},\, i\in\{K,L,B,S,gs\}.
\end{equation}
Formula (\ref{EQ:BAYES})  is a generalization for the of   the following well-known version of  Bayes's Rule
\[
P(B|A)=\frac{P(A|B)P(B)}{P(A)},\;\; \text{if } P(A)>0,
\]
where the events $A$, $B$ are replaced by the conditional events $A|H$, $B|K$, respectively, and the conditioning operator $|$ is replaced by the iterated conditioning $|_i$, $i\in\{K,L,B,S,{gs}\}$.
We also recall that,  given two events $A$ and $B$, it holds that
$
A=AB \vee A\no{B},
$
and hence $P(A)=P(AB)+P(A\no{B})=P(A|B)P(B)+P(A|\no{B})P(B)$. Then, we obtain 
 the following version of Bayes' Rule
\begin{equation}\label{EQ:BAYES2}
P(B|A)=\frac{P(A|B)P(B)}{\rule{0pt}{0.8em}P(A|B)P(B)+P(A|\no{B})P(\no{B})}.
\end{equation}
We now check, for each iterated conditioning $|_i$, $i\in\{K,L,B,S,{gs}\}$, 
  the validity of the following generalized version of  formula (\ref{EQ:BAYES2})
\begin{equation}
\label{EQ:BAYES2GEN}
\prev[(B|K)|(A|H)]=
    \frac{\prev{((A|H)|(B|K))P(B|K)}}{\prev{((A|H)|(B|K))P(B|K)}+\prev{((A|H)|(\no{B}|K))P(\no{B}|K)}}.
\end{equation}
Based on (\ref{EQ:COMPTHM}),for each 
$i\in\{K,L,B,S,gs\}$, it follows that 
\begin{equation}
\label{EQ:DISTR}
P{((A|H)\wedge_{i} (B|K))}+ P{((A|H)\wedge_{i} (\no{B}|K))}=\prev{((A|H)|_i(B|K))P(B|K)}+\prev{((A|H)|_i(\no{B}|K))P(\no{B}|K)}.
\end{equation}
We recall that when $A$, $B$, and $\no{B}$ in the equality $A=(A\wedge B)\vee (A\wedge \no{B})$ are replaced  by the conditional events $A|H, B|K$, and $\no{B}|K$,  respectively, and the conjunction $\wedge$ is replaced by $\wedge_i$, $i\in\{K,L,B,S\}$, it holds  that
the corresponding equality is not satisfied. Indeed we have that  
\[
(A|H)\neq [(A|H)\wedge_i (B|K)]\vee_i [(A|H)\wedge_i (\no{B}|K)], \;\; i\in\{K,L,B,S\}.
\]
More precisely, it holds that (see \cite[Section 4.1]{GiSa22})
\begin{itemize}
	\item $[(A|H)\wedge_K (B|K)]\vee_K [(A|H)\wedge_K (\no{B}|K)]=AHK|(AHK\vee \no{A}H) \neq A|H;$
	\item $[(A|H)\wedge_L (B|K)]\vee_L [(A|H)\wedge_L (\no{B}|K)]=AHK|(H\vee \no{K}) \neq A|H;$
	\item $[(A|H)\wedge_B (B|K)]\vee_B [(A|H)\wedge_B (\no{B}|K)]=A|(HK) \neq A|H;$
	\item $[(A|H)\wedge_S (B|K)]\vee_S [(A|H)\wedge_S (\no{B}|K)]=(A \vee \no{H})|(H\vee K) \neq A|H.$
\end{itemize}
Then, for each $i\in\{K,L,B,S\}$,  the conditional probability $P(A|H)$ does not necessarily coincide with 
\[
P{((A|H)\wedge_{i} (B|K))}+ P{((A|H)\wedge_{i} (\no{B}|K))}
\]
and hence, from (\ref{EQ:DISTR}), the probability  $P(A|H)$ does not necessarily coincide with 
\[
\prev{((A|H)|_i(B|K))P(B|K)}+\prev{((A|H)|_i(\no{B}|K))P(\no{B}|K)}.
\] Thus,  for each iterated conditioning  $|_i$, $i\in\{K,L,B,S\}$, the generalized version of Bayes's Rule given in formula  (\ref{EQ:BAYES2GEN}) is not satisfied.
However, we recall that (\cite{GiSa20})  
\[(A|H)=(A|H)\wedge_{gs} (B|K)+(A|H)\wedge_{gs} (\no{B}|K).\]
Then, by the linearity property of a coherent prevision and by (\ref{EQ:DISTR}), it follows that 
\[
P(A|H)=\prev{((A|H)\wedge_{gs} (B|K))}+ \prev{((A|H)\wedge_{gs} (\no{B}|K))}=\prev{((A|H)|_{gs}(B|K))P(B|K)}+\prev{((A|H)|_{gs}(\no{B}|K))P(\no{B}|K)}.
\]
 Hence, when $P(A|H)>0$, it holds that
\begin{equation}
\prev[(B|K)|_{gs}(A|H)]=
    \frac{\prev{((A|H)|_{gs}(B|K))P(B|K)}}{\prev{((A|H)|_{gs}(B|K))P(B|K)}+\prev{((A|H)|_{gs}(\no{B}|K))P(\no{B}|K)}}.
\end{equation}
 Therefore, the  generalization of the second version of  Bayes' Rule
only holds for $|_{gs}$ and does not hold for  $|_K, |_L, |_B,$ and  $|_S$.


\section{Generalized version of Modus Ponens and Two premise centering}
\label{SEC:MP}
In this section, 
for selected definitions of the iterated conditioning with possible value in $[0,1]$, based on the probability propagation rules obtained in the previous sections and exploited for checking the validity of property P4 (see Table \ref{TAB:LU-ALL}), we study the p-validity of the generalization of two  inference rules: Modus Ponens and two-premise centering.  We first recall the notions of p-consistency and  p-entailment for  conditional random quantities, which take values in a finite subset of $[0,1]$  (\cite{SPOG18}). These notion  are based  on the notions of p-consistency and p-entailment given for  conditional events by Adams
(\cite{adams75}) and also studied in the setting  of coherence (see, e.g., \cite{gilio02}).
\begin{definition}\label{DEF:PC}
	Let $\mathcal{F}_n = \{X_i|H_i \, , \; i=1,\ldots,n\}$ be  a family of $n$  conditional random quantities which take values in a finite subset of $[0,1]$. Then, $\mathcal{F}_n$ is  {\em p-consistent} if and only if,
	the (prevision) assessment $(\mu_1,\mu_2,\ldots,\mu_n)=(1,1,\ldots,1)$ on $\mathcal{F}_n$ is coherent.
\end{definition}
\begin{definition}\label{DEF:PE}
	A p-consistent family  $\mathcal{F}_n = \{X_i|H_i \, , \; i=1,\ldots,n\}$ {\em p-entails} a conditional random quantity $X|H$, which takes values in a finite subset of $[0,1]$, denoted by $\mathcal{F}_n \; \Rightarrow_p \; X|H$,
	if and only if  for any  coherent (prevision) assessment $(\mu_1,\ldots,\mu_n,z)$ on $\mathcal{F}_n \cup \{X|H\}$: if $\mu_1=\cdots=\mu_n=1$, then  $z=1$.
\end{definition}
We say that the inference from a p-consistent family of premises $\F_n$ to a conclusion $X|H$ is \emph{p-valid} if and only if  $\mathcal{F}_n$ p-entails $X|H$.

We recall that  as the iterated conditionals $(B|K)|_{C}(A|H)$, $(B|K)|_{dF}(A|H)$, and $(B|K)|_{F}(A|H)$ are conditional events, their indicators take values in the interval [0,1]. Moreover,  
from Table \ref{TAB:ITECMPPRB} we observe  that the iterated conditionals $(B|K)|_{K}(A|H)$, $(B|K)|_{L}(A|H)$ and $(B|K)|_{gs}(A|H)$ take  values in  $[0,1]$.  On the other hand,  the iterated conditionals $(B|K)|_{B}(A|H)$ and $(B|K)|_{S}(A|H)$  may take values outside the interval $[0,1]$. Thus, in order to examine the p-validity of generalized inference rules, we will only consider the iterated conditionals $(B|K)|_{i}(A|H)$,  $i\in\{C,dF,F,K,L,gs\}$.  
\subsection{Modus Ponens}
We recall that, given two events $A$ and $B$ the  Modus Ponens inference, with (p-consistent) premise set $\{A,B|A\}$ and conclusion $B$,   is p-valid (\cite{Wagner04}, see also \cite{GiPS20,hailperin96}), that is  
\begin{equation}
\label{EQ:MP}
\{A,B|A\}\Rightarrow_p B.
\end{equation}
For each $i\in\{C,dF,F,K,L,gs\}$, we will study the p-validity of the generalized version of Modus Ponens obtained when 
the events $A,B$ are replaced by the conditional events $A|H, B|K$ and the conditional event $B|A$ is replaced by the iterated conditional 
$(B|K)|_i(A|H)$\footnote{
The particular case where $K=\Omega$ and $|_i=|_{gs}$ has been studied in \cite{SaPG17}.}.
An example of this generalization is (see also \cite{gibbard81,SaPG17}): 
\begin{quote}
$\overbrace{\text{\em 
The cup is fragile when it is made of glass
}}^{A|H}$.\\
If $\overbrace{\text{\em 
the cup is fragile when it is made of glass
}}^{A|H}$, then $\overbrace{\text{\em the cup  breaks if dropped}}^{B|K}$.\\
Therefore, $\overbrace{\text{\em the cup  breaks if dropped}}^{B|K}$.
\end{quote}

Then, we
check the   p-validity of the  generalized version of Modus Ponens where the premise set is  $\{A|H, (B|K)|_i(A|H)\}$ and the conclusion is $B|K$, that is  
\begin{equation}
\label{EQ:GMP}
P(A|H)=1, \; \prev[(B|K)|_i(A|H)]=1\;\Longrightarrow P(B|K)=1.
\end{equation}
It is easy to verify (see Table \Ref{TAB:LU-ALL}) that the family  $\{A|H, (B|K)|_i(A|H)\}$ is p-consistent, $i\in\{C,dF,F,K,L,gs\}$.
We will show that 
the generalized version of Modus Ponens is p-valid for $|_i\in\{|_K,|_L,|_{gs}\}$.
Indeed, 
in these cases 
from Table \ref{TAB:LU-ALL} it follows that   $\prev[(B|K)|_i(A|H)]\leq \mu_i''=\min\{1,\frac{y}{x}\}=y=P(B|K)$, when  $P(A|H)=x=1$. Then, $P(B|K)$ is necessarily equal to 1, when $P(A|H)=1$ and  $\prev[(B|K)|_i(A|H)]=1$
and hence  formula  
(\ref{EQ:GMP}) is satisfied. Thus, 
\begin{equation}
\{A|H, \; (B|K)|_i(A|H)\}\;\Rightarrow_p B|K,\;\;
i\in\{K,L,{gs}\}.
\end{equation}

For $i\in\{C,dF,F\}$, as illustrated in   Table \ref{TAB:LU-ALL},
the interval of coherent extensions $[\mu_i',\mu_i'']$ on $(B|K)|_i(A|H)$ of the assessment $(1,0)$ on $\{A|H,B|K\}$ coincides with  the  unit interval $[0,1]$. In particular the assessment  
$(1,1,0)$ on $\{A|H,(B|K)|_i(A|H),(B|K)\}$  is coherent and hence   formula  
(\ref{EQ:GMP}) is not satisfied for $i \in \{C,dF,F\}$, that is 
 \begin{equation}
 	\{A|H, \; (B|K)|_i(A|H)\}\;\not \Rightarrow_p B|K,\;\;
 	i\in\{C,dF,F\}.
 \end{equation}

Concerning the example above,  
from  
$P(${\em the cup is fragile when it is made of glass}$)=1$ and 
$\prev($If 
{\em the cup is fragile when it is made of glass,}
 then  {\em the cup  breaks if dropped}$)=1$, 
 it follows  (as a natural result) that $P(${\em the cup  breaks if dropped}$)=1$   
when the iterated conditioning is interpreted as  $|_K,|_{L},|_{gs}$, because (\ref{EQ:GMP}) is satisfied in these cases. However, the same conclusion does not follow (which is strange) when the iterated conditioning is interpreted as  $|_C,|_{dF}$, and $|_F$ because (\ref{EQ:GMP}) does not hold in these cases.\\

\subsection{Two-premise centering}
We recall that the inference of {\em two-premise centering}, that is
inferring {\em if $A$ then $B$} from the two separate premises $A$ and $B$, is p-valid (\cite{SPOG18}). Indeed,   
 given two events $A$ and $B$  it holds that 
\begin{equation}
\label{EQ:2CENT}
\{A,B\}\Rightarrow_p B|A.
\end{equation}

For each $i\in\{C,dF,F,K,L,gs\}$, we will study the p-validity of the generalized version of  two-premise centering, where the unconditional events $A$, $B$ are replaced by the conditional events $A|H$,$B|K$, respectively, and the conditional event $B|A$ is replaced by the iterated conditional $(B|K)|_i(A|H)$.
An example of this generalization is

\begin{quote}
$\overbrace{\text{\em 
The cup is fragile when it is made of glass
}}^{A|H}$.\\
$\overbrace{\text{\em The cup  breaks if dropped}}^{B|K}$.

Therefore, 
if $\overbrace{\text{\em 
the cup is fragile when it is made of glass
}}^{A|H}$, then $\overbrace{\text{\em the cup  breaks if dropped}}^{B|K}$.\\
\end{quote}
For each $i$, we  consider as  premise set $\{A|H, B|K\}$, as conclusion the iterated conditional $(B|K)|_i(A|H)$ and we check the p-validity of generalized version of two-premise centering, that is 
\begin{equation}
\label{EQ:GTP}
P(A|H)=1,\, P(B|K)=1\, \Longrightarrow 
\prev[(B|K)|_i(A|H)]=1.
\end{equation}
Of course, as the assessment $(x,y)=(1,1)$ on $\{A|H,B|K\}$  is coherent, the premise set $\{A|H,B|K\}$ is p-consistent. 
We recall that  generalized two-premise centering is p-valid for  $|_{gs}$ (\cite[Equation (13)]{SPOG18}), that is 
\begin{equation}
	\{A|H,B|K\}\Rightarrow_p (B|K)|_{gs}(A|H)].
\end{equation}
 The previous  result can be also obtained  from Table \ref{TAB:LU-ALL}, by observing that    the lower bound  of the coherent extensions on $(B|K)|_{gs}(A|H)$ of the assessment $x=y=1$ is $\mu_i'=\max\{\frac{x+y-1}{x},0\}=1$. 
 Moreover, 
  for $i\in \{C,{dF},F,K,L\}$, 
  from Table \ref{TAB:LU-ALL} it follows that 
  that any value in $[0,1]$ is a coherent extension on  
 $(B|K)|_i(A|H)$
  of
 $P(A|H)=P(B|K)=1$.
 Thus, (\ref{EQ:GTP}) is not satisfied by $|_C,|_{dF},|_F,|_K,|_L,$ and hence
 \begin{equation}
 	\{A|H,B|K\}\not \Rightarrow_p (B|K)|_{i}(A|H)],\;\;\; i\in\{C,{dF},F,K,L\}.
 \end{equation}

Concerning the example above,  
from  
$P(${\em the cup is fragile when it is made of glass}$)=1$ and 
$P(${\em the cup  breaks if dropped}$)=1$ 
it follows that 
$\prev($If 
{\em the cup is fragile when it is made of glass,}
 then  {\em the cup  breaks if dropped}$)=1$, 
when the iterated conditioning is interpreted as  $|_{gs}$, because the generalized two-premise centering  is p-valid. However, the same conclusion does not follow  when the iterated conditioning is interpreted as  $|_C,|_{dF},|_F,|_K,|_L$, because the generalized two-premise centering  is p-invalid  in these cases.
\section{Conclusions}
\label{SEC:CONCLUSIONS}
We recalled the trivalent logics of Kleene-Lukasiewicz-Heyting-de Finetti, Lukasiewicz, Bochvar-Kleene, and  Soboci\'nski 
and the notion of  compound conditional as conditional random quantity.
We  considered four basic logical and probabilistic properties, P1-P4,
 valid for events and conditional events.
We generalized them by replacing events $A$ and $B$ with conditional events $A|H$ and $B|K$ and we checked their validity and the validity of the Import-Export principle for selected notions of iterated conditioning. 
In particular, we studied  the  iterated  conditioning  introduced in trivalent logics by 
Cooper-Calabrese $(|_C)$, de Finetti $(|_{dF})$, and Farrel $(|_F)$, by also focusing on the numerical representation of the truth-values. 
We observed that the
notions of  conjunction and disjunction of conditional events used by Cooper and Calabrese  coincide with  $\wedge_S$ and $\vee_S$, respectively. Farrell and  de Finetti   defined two different structures of iterated conditioning in the same trivalent logic where  conjunction and disjunction are $\wedge_K$ and  $\vee_K$, respectively.
We computed the set of coherent probability assessments on the families of events $\{A|H, (B|K)|_{C}(A|H), (A|H)\wedge_{S}(B|K)\}$,
$\{A|H, B|K, (B|K)|_{C}(A|H)\}$, as well as $\{A|H, (B|K)|_{i}(A|H), (A|H)\wedge_{K}(B|K)\}$ and  $\{A|H, B|K, (B|K)|_{i}(A|H)\}$ with $i\in \{dF, F\}$. 
In Table \ref{tab:recap} we summarize the results for  the properties given in Section \ref{SEC:3}. We observe that $|_C$, $|_{dF}$,  and  $|_F$ satisfy the Import-Export principle and none of these objects, defined in the framework of trivalent logics, satisfy the compound probability theorem (P3).

By exploiting  the 
 structure 
 $\Box|\bigcirc = \Box\wedge \bigcirc +\prev(\Box|\bigcirc)\no{\bigcirc}$
 used in order to define $|_{gs}$ from the conjunction $\wedge_{gs}$,
for each 
$i\in\{K,L,B,S\}$
we defined the iterated conditioning $|_i$ from the 
conjunction $\wedge_{i}$.
We observed that the iterated conditionals $(B|K|_i(A|H)$, $i\in\{K,L,B,S,gs\}$, 
 are all conditional  random quantities (not always in [0,1]) which satisfy  the compound prevision theorem (property P3). We also noticed that properties P1, P2 and the non validity of  the Import-Export principle are satisfied by all the iterated conditional
$|_{K},|_{L},|_{B},|_{S}, |_{gs}$ 
(see Table \ref{tab:recap}).
 However, property P4 is not satisfied by $|_K, |_L, |_B, |_S$.
We observed  that all the basic logical and probabilistic properties are  satisfied only by the  iterated conditioning $|_{gs}$. 
\def\No{}
\def\Yes{$\checkmark$ }
\begin{table}[h]
    \centering
    \begin{tabular}{|c|c|c|c|c|c|c|c|c|c|}
        \hline
        & {Property} & {$|_C$} & {$|_{dF}$} & {$|_F$} & {$|_K$} & {$|_L$} & {$|_B$} & {$|_S$} & {$|_{gs}$}\\
         \hline
        {No IE} & 
       $(B|K)|A\neq B|AK$ &
        {\No} & {\No} & {\No} & {\Yes} & {\Yes} & {\Yes} & {\Yes} & {\Yes}\\ 
        \hline
         {P1} & $(B|K)|(A|H)=[(A|H)\wedge (B|K)]|(A|H)$& {\No} & {\Yes} & {\Yes} & {\Yes} & {\Yes} & {\Yes} & {\Yes} & {\Yes}\\
         \hline
         {P2} & $(A|H)\wedge (B|K)\leq    (B|K)|(A|H)$ &{\No} & {\Yes} & {\Yes} & {\Yes} & {\Yes} & {\Yes} & {\Yes} & {\Yes}\\ 
        \hline
        {P3} & $\prev[(A|H)\wedge (B|K)]=\prev[(B|K)|(A|H)]P(A|H)$& {\No} & {\No} & {\No} & {\Yes} & {\Yes} & {\Yes} & {\Yes} & {\Yes}\\ 
        \hline
        {P4} & 
Lower and upper bounds for $(B|K)|(A|H)$& {\No} & {\No} & {\No} & {\No} & {\No} & {\No} & {\No} & {\Yes}\\ 
     
        \hline
    \end{tabular}
    \caption{
    Properties P1-P4, the non validity of Import-Export Principle (No IE), and  iterated conditionals in their own logic. The symbol \Yes  means that the property is satisfied. The blank space means that the property is not satisfied.}
    \label{tab:recap}
\end{table}
We also showed that a generalized version of the Bayes' Rule $P(B|A)=\frac{P(A|B)P(B)}{P(A)}$ for iterated conditionals is satisfied by $|_i$ with $i\in\{K,L,B,S,gs\}$. However, a generalized version of the  Bayes' Rule in the following version  
$P(B|A)=\frac{P(A|B)P(B)}{\rule{0pt}{0.8em}P(A|B)P(B)+P(A|\no{B})P(\no{B})}
$
for iterated conditionals, only holds for  $|_{gs}$ and does not hold for  $|_K, |_L, |_B, |_S$.

We discussed the implications of the  obtained results on the probability propagation rules from $\{A|H,B|K\}$ to the iterated conditional $(B|K)|_i(A|H)$, 
by studying  the p-validity of the generalized versions of  Modus Ponens and two-premise centering, which involve conditional events and iterated conditionals.
We observed that for $|_B$ and $|_S$, as the iterated conditionals $(B|K)|_B(A|H)$ and $(B|K)|_S(A|H)$ can take value outside the interval [0,1], the study of p-validity is meaningless.
For the remaining ones $|_C,|_{dF},|_F,|_K,|_L$, and $|_{gs}$, we showed that   only for the iterated conditioning $|_{gs}$  turns out that both generalized versions of the inference rules are p-valid. 
We also examined two examples of the generalized  inference rules in natural language  by observing that, when 
we adopt the notion of  iterated conditioning $|_i$, with  $|_i\in\{|_C,|_{dF},|_F,|_K,|_L\}$,   some results which are  counterintuitive  in commonsense reasoning can be obtained.

Therefore, based on the results illustrated above, only the iterated conditional $|_{gs}$, which is  based on the conjunction $\wedge_{gs}$ introduced in the framework of conditional random quantities, preserves all the  basic logical and probabilistic properties;  moreover,  Lewis’ triviality results  are avoided in particular because the Import-Export Principle is not satisfied. Differently, by using the other iterated conditionals, which are based on different logics,   certain probabilistic properties are not preserved. Then, the `probability' of these iterated conditionals does not properly allow to represent uncertainty  in conditional sentences of human reasoning and the study of uncertainty  is necessary for understanding human and artificial rationality in general.

Then, the  results obtained in this paper can be useful in AI in order to build a theory of formal reasoning which properly manage the uncertainty present in conditional or compound conditional sentences. 
Once it is described how  the conditionals and the logical operations among them are interpreted, by means of our results, it is possible, for instance, to  understand how a (coherent) agent propagates the uncertainty 
present in the conditionals $A|H$ and $B|K$ to the iterated conditional $(B|K)|_i(A|H)$.
Then,  it is crucial to know, for each iterated conditioning, which  (desirable) basic logical and probabilistic properties are preserved.
In particular,  our results allow to know that when the iterated conditioning  $|_i$ belongs to  $\{|_C,|_{dF},|_F,|_K,|_L\}$, as the two-premise centering is not p-valid, then  the agent   agrees with the following  probabilistically non-informative inference: from $P(A|H)=1$ and $P(B|K)=1$ infer  that every 
$P[(B|K)|_i(A|H)]\in [0,1]$ is coherent.

As already done for the iterated $|_{gs}$ in \cite{GiSa21ECSQARU}, future work will concern the (possible) characterization of the p-entailment of Adams  in the setting of coherence by means of the iterated conditioning $|_{K},|_{L},|_{B},$ and $|_{S}$. We 
will also study the different notions of iterated conditioning  in the framework of nonmonotonic reasoning in System P 
and  in  other non-classical logics, like connexive logic \cite{PfSa23}.  Finally, 
 as done in \cite{GiSa21ECSQARU} for the case of  conjoined conditionals,  in the more general theory of compound conditionals as suitable conditional random quantities (\cite{KR2022-12}), based on the  
 structure 
 $\Box|\bigcirc = \Box\wedge \bigcirc +\prev(\Box|\bigcirc)\no{\bigcirc}$,  future work  for this group of  research could   be devoted to the study of the (extended) iterated conditioning $|_{gs}$ to the case  where $\Box$ and $\bigcirc$ are compound conditionals.

\appendix
\section{}
\label{SEC:APP}
\subsection{Proof of  Theorem \ref{ThdeFinetti}.}
{\bf Theorem} \ref{ThdeFinetti}\label{SEC:APPdf1}:
\emph{Let $A$, $B$, $H$, $K$ be any logically independent events.
A probability assessment $\P=(x, y, z)$ on the family of conditional events $\F=\{A|H, (B|K)|_{dF}(A|H), (A|H)\wedge_{K}(B|K)\}$ is coherent if and only if 
$(x, y)\in[0,1]^2$ and $z \in [z', z'']$, where $z'=0$ and $z''=xy$.}
\begin{proof}
The constituents $C_h$'s and the point $Q_h$'s associated with the assessment $\P=(x,y,z)$ on $\F$ are
\[
\begin{array}{l}
C_1=AHBK, C_2=\no{A}HBK\vee \no{A}H\no{B}K\vee \no{A}H\no{K}=\no{A}H,\\ C_3=AH\no{B}K, C_4=\no{H}\no{B}K, C_5=AH\no{K}, C_0= \no{H}BK \vee \no{H}\,\no{K},
\end{array}
\]
and
\[
\begin{array}{l}
Q_1=(1,1,1), Q_2=(0,y,0), Q_3=(1,0,0),\\ Q_4=(x,y,0), Q_5=(1,y,z), \P=Q_0=(x,y,z).
\end{array}
\]
\begin{table}[h]
    \centering
    \begin{tabular}{|l|c|c c c| c|}
    \hline
        {} & {$C_h$} & {$A|H$} & {$(B|K)|_{dF} (A|H)$} & {$(A|H)\wedge_{K}(B|K)$} & {$Q_h$}\\
        \hline
        {$C_1$} & {$AHBK$} & {1} & {1} & {1} & {$Q_1$} \\
        {$C_2$} & {$\no{A}H$} & {0} & {$y$} & {0} & {$Q_2$} \\
        {$C_3$} & {$AH\no{B}K$} & {1} & {0} & {0} & {$Q_3$} \\
        {$C_4$} & {$\no{H}\no{B}K$} & {$x$} & {$y$} & {0} & {$Q_4$} \\
        {$C_5$} & {$AH\no{K}$} & {$1$} & {$y$} & {$z$} & {$Q_5$} \\
        {$C_0$} & {$\no{H}BK \vee \no{H}\,\no{K}$} & {$x$} & {$y$} & {$z$} & {$Q_0$} \\
    \hline
    \end{tabular}
  \caption{Constituents $C_h$ and points $Q_h$ associated with $\F=\{A|H, (B|K)|_{dF}(A|H), (A|H)\wedge_{K}(B|K)\}$, the probability assessment $\P=(x, y, z)$.}    
\end{table}  
The system $(\Sigma)$ in (\ref{SYST-SIGMA}) associated with the pair $(\F,\P)$  becomes\\
\begin{equation}
\left\{
    \begin{array}{llll}
     \lambda_1+\lambda_3+x\lambda_4+\lambda_5=x, \\
     \lambda_1+y\lambda_2+y\lambda_4+y\lambda_5=y,\\
     \lambda_1+z\lambda_5=z,\\
          \lambda_1+\lambda_2+\lambda_3+\lambda_4+\lambda_5=1, \\ \lambda_i\geq0 \ \forall i=1, \dots, 5.\\
    \end{array}
    \label{sis_df}
\right.    
\end{equation}
\subparagraph{Lower bound}
We  first prove that the assessment $(x, y, 0)$
is coherent for every 
$(x,y)\in[0,1]^2$.
We observe that $\P=(x,y,0)=Q_4$, so a solution of (\ref{sis_df}) is given by $\Lambda=(0,0,0,1,0)$.\\
Then, by considering the  function  $\phi$ as defined in (\ref{EQ:I0}),  it holds that 
\[
\phi_{1}(\Lambda)=\sum_{h:C_h\subseteq H}\lambda_h=\lambda_1+\lambda_2+\lambda_3+\lambda_5=0,
\]
\[
\phi_{2}(\Lambda)=\sum_{h:C_h\subseteq (AHK)}\lambda_h=\lambda_1+\lambda_3=0,
\]
\[
\phi_{3}(\Lambda)=\sum_{h:C_h\subseteq (AHBK\lor \bar{A}H \lor \bar{B}K)}\lambda_h=\lambda_1+\lambda_2+\lambda_3+\lambda_4=1>0.
\]

Let $\mathcal{S}'=\{(0,0,0,1,0)\}$  denote a subset of the set $\mathcal{S}$ of all solutions of (\ref{sis_df}). We have that $M_1'=0$, $M_2'=0$, $M_3'=1$  (as defined in (\ref{EQ:I0'})). 
Then $I_0'=\{1,2\}$ and we set $\K=\F_0'=\{A|B, (B|K)|_{dF}(A|H)\}$ and $\V=\P_0'=(x, y)$. 
The constituents $C_h$'s and the point $Q_h$'s associated with the assessment $\P_0$ on $\F_0$ are
\[
\begin{array}{l}
C_1=AHBK, C_2=\no{A}HBK\vee \no{A}H\no{B}K\vee \no{A}H\no{K}=\no{A}H,\\ C_3=AH\no{B}K,  C_4=AH\no{K}, C_0=\no{H}\,\no{K} \lor \no{H}BK \lor \no{H}\no{B}K=\no{H},
\end{array}
\]
and
\[
\begin{array}{l}
Q_1=(1,1),\, Q_2=(0,y),\, Q_3=(1,0),\, Q_4=(1,y), \P=Q_0=(x,y).
\end{array}
\]
The system $(\Sigma)$ in (\ref{SYST-SIGMA}) associated with the pair $(\K,\V)$  becomes\\
\begin{equation}
\left\{
    \begin{array}{llll}
     \lambda_1+\lambda_3+\lambda_4=x, \\
     \lambda_1+y\lambda_2+y\lambda_4=y,\\
    \lambda_1+\lambda_2+\lambda_3+\lambda_4+=1, \\ \lambda_i\geq0 \ \forall i=1, \dots, 4.\\
    \end{array}
    \label{sis_df2}
\right.    
\end{equation}
We observe that $(x,y)=(1-x)Q_2+xQ_4$ so a solution of (\ref{sis_df2} is $\Lambda=(0,1-x,0,x)$. \\
By considering the  function  $\phi$ as defined in (\ref{EQ:I0}),  it holds that 
\[
\phi_{1}(\Lambda)=\sum_{h:C_h\subseteq H}\lambda_h=\lambda_1+\lambda_2+\lambda_3+\lambda_4=1,
\]
\[
\phi_{2}(\Lambda)=\sum_{h:C_h\subseteq (AHK)}\lambda_h=\lambda_1+\lambda_3=0,
\]
Let $\mathcal{S}'=\{(0,1-x,0,x)\}$  denote a subset of the set $\mathcal{S}$ of all solutions of (\ref{sis_df2}). We have that $M_1'=0$, $M_2'=1$,  (as defined in (\ref{EQ:I0'})). 
Then, the set $I_0'$ associated with $(\K,\V)$ is $I_0'=\{2\}$. 
We observe that  the sub-assessment $y$ on $\{(B|K)|_{dF}(A|H)\}$ is coherent for every $y\in[0,1]$. Then, by Theorem \ref{CNES-PREV-I_0-INT}, the assessment $(x,y,0)$  on $\F$ is coherent $\forall (x,y)\in[0,1]^2$.
\subparagraph{Upper bound}
We  verify that the assessment  $(x,y,xy)$  on $\F$ is coherent  for every  $(x,y)\in[0,1]^2$. 
Moreover, we show that  $z''=xy$ is the upper bound for $z=P((A|H)\wedge_{K}(B|K))$ by showing that any assessment $(x,y,z)$ on $\F$ with $(x,y)\in[0,1]^2$ and  $z>xy$ is not coherent.\\
We observe that 
$$(x,y,xy)=xyQ_1+(1-x)Q_2+x(1-y)Q_3.$$
Then, the vector $\Lambda=(xy, 1-x, x(1-y), 0,0)$
is a solution of (\ref{sis_df}). Moreover, it holds that 
\[
\phi_{1}(\Lambda)=\sum_{h:C_h\subseteq H}\lambda_h=\lambda_1+\lambda_2+\lambda_3+\lambda_5=1>0,
\]
\[
\phi_{2}(\Lambda)=\sum_{h:C_h\subseteq (AHK)}\lambda_h=\lambda_1+\lambda_3=xy+1-x,
\]
\[
\phi_{3}(\Lambda)=\sum_{h:C_h\subseteq (AHBK\lor \bar{A}H \lor \bar{B}K)}\lambda_h=\lambda_1+\lambda_2+\lambda_3+\lambda_4=1>0.
\]
Let $\mathcal{S}'=\{(xy, 1-x, x(1-y), 0,0)\}$  denote a subset of the set $\mathcal{S}$ of all solutions of (\ref{sis_df}). We have that $M_1'=1$, $M_2'=xy+1-x$, $M_3'=1$  (as defined in (\ref{EQ:I0'})). 
We distinguish two cases: $(i)$ $(x\neq 1)\lor (y\neq0)$, $(ii)$ $(x=1)\wedge (y=0)$. In the case $(i)$ we get 
$M_1'>0$, $M_2'>0$, $M_3'>0$ and hence
$I'_0=\emptyset$. By Theorem \ref{CNES-PREV-I_0-INT}, the assessment $(x,y,xy)$ is coherent $\forall (x,y)\in[0,1]^2$. 
In the case $(ii)$ we get $M_1'>0$, $M_2'=0$, $M_3'>0$, then $I_0'=\{2\}$. We observe that  the sub-assessment $\P_0'y$ on $\F_0'\{(B|K)|_{dF}(A|H)\}$ is coherent for every $y\in[0,1]$. Then, by Theorem \ref{CNES-PREV-I_0-INT}, the assessment $(x,y,xy)$  on $\F$ is coherent $\forall (x,y)\in[0,1]^2$.
\\
We verify that $z''=xy$ is the upper bound for $z$, by showing that the assessment $(x,y,z)$ is incoherent when $z>xy$.\\ Let $z>xy$.  
We distinguish the following cases: $(a)$ $y\neq 0$; $(b)$ $y=0$.
\\
Case $(a)$. We observe that the points $Q_1, Q_2, Q_3$ belong to the plane $\pi: yX+Y-Z=y$, where $X, Y, Z$ are the axes coordinates. We set  $f(X,Y,Z)= yX+Y - Z$ and we  observe that  $f(\P)=f(x,y,z)=yx+y-z$. 
For each $h=1,2,3,4,5$, we compute the difference $f(Q_h)-f(\P)$. We have
\[
\begin{array}{ll}
f(Q_1)-f(\P)=f(Q_2)-f(\P)=
f(Q_3)-f(\P)=z-xy;\\
f(Q_4)-f(\P)=z;\;\;f(Q_5)-f(\P)=y-xy=y(1-x);\\
\end{array}
\]
We consider the sub-cases: $(i)$ $x<1$;  $(ii)$ $x=1$.\\
$(i)$ We recall that, by setting the stakes $s_1=y$, $s_2=1$, $s_3=-1$, it holds that  $g_h=f(Q_h)-f(\P)$, $h=1,\dots,5$, where  $g_h$ is the value of the random gain $G$ associated with the constituent $C_h\subseteq \H$. As  $x<1$, it follows that  $g_h=f(Q_h)-f(\P)>0$, $h=1,\dots,5$. Thus, as there exist $(s_1,s_2,s_3)$ such that 
 $\min G_{\H} \max G_{\H}>0$, it follows that the assessment $(x,y,z)$ is not coherent  when $z>xy$ and $x<1$.\\
$(ii). $ In this case, as $x=1$, it follows that  $\P=(x, y, z)=(1,y,z)=Q_5 \in \I$. Then, the vector $\Lambda=(0,0,0,0,1)$ is a solution of (\ref{sis_df}). Then, by considering the  function  $\phi$ as defined in (\ref{EQ:I0}),  it holds that 
\[
\phi_{1}(\Lambda)=\sum_{h:C_h\subseteq H}\lambda_h=\lambda_1+\lambda_2+\lambda_3+\lambda_5=1>0,
\]
\[
\phi_{2}(\Lambda)=\sum_{h:C_h\subseteq (AHK)}\lambda_h=\lambda_1+\lambda_3=0,
\]
\[
\phi_{3}(\Lambda)=\sum_{h:C_h\subseteq (AHBK\lor \bar{A}H \lor \bar{B}K)}\lambda_h=\lambda_1+\lambda_2+\lambda_3+\lambda_4=0.
\]

Then, we get $\I_0\subseteq \{2,3\}$.
and  we study the coherence of the sub-assessment $\P_0=(y, z)$ on 
 $\F_0=\{(B|K)|_{dF}(A|H), (A|H)\wedge_{K}(B|K)\}$.\\
The constituents $C_h$'s and the point $Q_h$'s associated with the pair  $(\P_0,\F_0)$ are
\[
\begin{array}{l}
C_1=AHBK,\ C_2=\no{A}H \lor \no{H}\no{B}K,\\ C_3=AH\no{B}K, \ C_0=\no{H}\,\no{K} \lor AH\no{K} \lor \no{H}BK,
\end{array}
\]
with $C_h\subseteq \H=HBK\vee \no{A}H\vee \no{B}K$, $h=1,2,3$. The associated points $Q_h$'s are
\[
\begin{array}{l}
Q_1=(1,1),\ Q_2=(y,0),\ Q_3=(0,0),\ \P=Q_0=(x,y,z).
\end{array}
\]
As shown in Figure \ref{FIG:fig_it_df_2d}  the convex hull $\I$ of the points $Q_1,Q_2$, and  $Q_3$ is the triangle of vertices $(1,1), (y,0)$, and $(0,0)$. We observe that $\P_0=(y,z) \notin \I$ because $z>y$.
 Then, the sub-assessment $\P_0$  on  $\F_0$ is not coherent. Therefore, by Theorem ... the assessment $\P=(1,y,z)$ on $\F$ is not coherent too.
 
 Notice that, if $y=0$, $(0,z)\in \I$ $\iff z=0$, so if $z>0$ the assessment is not coherent.
 
\begin{figure}[tbph]
\centering
\includegraphics[width=0.75\linewidth]{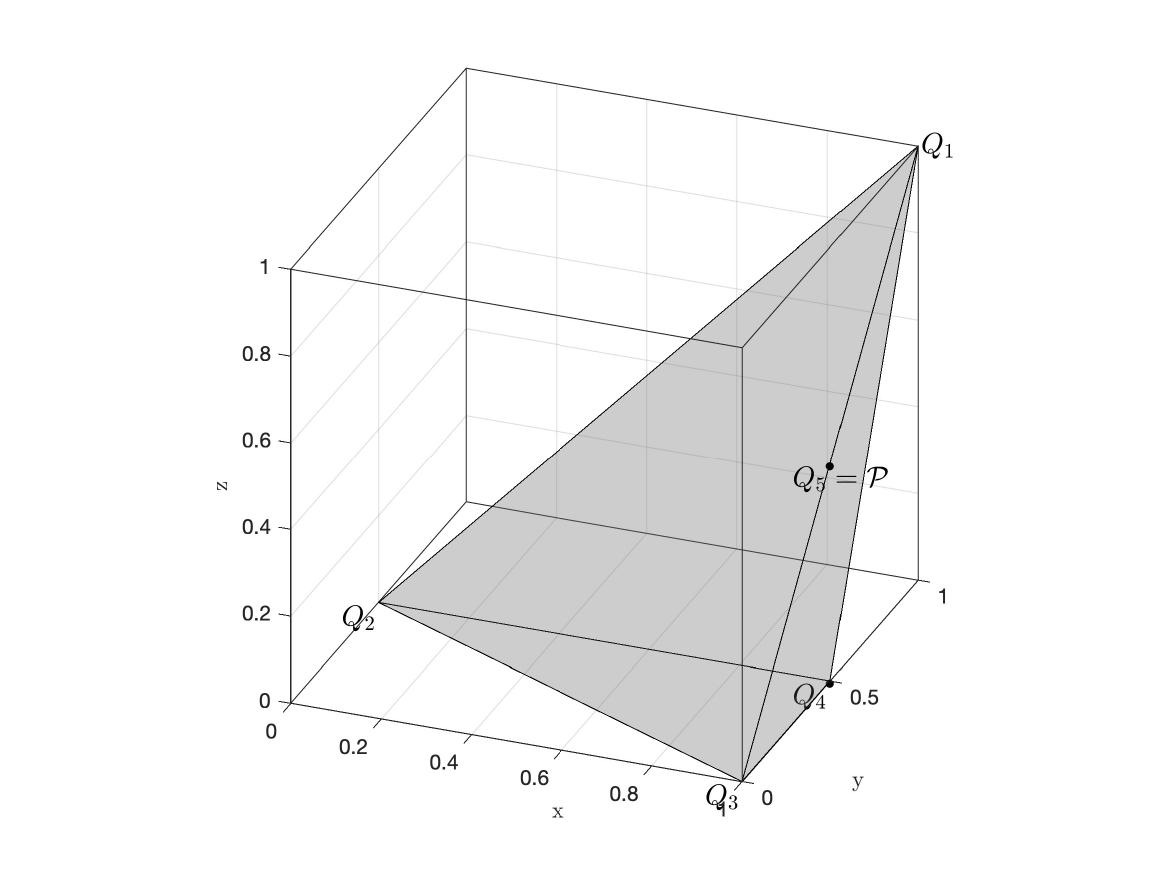}
\caption{Convex hull of  the points $Q_1, Q_2,Q_3, Q_4, Q_5$ associated with the pair $(\F,\P)$, where  $\F=\{A|H, (B|K)|_{dF} (A|H), (A|H)\wedge_{K}(B|K)\}$ and $\P=(x,y,z)$.  
In the figure the numerical  values are: $x=1$, $y=0.5$, $z=0.5$. 
Notice that $Q_5=\P$.}
\label{FIG:fig_it_df}
\end{figure}

\begin{figure}[tbph]
\centering
\includegraphics[width=0.45\linewidth]{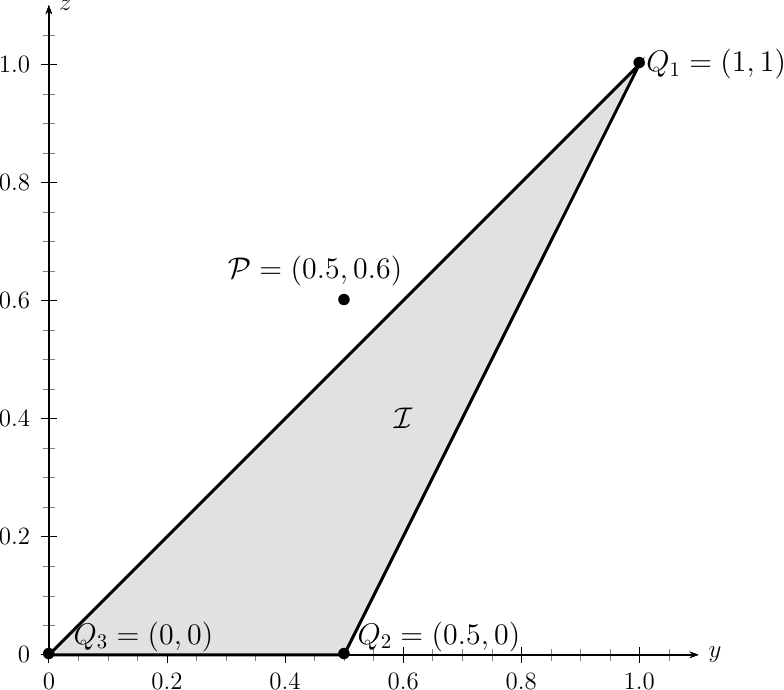}
\caption{Convex hull of  the points $Q_1, Q_2,Q_3$ associated with the pair $(\P_0,\F_0)$, where $\P_0=(y,z)$ and $\F_0=\{(B|K)|_{dF}(A|H), (A|H)\wedge_K (B|K)\}$.  
In the figure the numerical  values are: $y=0.5$ and  $z=0.6$.}
\label{FIG:fig_it_df_2d}
\end{figure}

So in this case we have that the probability assessment is coherent if and only if  $z \leq xy$.

Case $(b)$. In this case $y=0$. We have already analyze the situation $(x=1)\wedge(y=0)$ in the precedent case, so let's assume $x\neq1$.\\
 We observe that the points $Q_1, Q_2, Q_5$ belong to the plane $\pi: zX+(1-z)Y-Z=0$, where $X, Y, Z$ are the axes coordinates. We set  $f(X,Y,Z)=zX+(1-z)Y-Z$ and we  observe that  $f(\P)=f(x,y,z)=-(1-x)z$. 
For each $h=1,2,3,4,5$, we consider the quantity $f(Q_h)-f(P)$. Then, we obtain
\[
\begin{array}{ll}
f(Q_1)-f(\P)=z(1-x);\\
f(Q_2)-f(\P)=z(1-x);\\
f(Q_3)-f(\P)=z(2-x);\\
f(Q_4)-f(\P)=z;\\
f(Q_5)-f(\P)=z(1-x);\\
\end{array}
\]
We recall that, by setting the stakes $s_1=z$, $s_2=1-z$, $s_3=-1$, it holds that  $g_h=f(Q_h)-f(\P)$, $h=1,\dots,5$, where  $g_h$ is the value of the random gain $G$ associated with the constituent $C_h\subseteq \H$. As  $x<1$, it follows that  $g_h=f(Q_h)-f(\P)>0$, $h=1,\dots,5$. Thus, as there exist $(s_1,s_2,s_3)$ such that 
 $\min G_{\H} \max G_{\H}>0$, it follows that the assessment $(x,y,z)$ is not coherent  when $z>xy$ and $y=0$.\\
We conclude that $z''=xy$ is the upper bound for $z$.

 \end{proof}

\paragraph{Proof of  Theorem \ref{ThdeFinetti2}}\ \\
{\bf Theorem} \ref{ThdeFinetti2}\label{SEC:APPdf2}:
\emph{Let $A$, $B$, $H$, $K$ be any logically independent events.
The probability assessments  $\P=(x, y, z)$ on the family of conditional events $\F=\{A|H, B|K, (B|K)|_{dF}(A|H)\}$
is  coherent for every  $(x,y,z)\in[0,1]^3$.}
\begin{proof}
The constituents $C_h$'s and the points $Q_h$'s associated with the assessment $\P=(x,y,z)$ on $\F$ are (see also Table)
\[
\begin{array}{l}
C_1=AHBK, C_2=AH\no{B}K, C_3=AH\no{K}, C_4=\no{A}HBK, C_5=\no{A}H\no{B}K,\\ C_6=\no{A}H\no{K}, C_7=\no{H}BK, C_8=\no{H}\no{B}K, C_0=\no{H}\,\no{K},
\end{array}
\]
and
\[
\begin{array}{l}
Q_1=(1,1,1), Q_2=(1,0,0), Q_3=(1,y,z), Q_4=(0,1,z), Q_5=(0,0,z), \\
Q_6=(0,y,z), Q_7=(x,1,z), Q_8=(x,0,z), \P=Q_0=(x,y,z).
\end{array}
\]
\begin{table}[h]
    \centering
    \begin{tabular}{|l|c|c c c| c|}
    \hline
        {} & {$C_h$} & {$A|H$} & {$B|K$} & {$(B|K)|_{dF}(A|H)$} & {$Q_h$}\\
        \hline
        {$C_1$} & {$AHBK$} & {1} & {1} & {1} & {$Q_1$} \\
        {$C_2$} & {$AH\no{B}K$} & {1} & {0} & {0} & {$Q_2$} \\
        {$C_3$} & {$AH\no{K}$} & {1} & {$y$} & {$z$} & {$Q_3$} \\
        {$C_4$} & {$\no{A}HBK$} & {0} & {1} & {$z$} & {$Q_4$} \\
        {$C_5$} & {$\no{A}H\no{B}K$} & {0} & {0} & {$z$} & {$Q_5$} \\
        {$C_6$} & {$\no{A}H\no{K}$} & {0} & {$y$} & {$z$} & {$Q_6$} \\
        {$C_7$} & {$\no{H}BK$} & {$x$} & {$1$} & {$z$} & {$Q_7$} \\
        {$C_8$} & {$\no{H}\no{B}K$} & {$x$} & {$0$} & {$z$} & {$Q_8$} \\
        {$C_0$} & {$\no{H}\,\no{K}$} & {$x$} & {$y$} & {$z$} & {$Q_0$} \\
    \hline
    \end{tabular}
\caption{Constituents and points $Q_h$' associated with $\F=\{A|H, B|K, (B|K)|_{dF}(A|H)\}$, the probability assessment $\P=(x, y, z)$.}    
\end{table}
The system $(\Sigma)$ in (\ref{SYST-SIGMA}) associated with the pair $(\F,\P)$  becomes\\
\begin{equation}
\left\{
    \begin{array}{llll}
     \lambda_1+\lambda_2+\lambda_3+x\lambda_7+x\lambda_8=x, \\
     \lambda_1+y\lambda_3+\lambda_4+y\lambda_6+\lambda_7=y,\\
     \lambda_1+z\lambda_3+z\lambda_4+z\lambda_5+z\lambda_6+z\lambda_7+z\lambda_8=z,\\
          \lambda_1+\lambda_2+\lambda_3+\lambda_4+\lambda_5+\lambda_6+\lambda_7+\lambda_8=1, \\ \lambda_i\geq0 \ \forall i=1, \dots, 8.\\
    \end{array}
    \label{sdf}
\right.    
\end{equation}
We observe that $\P$ belongs to the segment with end points $Q_3$, $Q_6$; indeed $(x,y,z)=xQ_3+(1-x)Q_6=x(1,y,z)+(1-x)(0,y,z)$.\\
$\P$ also belongs to the segment $\no{Q_7Q_8}$, $(x,y,z)=yQ_7+(1-y)Q_8=y(x,1,z)+(1-y)(x,0,z)$.\\
Then 
$$(x,y,z)=\frac{x}{2}(1,y,z)+\frac{1-x}{2}(0,y,z)+\frac{y}{2}(x,1,z)+\frac{1-y}{2}(x,0,z),$$
so the vector $\Lambda=(0,0,\frac{x}{2},0,0,\frac{1-x}{2}, \frac{y}{2}, \frac{1-y}{2})$ is a solution of (\ref{sdf}), with
\[
\phi_{1}(\Lambda)=\sum_{h:C_h\subseteq H}\lambda_h=\lambda_1+\lambda_2+\lambda_3+\lambda_4+\lambda_5\lambda_6=\frac{1}{2}>0,
\]
\[
\phi_{2}(\Lambda)=\sum_{h:C_h\subseteq K}\lambda_h=\lambda_1+\lambda_2+\lambda_4+\lambda_5+\lambda_7+\lambda_8=\frac{1}{2},
\]
\[
\phi_{3}(\Lambda)=\sum_{h:C_h\subseteq AHK}\lambda_h=\lambda_1+\lambda_2=0.
\]
Let $\mathcal{S}'=\{(0,0,\frac{x}{2},0,0,\frac{1-x}{2}, \frac{y}{2}, \frac{1-y}{2})\}$  denote a subset of the set $\mathcal{S}$ of all solutions of (\ref{sdf}). We have that $M_1'=\frac{1}{2}$, $M_2'=\frac{1}{2}$, $M_3'=0$ (\ref{EQ:I0'}).
We get $\I'_0= \{3\}$. We observe that  the sub-assessment $z$ on $\{(B|K)|_{dF}(A|H)\}$ is coherent for every $z\in[0,1]$. Then, by Theorem \ref{CNES-PREV-I_0-INT}, the assessment $(x,y,z)$  on $\F$ is coherent $\forall (x,y,z)\in[0,1]^3$.
\end{proof}

\subsection{Proof of  Theorem \ref{ThFarrell}}

\label{SEC:APPF}
{\bf Theorem \ref{ThFarrell}}\label{SEC:APPF1}:
\emph{
    Let $A$, $B$, $H$, $K$, be any logically independent events. A probability assessment $\P=(x, y,z)$ on the family of conditional events $\F=\{A|H, (B|K)|_{F}(A|H), (A|H)\wedge_{K} (B|K)\}$ is coherent if and only if  $(x,y)\in[0,1]^2$ and $z\in[z', z'']$, where $z'=0$ and  $z''=T^{H}_0(x,y)$,
 where 
 \[
T^{H}_0(x,y)=\begin{cases}
                0, & \text{if } x=0 \text{ or } y=0\\
             \frac{xy}{x+y-xy}, &  \text{if } x\neq 0 \text{ and } y\neq 0,\\
         \end{cases}
 \]   
 is the Hamacher t-norm with parameter $\lambda=0$.
}
\begin{proof}
    The constituents $C_h$'s and the point $Q_h$'s associated with the assessment $\P=(x,y,z)$ on $\F$ are  (see Table \ref{TAB:F1})
\[
\begin{array}{l}
C_1=AHBK, \;C_2=AH\no{B}K, \;C_3=AH\no{K},C_4=\no{A}HBK\vee \no{A}H\no{B}K\vee \no{A}H\no{K}=\no{A}H,\; C_5=\no{H}\no{B}K, \;C_0= \no{H}BK \vee \no{H}\,\no{K},
\end{array}
\]
and
\[
\begin{array}{l}
Q_1=(1,1,1),\; Q_2=(1,0,0),\; Q_3=(1,y,z),\; Q_4=(0,y,0),\; Q_5=(x,0,0),\; \P=Q_0=(x,y,z).
\end{array}
\]
\begin{table}[h]
    \centering
    \begin{tabular}{|l|c|c c c| c|}
    \hline
        {} & {$C_h$} & {$A|H$} & {$(B|K)|_{F} (A|H)$} & {$(A|H)\wedge_{K}(B|K)$} & {$Q_h$}\\
        \hline
        {$C_1$} & {$AHBK$} & {1} & {1} & {1} & {$Q_1$} \\
        {$C_2$} & {$AH\no{B}K$} & {$1$} & {$0$} & {$0$} & {$Q_2$} \\
        {$C_3$} & {$AH\no{K}$} & {$1$} & {$y$} & {$z$} & {$Q_3$} \\
        {$C_4$} & {$\no{A}H$} & {$0$} & {$y$} & {$0$} & {$Q_4$} \\
        {$C_5$} & {$\no{H}\no{B}K$} & {$x$} & {$0$} & {$0$} & {$Q_5$} \\
        {$C_0$} & {$\no{H}BK \vee \no{H}\,\no{K}$} & {$x$} & {$y$} & {$z$} & {$Q_0$} \\
    \hline
    \end{tabular}
  \caption{Constituents  and points $Q_h$'s associated with $\F=\{A|H, (B|K)|_{F}(A|H), (A|H)\wedge_{K}(B|K)\}$ and  $\P=(x, y, z)$.}    
\label{TAB:F1}
\end{table}
We observe that $\H_3=C_1\vee \cdots \vee C_5=H \vee \no{H}\no{B}K$.
The system $(\Sigma)$ in (\ref{SYST-SIGMA}) associated with the pair $(\F,\P)$  becomes\\
\begin{equation}
\left\{
    \begin{array}{llll}
     \lambda_1+\lambda_2+\lambda_3+x\lambda_5=x, \\
     \lambda_1+y\lambda_3+y\lambda_4=y,\\
     \lambda_1+z\lambda_3=z,\\
          \lambda_1+\lambda_2+\lambda_3+\lambda_4+\lambda_5=1, \\ \lambda_i\geq0 \ \forall i=1, \dots, 5.\\
    \end{array}
    \label{sis_F}
\right.    
\end{equation}
\subparagraph{Lower bound}
We  first prove that the assessment $(x, y, 0)$ is coherent for every $(x,y)\in[0,1]^2$.
We observe that $\P=(x,y,0)=xQ_3+(1-x)Q_4$, so a solution of (\ref{sis_F}) is given by $\Lambda=(0,0,x,1-x,0)$.\\
Then, by considering the  function  $\phi$ as defined in (\ref{EQ:I0}),  it holds that 
\[
\phi_{1}(\Lambda)=\sum_{h:C_h\subseteq H}\lambda_h=\lambda_1+\lambda_2+\lambda_3+\lambda_4=1>0,
\]
\[
\phi_{2}(\Lambda)=\sum_{h:C_h\subseteq (AHBK \vee AH\no{B}K \vee \no{H}\no{B}K)}\lambda_h=\lambda_1+\lambda_2+\lambda_5=0,
\]
\[
\phi_{3}(\Lambda)=\sum_{h:C_h\subseteq (AHBK\lor \bar{A}H \lor \bar{B}K)}\lambda_h=\lambda_1+\lambda_2+\lambda_4+\lambda_5=1-x.
\]
Let $\mathcal{S}'=\{(0,0,x,1-x,0)\}$  denote a subset of the set $\mathcal{S}$ of all solutions of (\ref{sis_F}). 
We have that $M_1'=1$, $M_2'=0$, $M_3'=1-x$ (as defined in (\ref{EQ:I0'})). 
We distinguish two cases: $(i)$ $x\neq 1$; $(ii)$ $x=1$.
\begin{itemize}
    \item [$(i)$]  
    In this case,  $M_1'>0$, $M_2'=0$, $M_3'>0$ and hence  $I_0'=\{2\}$.
    We observe that  the sub-assessment $\P_0'=y$ on $\F_0'=\{(B|K)|_{F}(A|H)\}$ is coherent for every $y\in[0,1]$. Then, by Theorem \ref{CNES-PREV-I_0-INT'}, the assessment $(x,y,0)$  on $\F$ is coherent $\forall (x,y)\in[0,1]^2$;
    \item[$(ii)$]
   We have that  $M_1'>0$, $M_2'=M_3'=0$ and hence  $I_0'=\{2,3\}$.
   We set $\K=\F_0'=\{(B|K)|_{F}(A|H), (A|H)\wedge_{K}(B|K)\}$ and $\V=\P_0'=(y,0)$. 
   By Theorem \ref{CNES-PREV-I_0-INT'},
 as (\ref{sis_F}) is solvable and $I_0'=\{2,3\}$,
  it is sufficient to check the coherence of 
 the sub-assessment  $\V$ on $\K$ in order to check the coherence of $(x,y,0)$.
   The constituents $C_h$'s and the point $Q_h$'s associated with the assessment $\V$ on $\K$ are
\[
\begin{array}{l}
C_1=AHBK, \;C_2=AH\no{B}K \vee \no{H}\no{B}K, \;\;C_3=\no{A}H, C_0=AH\no{K} \vee \no{H}\,\no{K} \lor \no{H}BK,
\end{array}
\]
and
\[
\begin{array}{l}
Q_1=(1,1),\, Q_2=(0,0),\, Q_3=(y,0),\, \P_0'=Q_0=(y,0).
\end{array}
\]
We have that $\H_2=C_1\vee \cdots \vee C_3=AHBK\vee \no{A}H\vee \no{B}K$.
The system $(\Sigma)$ in (\ref{SYST-SIGMA}) associated with the pair $(\K,\V)$  becomes\\
\begin{equation}
\left\{
    \begin{array}{llll}
     \lambda_1+y\lambda_3=y, \\
     \lambda_1=0,\\
    \lambda_1+\lambda_2+\lambda_3=1, \\ \lambda_i\geq0 \ \forall i=1, \dots, 3.\\
    \end{array}
    \label{sis_F2}
\right.    
\end{equation}
We observe that $\V=(y,0)=Q_3$ so a solution of (\ref{sis_F2} is $\Lambda=(0,0,1)$. \\
By considering the  function  $\phi$ as defined in (\ref{EQ:I0}),  it holds that 
\[
\phi_{1}(\Lambda)=\sum_{h:C_h\subseteq (AHBK \vee AH\no{B}K \vee \no{H}\no{B}K)}\lambda_h=\lambda_1+\lambda_2=0,
\]
\[
\phi_{2}(\Lambda)=\sum_{h:C_h\subseteq (AHBK \vee \no{A}H\vee \no{B}K)}\lambda_h=\lambda_1+\lambda_2+\lambda_3=1.
\]
Let $\mathcal{S}'=\{(0,0,1)\}$  denote a subset of the set $\mathcal{S}$ of all solutions of (\ref{sis_F}). We have that $M_1'=0$, $M_2'=1$,  (as defined in (\ref{EQ:I0'})). 
Then, the set $I_0'$ associated with $(\K,\V)$ is $I_0'={1}$. We observe that  the sub-assessment $y$ on $\{(B|K)|_{dF}(A|H)\}$ is coherent for every $y\in[0,1]$.  
Then, by Theorem \ref{CNES-PREV-I_0-INT'}, 
the sub-assessment $(y,0)$ on $\K$ is coherent for every $y\in[0,1]$ and hence 
the assessment $(x,y,0)$  on $\F$ is coherent $\forall (x,y)\in[0,1]^2$. Thus, $z'=0$ is the lower bound of $z=P((A|H)\wedge_K(B|K))$, for every $(x,y)\in[0,1]$ on $\{A|H,(B|K)|_F(A|H)\}$.
\end{itemize}

\subparagraph{Upper bound}
To study the upper bound $z''$ of $z=P((A|H)\wedge_K(B|K))$ we distinguish the following  cases: $(a)$  $x=0$; $(b)$  $y=0$; 
$(c)$ $x\neq 0$ and $y\neq 0$.
\\ 
 $(a)$. In this case system (\ref{sis_F}) associated to the pair  $(\F,(0,y,z))$ becomes
    \begin{equation}
\left\{
    \begin{array}{llll}
     \lambda_1+\lambda_2+\lambda_3=0, \\
     \lambda_1+y\lambda_3+y\lambda_4=y,\\
     \lambda_1+z\lambda_3=z,\\
    \lambda_1+\lambda_2+\lambda_3+\lambda_4+\lambda_5=1, \\ \lambda_i\geq0 \ \forall i=1, \dots, 5.\\
    \end{array}
\right.    
\Longleftrightarrow    
\left\{
    \begin{array}{llll}
     \lambda_1=\lambda_2=\lambda_3=0, \\
     y\lambda_4=y
     \\
z=0,\\  \lambda_4+\lambda_5=1, \\ \lambda_i\geq0 \ \forall i=4, 5,\\
    \end{array}
\right.    
\label{sis_F3}
\end{equation}
which  is solvable  only if $z=0$. We 
have already  verified  that $\P=(x,y,0)$ is coherent on $\F$ for all $(x,y)\in [0,1]^2$, hence, in particular, $(0,y,0)$ is coherent.  We observe that  system (\ref{sis_F3}) with $z>0$ is not solvable, then the assessment $(0,y,z)$, with $z>0$ is not coherent. Thus,  $z''=0$ is the upper bound for $z$ when $x=0$ and  $y\in[0,1]$.
\\ \ \\
$(b)$. 
 We recall that the pair $(\wedge_K, |_F)$ satisfies property P2, that is
$z\leq y$. As $y=0$, it follows that $z=0$. We 
have already  verified  that $\P=(x,y,0)$ is coherent on $\F$ for all $(x,y)\in [0,1]^2$, hence, in particular, $(x,0,0)$ is coherent.  We also observe that   the assessment $(x,0,z)$, with $z>0$ is not coherent,
 because the pair $(y=0,z)$ violates property P2.  Thus,  $z''=0$ is the upper bound for $z$ when $x\in[0,1]$ and $y=0$.
\\ \ \\
(c) Fist of all, we verify that the assessment $\P=(x,y,\frac{xy}{x+y-xy})$ on $\F$ is coherent. Then, we verify that $z''=\frac{xy}{x+y-xy}$ is the upper bound for $z$, by showing  that $(x,y,z)$, with $z>z''$
, is incoherent for every $(x,y)\in[0,1]^2$. \\ We observe that
\[
\P=(x,y,\frac{xy}{x+y-xy})=\frac{xy}{x+y-xy}Q_1+\frac{y(1-x)}{x+y-xy}Q_4+\frac{x(1-y)}{x+y-xy}Q_5.
\]
 Then a solution of (\ref{sis_F}) is $\Lambda=(\frac{xy}{x+y-xy}, 0,0,\frac{y(1-x)}{x+y-xy},\frac{x(1-y)}{x+y-xy})$ and it holds that 
 \[
\phi_{1}(\Lambda)=\sum_{h:C_h\subseteq H}\lambda_h=\lambda_1+\lambda_2+\lambda_3+\lambda_4=\frac{y}{x+y-xy},
\]
\[
\phi_{2}(\Lambda)=\sum_{h:C_h\subseteq (AHBK \vee AH\no{B}K \vee \no{H}\no{B}K)}\lambda_h=\lambda_1+\lambda_2+\lambda_5=\frac{x}{x+y-xy},
\]
\[
\phi_{3}(\Lambda)=\sum_{h:C_h\subseteq (AHBK\lor \bar{A}H \lor \bar{B}K)}\lambda_h=\lambda_1+\lambda_2+\lambda_4+\lambda_5=1>0.
\]
Let $\mathcal{S}'=\{(\frac{xy}{x+y-xy}, 0,0,\frac{y(1-x)}{x+y-xy},\frac{x(1-y)}{x+y-xy})\}$  be a subset $\mathcal{S}$ of all solutions of (\ref{sis_F}). We have that $M_1'>0$, $M_2'>0$, $M_3'>0$  (as defined in (\ref{EQ:I0'})).   Then, by Theorem\ref{CNES-PREV-I_0-INT'} the assessment $(x,y,\frac{xy}{x+y-xy})$ is coherent for every $x\neq 0$ and $y\neq 0$.\\
Now we  prove that in this case $z''=\frac{xy}{x+y-xy}$ is the upper bound for $z$.\\
We distinguish two sub-cases: $(i)$ $x=1$, $(ii)$ $x\neq 1$.\\
$(i)$ We recall that the pair $(\wedge_K, |_F)$ satisfies property P2 and hence  $z\leq y$. Then, when $x=1$, as $\frac{xy}{x+y-xy}=y$, we have that any assessment $z>\frac{xy}{x+y-xy}=y$ is not coherent. Therefore, $z''=\frac{xy}{x+y-xy}=y$ is the upper bound for $z=P((A|H)\wedge_{K}(B|K))$.

$(ii)$ We recall that in this case $x\in (0,1)$ and $y\neq 0$. We observe that the points $Q_1, Q_4, Q_5$ belong to the plane $\pi: yX+xY-(x+y-xy)Z=xy$, where $X, Y, Z$ are the axes coordinates. We set  $f(X,Y,Z)=yX+xY-(x+y-xy)Z-xy$, we
choose $\P=(x,y,z)$ with $z>\frac{xy}{x+y-xy}$, i.e. $(x+y-xy)z-xy>0$, and we observe that  $f(\P)=xy-(x+y-xy)z$. 
For each $h=1,\ldots, 5$, we consider the quantity $f(Q_h)-f(P)$. Then, we obtain
\begin{equation}\label{EQ:ghFarrell}
\begin{array}{ll}
f(Q_1)-f(\P)=(x+y-xy)z-xy>0;\\
f(Q_2)-f(\P)=y(1-x)+(x+y-xy)z-xy>0;\\
f(Q_3)-f(\P)=y-xy=y(1-x)>0;\\
f(Q_4)-f(\P)=(x+y-xy)z-xy>0;\\
f(Q_5)-f(\P)=(x+y-xy)z-xy>0.\\
\end{array}
\end{equation}
We recall that, by setting the stakes $s_1=y$, $s_2=x$, $s_3=-x-y+xy$, it holds that  $g_h=f(Q_h)-f(\P)$, $h=1,\dots,5$, where  $g_h$ is the value of the random gain $G$ associated with the constituent $C_h\subseteq \H_3$. 
From (\ref{EQ:ghFarrell}) we obtain that 
$g_h=f(Q_h)-f(\P)>0$, $h=1,\dots,5$. Thus, as there exist $(s_1,s_2,s_3)$ such that 
 $\min G_{\H} \max G_{\H}>0$, it follows that the assessment $(x,y,z)$ is not coherent  when $z>xy$.\\
We conclude that $z''=\frac{xy}{x+y-xy}$ is the upper bound for $z$. \\ 
Therefore $z''=T_0^H(x,y)$ is the upper bound for $z$, for  every $(x,y)\in[0,1]^2$.
\end{proof}

\subsection{Proof of  Theorem \ref{ThFarrell2}}
{\bf Theorem} \ref{ThFarrell2}\label{SEC:APPF2}:
\emph{
Let $A$, $B$, $H$, $K$ be any logically independent events.
The probability assessments  $\P=(x, y, z)$ on the family of conditional events $\F=\{A|H, B|K, (B|K)|_{F}(A|H)\}$
is  coherent for every  $(x,y,z)\in[0,1]^3$.
}
\begin{proof}
The constituents $C_h$'s and the points $Q_h$'s associated with the assessment $\P=(x,y,z)$ on $\F$ are (see also Table\ref{tab:FP4})
\[
\begin{array}{l}
C_1=AHBK, C_2=AH\no{B}K, C_3=AH\no{K}, C_4=\no{A}HBK, C_5=\no{A}H\no{B}K,\\ C_6=\no{A}H\no{K}, C_7=\no{H}BK, C_8=\no{H}\no{B}K, C_0=\no{H}\,\no{K},
\end{array}
\]
and
\[
\begin{array}{l}
Q_1=(1,1,1), Q_2=(1,0,0), Q_3=(1,y,z), Q_4=(0,1,z), Q_5=(0,0,z), \\
Q_6=(0,y,z), Q_7=(x,1,z), Q_8=(x,0,0), \P=Q_0=(x,y,z).
\end{array}
\]
\begin{table}[h]
    \centering
    \begin{tabular}{|l|c|c c c| c|}
    \hline
        {} & {$C_h$} & {$A|H$} & {$B|K$} & {$(B|K)|_{F}(A|H)$} & {$Q_h$}\\
        \hline
        {$C_1$} & {$AHBK$} & {1} & {1} & {1} & {$Q_1$} \\
        {$C_2$} & {$AH\no{B}K$} & {1} & {0} & {0} & {$Q_2$} \\
        {$C_3$} & {$AH\no{K}$} & {1} & {$y$} & {$z$} & {$Q_3$} \\
        {$C_4$} & {$\no{A}HBK$} & {0} & {1} & {$z$} & {$Q_4$} \\
        {$C_5$} & {$\no{A}H\no{B}K$} & {0} & {0} & {$z$} & {$Q_5$} \\
        {$C_6$} & {$\no{A}H\no{K}$} & {0} & {$y$} & {$z$} & {$Q_6$} \\
        {$C_7$} & {$\no{H}BK$} & {$x$} & {$1$} & {$z$} & {$Q_7$} \\
        {$C_8$} & {$\no{H}\no{B}K$} & {$x$} & {$0$} & {$0$} & {$Q_8$} \\
        {$C_0$} & {$\no{H}\,\no{K}$} & {$x$} & {$y$} & {$z$} & {$Q_0$} \\
    \hline
    \end{tabular}
\caption{Constituents and points $Q_h$'s associated with $\F=\{A|H, B|K, (B|K)|_{F}(A|H)\}$ and  $\P=(x, y, z)$.}    
    \label{tab:FP4}
\end{table}
We observe that $C_1\vee \cdots \vee C_8=H\vee K=\H_3$.
The system $(\Sigma)$ in (\ref{SYST-SIGMA}) associated with the pair $(\F,\P)$  becomes\\
\begin{equation}
\left\{
    \begin{array}{llll}
     \lambda_1+\lambda_2+\lambda_3+x\lambda_7+x\lambda_8=x, \\
     \lambda_1+y\lambda_3+\lambda_4+y\lambda_6+\lambda_7=y,\\
     \lambda_1+z\lambda_3+z\lambda_4+z\lambda_5+z\lambda_6+z\lambda_7=z,\\
          \lambda_1+\lambda_2+\lambda_3+\lambda_4+\lambda_5+\lambda_6+\lambda_7+\lambda_8=1, \\ \lambda_i\geq0 \ \forall i=1, \dots, 8.\\
    \end{array}
    \label{sF}
\right.    
\end{equation}
We observe that $\P$ belongs to the segment with end points $Q_3$, $Q_6$; indeed $(x,y,z)=xQ_3+(1-x)Q_6=x(1,y,z)+(1-x)(0,y,z)$. Then, the vector $\Lambda=(0,0,x,0,0,1-x,0,0)$ is a solution of (\ref{sF}), with
\[
\phi_{1}(\Lambda)=\sum_{h:C_h\subseteq H}\lambda_h=\lambda_1+\lambda_2+\lambda_3+\lambda_4+\lambda_5+\lambda_6=1,
\]
\[
\phi_{2}(\Lambda)=\sum_{h:C_h\subseteq K}\lambda_h=\lambda_1+\lambda_2+\lambda_4+\lambda_5+\lambda_7+\lambda_8=0,
\]
\[
\phi_{3}(\Lambda)=\sum_{h:C_h\subseteq [AHBK \vee AH\no{B}K \vee \no{H}\no{B}K]}\lambda_h=\lambda_1+\lambda_2+\lambda_8=0.
\]
Let $\mathcal{S}'=\{(0,0,x,0,0,1-x,0,0)\}$  be a subset $\mathcal{S}$ of all solutions of (\ref{sF}). We have that $M_1'>0$, $M_2'=0$, $M_3'=0$  (as defined in (\ref{EQ:I0'})).   Then,  $I_0'=\{2,3\}$. We check the coherence of the sub-assessment $\P_0'=(y,z)$ on $\F_0'=\{B|K, (B|K)|_{F}(A|H)\}$.
The constituents $C_h$'s and the point $Q_h$'s associated with the assessment $\P_0'$ on $\F_0'$ are
\[
\begin{array}{l}
C_1=AHBK, C_2=AH\no{B}K \vee \no{H}\no{B}K,\\ C_3=\no{A}HBK\vee \no{H}BK,  C_4=\no{A}H\no{B}K, C_0=\no{K},
\end{array}
\]
and
\[
\begin{array}{l}
Q_1=(1,1),\, Q_2=(0,0),\, Q_3=(1,z),\, Q_4=(0,z), \P=Q_0=(y,z).
\end{array}
\]
We observe that $C_1\vee \cdots \vee C_4=K=\H_2$.
The system $(\Sigma)$ in (\ref{SYST-SIGMA}) associated with the pair $(\F_0,\P_0)$  becomes\\
\begin{equation}
\left\{
    \begin{array}{llll}
     \lambda_1+\lambda_3=y, \\
     \lambda_1+z\lambda_3+z\lambda_4=z,\\
    \lambda_1+\lambda_2+\lambda_3+\lambda_4=1, \\ \lambda_i\geq0 \ \forall i=1, \dots, 4.\\
    \end{array}
    \label{sF2}
\right.    
\end{equation}
We observe that $(y,z)=yQ_3+(1-y)Q_4$ so a solution of (\ref{sF2} is $\Lambda_0=(0,0,y,1-y)$.
By considering the  function  $\phi'$ it holds that
\[
\phi_{1}(\Lambda)=\sum_{h:C_h\subseteq K}\lambda_h=\lambda_1+\lambda_2+\lambda_3+\lambda_4=1,
\]
\[
\phi_{2}(\Lambda)=\sum_{h:C_h\subseteq [AHBK \vee AH\no{B}K \vee \no{H}\no{B}K]}\lambda_h=\lambda_1+\lambda_2=0.
\] 
Let $\mathcal{S}'=\{(0,0,y,1-y)\}$  be a subset $\mathcal{S}$ of all solutions of (\ref{sF2}). We have that $M_1'>0$, $M_2'=0$  (as defined in (\ref{EQ:I0'})).   Then,  $I_0'=\{2\}$.  We observe that  the sub-assessment $z$ on $\{(B|K)|_{F}(A|H)\}$ is coherent for every $z\in[0,1]$. Then, by Theorem \ref{CNES-PREV-I_0-INT'}, the assessment $(x,y,z)$  on $\F$ is coherent $\forall (x,y,z)\in[0,1]^3$.
\end{proof}

\subsection{Proof of  Theorem \ref{THM:PIL}}
\label{SEC:APPL} 
 {\bf Theorem} \ref{THM:PIL}
\emph{Let $A$, $B$, $H$, $K$ be any logically independent events. The set $\Pi$ of all the coherent assessment $(x,y,z,\mu)$ on the family $\F=\{A|H, B|K, (A|H)\wedge_{L}(B|K), (B|K)|_{L}(A|H)\}$ is $\Pi=\Pi' \cup \Pi''$, where 
$\Pi'=\{(x,y,z,\mu): x\in(0,1], y\in [0,1], z\in[z', z''], \mu=\frac{z}{x}\}$
 with $ z'=0$,    $z''=min\{x,y\}$, and
$\Pi''=\{(0,y,0,\mu): (y,\mu)\in[0,1]^2\}.$}
\begin{proof}
It is well-known that the assessment $(x,y)$ on $\{A|H, B|K\}$ is coherent for every $(x,y)\in [0,1]^2$. By Table (\ref{tab:LUB}), the assessment $z=P((A|H)\wedge_{L}(B|K))$ is a coherent extension of $(x,y)$ if and only if $z \in [z',z'']$ where $z'=0$ and $z''=min\{x,y\}$. Assuming $x>0$, from (\ref{prevision}) it holds that $\mu=\frac{z}{x}$. Then, every $(x,y,z,\mu)\in \Pi'$ is coherent, i.e., $\Pi'\subseteq \Pi$. Of course, if $x>0$ and $(x,y,z,\mu)\notin \Pi'$, then the assessment $(x,y,z,\mu)$ is not coherent, i.e. $(x,y,z,\mu)\notin \Pi$.\\
  Let us consider now the case $x=0$, so that $z'=0$ and $z''=0$. We show that the assessment $(0,y,0,\mu)$ is coherent if and only if  $(y,\mu)\in [0,1]^2$, that is $\Pi''\subseteq \Pi$.
As $x=0$, from (\ref{EQ:iteL'}), it holds that 
\begin{equation}\label{EQ:iteL'x=0}
\begin{array}{l}
(B|K)|_{L}(A|H) 
= \begin{cases}
    1,               & {\ \text{if} \ AHBK \ \text{is true}},\\
    0,               & {\ \text{if} \ AH\no{B}K \vee AH\no{K} \ \text{is true}},\\
    \mu_L,             & {\ \text{if} \ \no{A}H\vee \no{H}BK \vee \no{H}\no{B}K\lor \no{H}\no{K}\ \text{is true}}.
    \end{cases}
\end{array}
\end{equation}
that is $(B|K)|_{L}(A|H)=BK|AH$ (see proof of Theorem \ref{THM:PIK}).  Therefore $\F=\{A|H, B|K, (A|H)\wedge_{L}(B|K), (B|K)|_{L}(A|H)\}=\{A|H, B|K, AHBK|(AHBK\vee \overline{A}H\vee \overline{B}K \vee \overline{H}\,\overline{K}), BK|AH\}$.
The constituents $C_h's$ and the points $Q_h's$ associated with $(\F, \P)$, where $\P=(0,y,0,\mu)$ are the following (see also Table \ref{tabxL}):
 \[
 \begin{array}{l}
 C_1=AHBK, C_2=AH\no{B}K, C_3=AH\no{K}, C_4=\no{A}HBK, C_5=\no{A}H\no{B}K,\\ C_6=\no{A}H\no{K}, C_7=\no{H}BK, C_8=\no{H}\no{B}K, C_9=\no{H}\,\no{K},
 \end{array}
 \]
 and
 \[
 \begin{array}{l}
 Q_1=(1,1,1,1), Q_2=(1,0,0,0), Q_3=(1,y,0,0), Q_4=(0,1,0,\mu), 
 Q_5=(0,0,0,\mu),\\ Q_6=(0,y,0,\mu),Q_7=(0,1,0,\mu), Q_8=(0,0,0,\mu), 
 Q_9=(0,y,0,\mu).
 \end{array}
 \]
 \begin{table}[h]
     \centering
     \begin{tabular}{|l|c|c c c c| c|}
     \hline
         {} & {$C_h$} & {$A|H$} & {$B|K$} & {$(A|H)\wedge_{L}(B|K)$} &{$(B|K)|_{L}(A|H)$} & {$Q_h$}\\
         \hline
         {$C_1$} & {$AHBK$} & {1} & {1} & {1} & {1} &{$Q_1$} \\
         {$C_2$} & {$AH\no{B}K$} & {1} & {0} & {0}& {0} & {$Q_2$} \\
         {$C_3$} & {$AH\no{K}$} & {1} & {$y$} & {0}& {0} & {$Q_3$} \\
         {$C_4$} & {$\no{A}HBK$} & {0} & {1} & {0} & {$\mu$} & {$Q_4$} \\
         {$C_5$} & {$\no{A}H\no{B}K $} & {0} & {0} & {0} & {$\mu$} & {$Q_5$} \\
         {$C_6$} & {$\no{A}H\no{K}$} & {0} & {$y$} & {0} & {$\mu$} & {$Q_6$} \\
         {$C_7$} & {$\no{H}BK$} & {0} & {1} & {0} & {$\mu$} & {$Q_7$} \\
         {$C_8$} & {$\no{H}\no{B}K $} & {0} & {0} & {0} & {$\mu$} & {$Q_8$} \\
         {$C_9$} & {$\no{H}\,\no{K}$} & {0} & {$y$} & {0} & {$\mu$} & {$Q_9$} \\
     \hline
     \end{tabular}
 \caption{Constituents and points $Q_h$' associated with $\F = \{A|H, B|K, (A|H)\wedge_{L}(B|K), (B|K)|_{L}(A|H)\}$ and  $\P=(0, y, 0, \mu)$.}  
 \label{tabxL}
 \end{table}
We notice that $C_1\vee\cdots \vee  C_9=\Omega.$
The system $(\Sigma)$ in (\ref{SYST-SIGMA}) associated with the pair $(\F,\P)$  is\\
 \begin{equation}
 \left\{
    \begin{array}{llll}
       \lambda_1+\lambda_2+\lambda_3=0, \\
       \lambda_1+y\lambda_3+\lambda_4+y\lambda_6+\lambda_7+y\lambda_9=y,\\
      \lambda_1=0,\\
      \lambda_1+\mu \lambda_4+\mu \lambda_5+\mu \lambda_6+\mu \lambda_7+\mu \lambda_8+\mu \lambda_9=\mu\\
      \lambda_1+\lambda_2+\lambda_3+\lambda_4+\lambda_5+\lambda_6+\lambda_7+\lambda_8+\lambda_9=1, \\ \lambda_i\geq0 \ \forall i=1, \dots, 9.\\
     \end{array}
     \label{sL}
 \right.    
 \end{equation}
 We observe that $\P$ belongs to the segment with end points $Q_4$, $Q_5$; indeed $(0,y,0,\mu)=yQ_4+(1-y)Q_5=y(0,1,0,\mu)+(1-y)(0,0,0,\mu)$.
 Then the vector $\Lambda=(0,0,0,y,1-y,0,0,0,0)$ is a solution of (\ref{sL}), with
\[
 \phi_{1}(\Lambda)=\sum_{h:C_h\subseteq H}\lambda_h=\lambda_1+\lambda_2+\lambda_3+\lambda_4+\lambda_5=1>0,
 \]
 \[
\phi_{2}(\Lambda)=\sum_{h:C_h\subseteq K}\lambda_h=\lambda_1+\lambda_2+\lambda_4+\lambda_5=1,
 \]
 \[
 \phi_{3}(\Lambda)=\sum_{h:C_h\subseteq (AHBK\lor\no{A}H \lor \no{B}K \lor \no{H}\no{K})}\lambda_h=\lambda_1+\lambda_2+\lambda_4+\lambda_5+\lambda_6+\lambda_8+\lambda_9=1.
 \]
 \[
 \phi_{4}(\Lambda)=\sum_{h:C_h\subseteq AH}\lambda_h=\lambda_1+\lambda_2+\lambda_3=0.
 \]
Let $\mathcal{S}'=\{(0,0,0,y,1-y,0,0,0,0)\}$  denote a subset of the set $\mathcal{S}$ of all solutions of (\ref{sL}). We have that $M_1'=M_2'=M_3'=1$ and $M_4'=0$   (as defined in (\ref{EQ:I0'})). 
It follows that $I'_0=\{4\}$. As the sub-assessment $\P'_0=\mu$ on $\F'_0=\{BK|AH\}$ is coherent $\forall \mu \in [0,1]$,  by Theorem \ref{CNES-PREV-I_0-INT}, it follows that 
the assessment  $(0,y,0,\mu)$ on $\F=\{A|H, B|K, (A|H)\wedge_{L}(B|K), (B|K)|_{L}(A|H)\}$ is coherent for every $(y,\mu)\in[0,1]^2$, that is 
$(0,y,0,\mu)\in \Pi''$. Thus $\Pi''\subseteq \Pi$.
 Of course, if $(0,y,z,\mu)\notin\Pi''$ the assessment $(0,y,z,\mu)$ is not coherent and hence $(0,y,z,\mu)\notin\Pi$.
Therefore
$\Pi=\Pi'\cup \Pi'' $.
\end{proof}

\subsection{Proof of  Theorem \ref{THM:PIB}}
\label{SEC:APPB}
{\bf Theorem} \ref{THM:PIB}:
\emph{Let $A$, $B$, $H$, $K$ be any logically independent events. The set $\Pi$ of all the coherent assessment $(x,y,z,\mu)$ on the family $\F=\{A|H, B|K, (A|H)\wedge_{B}(B|K), (B|K)|_{B}(A|H)\}$ is $\Pi=\Pi' \cup \Pi''$, where 
$\Pi'=\{(x,y,z,\mu): x\in(0,1], y\in [0,1], z\in[z', z''], \mu=\frac{z}{x}\}$
 with $ z'=0$,    $z''=1$, and
$\Pi''=\{(0,y,0,\mu): (y,\mu)\in[0,1]^2\}.$}
\begin{proof}
Of course $(x,y)$ on $\{A|H, B|K\}$ is coherent for every $(x,y)\in [0,1]^2$. By Table (\ref{tab:LUB}), the assessment $z=P[(A|H)\wedge_{B}(B|K)]$ is a coherent extension of $(x,y)$ if and only if $z \in [z',z'']$ where $z'=0$ and $z''=1$. Assuming $x>0$, from (\ref{prevision}) it holds that $\mu=\frac{z}{x}$. Then, every $(x,y,z,\mu)\in \Pi'$ is coherent, i.e., $\Pi'\subseteq \Pi$. Of course, if $x>0$ and $(x,y,z,\mu)\notin \Pi'$, then the assessment $(x,y,z,\mu)$ is not coherent, i.e. $(x,y,z,\mu)\notin \Pi$.\\
\\
Let us consider now the case $x=0$.
We recall that it is coherent  to assess $(0,y,z)$, with $z>0$, on $\{A|H,B|K,(A|H)\wedge_B(B|K)\}$ (see  Table \ref{tab:LUB}). However, for the further object 
$(B|K)|_B(A|H)$ coherence  also requires that $z=\mu x$ (see \ref{prevision}). Then, 
coherence requires that $z=0$ when we consider the assessment $(0,y,z,\mu)$ on
$\{A|H, B|K,(A|H)\wedge_{B}(B|K),(B|K)|_{B}(A|H)\}$.
We show that the assessment $(0,y,0,\mu)$ is coherent if and only if  $(y,\mu)\in [0,1]^2$, that is $\Pi''\subseteq \Pi$.
As $x=0$, from (\ref{EQ:iteB'}), it holds that 
\begin{equation} \label{EQ:iteB'x=0}
\begin{array}{l}
(B|K)|_{B}(A|H) 
=\begin{cases}
    1,               & {\ \text{if} \ AHBK \ \text{is true}},\\
    0,               & {\ \text{if} \ AH\no{B}K\vee AH\no{K} \ \text{is true}},\\
    \mu_B,             & {\ \text{if} \ \no{H}\vee \no{A}HBK \lor \no{A}H\no{B}K\vee \no{A}H\no{K} \ \text{is true}},
    \end{cases}
\end{array}
\end{equation}
that is $(B|K)|_{B}(A|H)=BK|AH$ (see proof Theorem \ref{THM:PIK}.
As $x=0$ and $z=0$, it holds that  $(B|K)|_{B}(A|H)=BK|AH$. We show that the assessment $\P=(0,y,0,\mu)$ on $\F =\{A|H, B|K, (A|H)\wedge_{B}(B|K), (B|K)|_{B}(A|H)\}=\{A|H, B|K, AHBK|HK, BK|AH\}$ is coherent if and only if  $(y,\mu)\in [0,1]^2$, that is $\Pi''\subseteq \Pi$.
The constituents $C_h's$ and the points $Q_h's$ associated with $(\F, \P)$, where $\P=(0,y,0,\mu)$ are the following:
\[
 \begin{array}{l}
 C_1=AHBK, C_2=AH\no{B}K, C_3=AH\no{K}, C_4=\no{A}HBK, C_5=\no{A}H\no{B}K,\\ C_6=\no{A}H\no{K}, C_7=\no{H}BK, C_8=\no{H}\no{B}K, C_0=\no{H}\,\no{K},
 \end{array}
 \]
 and
 \[
 \begin{array}{l}
 Q_1=(1,1,1,1), Q_2=(1,0,0,0), Q_3=(1,y,0,0), Q_4=(0,1,0,\mu), Q_5=(0,0,0,\mu), \\
 Q_6=(0,y,0,\mu), Q_7=(0,1,0,\mu), Q_8=(0,0,0,\mu), \P=Q_0=(0,y,0,\mu).
 \end{array}
 \]
We observe that $C_1\vee \cdots \vee C_8=H\vee K$.
 The system $(\Sigma)$ in (\ref{SYST-SIGMA}) associated with the pair $(\F,\P)$  is\\
\begin{equation}
 \left\{
     \begin{array}{llll}
      \lambda_1+\lambda_2+\lambda_3=0, \\
      \lambda_1+y\lambda_3+\lambda_4+y\lambda_6+\lambda_7=y,\\
      \lambda_1=0,\\
      \lambda_1+\mu \lambda_4+\mu \lambda_5+\mu \lambda_6+\mu \lambda_7+\mu \lambda_8=\mu\\
     \lambda_1+\lambda_2+\lambda_3+\lambda_4+\lambda_5+\lambda_6+\lambda_7+\lambda_8=1, \\ \lambda_i\geq0 \ \forall i=1, \dots, 8.\\
     \end{array}
     \label{sB}
 \right.    
\end{equation}
We observe that $\P$ belongs to the segment with end points $Q_4$, $Q_5$; indeed $(0,y,0,\mu)=yQ_4+(1-y)Q_5=y(0,1,0,\mu)+(1-y)(0,0,0,\mu)$.\\
 Then the vector $\Lambda=(0,0,0,y,1-y,0,0,0)$ is a solution of (\ref{sB}), with
 \[
 \phi_{1}(\Lambda)=\sum_{h:C_h\subseteq H}\lambda_h=\lambda_1+\lambda_2+\lambda_3+\lambda_4+\lambda_5+\lambda_6=1>0,
 \]
 \[
 \phi_{2}(\Lambda)=\sum_{h:C_h\subseteq K}\lambda_h=\lambda_1+\lambda_2+\lambda_4+\lambda_5+\lambda_7+\lambda_8=1,
 \]
 \[
 \phi_{3}(\Lambda)=\sum_{h:C_h\subseteq HK}\lambda_h=\lambda_1+\lambda_2+\lambda_4+\lambda_5=1,
 \]
 \[
 \phi_{4}(\Lambda)=\sum_{h:C_h\subseteq AH}\lambda_h=\lambda_1+\lambda_2+\lambda_3=0.
 \]
Let $\mathcal{S}'=\{(0,0,0,y,1-y,0,0,0)\}$  denote a subset of the set $\mathcal{S}$ of all solutions of (\ref{EQ:SIGMAK}). We have that $M_1'=M_2'=M_3'=1$ and $M_4'=0$   (as defined in (\ref{EQ:I0'})). 
It follows that $I'_0=\{4\}$. As the sub-assessment $\P'_0=\mu$ on $\F'_0=\{BK|AH\}$ is coherent $\forall \mu \in [0,1]$,  by Theorem \ref{CNES-PREV-I_0-INT}, it follows that 
the assessment  $(0,y,0,\mu)$ on $\F=\{A|H, B|K, (A|H)\wedge_{B}(B|K), (B|K)|_{B}(A|H)\}$ is coherent for every $(y,\mu)\in[0,1]^2$, that is 
$(0,y,0,\mu)\in \Pi''$. Thus $\Pi''\subseteq \Pi$.
Of course, if $(0,y,z,\mu)\notin\Pi''$ the assessment $(0,y,z,\mu)$ is not coherent and hence $(0,y,z,\mu)\notin\Pi$.
Therefore
$\Pi=\Pi'\cup \Pi'' $.
 \end{proof}

\subsection{Proof of Theorem \ref{THM:LU_B}} \ \label{SEC:APPB2}\\
{\bf Theorem} \ref{THM:LU_B}.
	Let $A,B,H,K$ be any logically independent events. Given a coherent assessment $(x,y)$ on $\{A|H,B|K\}$, for the iterated conditional $(B|K)|_{B}(A|H)$ the  extension $\mu_B=\prev((B|K)|_{B}(A|H))$ is coherent if and only if 
	$\mu_B \in [\mu_B', \mu_B'']$, where	\begin{equation*}
	\begin{array}{ll}
	\mu_B'=0,\;\;
	\mu_B''=\left\{
	\begin{array}{ll}
	\frac{1}{x},& \text{ if } 0<x<1,\\
	1,& \text{ if } x=0 \vee x=1.\\
	\end{array}
	\right.
	\end{array}	\end{equation*} 
\begin{proof}
Assume that $x>0$. We simply write $\mu$ instead of $\mu_B$. From Theorem \ref{THM:PIB} it follows that the set of all coherent assessments $(x,y,z,\mu)$ on $\F$ is $\Pi'=\{(x,y,z,\mu): 0<x\leq 1, 0\leq y\leq 1, z' \leq z\leq z'', \mu=\frac{z}{x}\} $, where $z'=0$ and $z''=1$ (see Table \ref{tab:LUB}). Then,  $\mu$ is a coherent extension of $(x,y)$ if and only if $\mu \in [\mu', \mu'']$, where $\mu'=0$ and $\mu''=\frac{z''}{x}=\frac{1}{x}$.
Assume that  $x=0$. From Theorem \ref{THM:PIB} it follows that the set of all coherent assessments $(x,y,z,\mu)$ on $\F=\{A|H,B|K$, $(A|H)\wedge_B (B|K), (B|K)|_B(A|H)\}$ is  $\Pi''=\{(0,y,0,\mu): (y,\mu)\in[0,1]^2\}$. Then, 
$\mu$ is a coherent extension of $(x,y)$ if and only if $\mu \in[\mu',\mu'']$, where $\mu'=0$ and $\mu''=1$.
\end{proof}

\subsection{Proof of Theorem \ref{THM:PIS}}
\label{SEC:APPS} 
{\bf Theorem} \ref{THM:PIS}: \emph{
Let $A$, $B$, $H$, $K$ be any logically independent events. The set $\Pi$ of all the coherent assessment $(x,y,z,\mu)$ on the family $\F=\{A|H, B|K, (A|H)\wedge_{S}(B|K), (B|K)|_{S}(A|H)\}$ is $\Pi=\Pi' \cup \Pi''$, where 
\begin{equation*}
\begin{array}{ll}
\Pi'=\{(x,y,z,\mu): x\in(0,1], y\in [0,1], z\in[z', z''], \mu=\frac{z}{x}\},\\
 \text{with } z'=\max\{x+y-1, 0\},  
z''=\left\{
\begin{array}{ll}
	\frac{x+y-2xy}{1-xy}, &\text{ if } (x,y)\neq (1,1),\\
	1, &\text{ if } (x,y)= (1,1),\\
\end{array}
\right.
\end{array}
\end{equation*} 
 and
\begin{equation*}
\Pi''=\{(0,y,0,\mu): y\in[0,1], \mu \geq 0\}.
\end{equation*}
}
 \begin{proof}
We know that $(x,y)$ on $\{A|H, B|K\}$ is coherent for every $(x,y)\in [0,1]^2$. We also recall that   the assessment $z=P((A|H)\wedge_{S}(B|K))$ is a coherent extension of $(x,y)$ if and only if $z \in [z',z'']$, where (see Table \ref{tab:LUB})
\[
z'=\max\{x+y-1, 0\}\,\text{and } z''=\left\{
\begin{array}{ll}
	\frac{x+y-2xy}{1-xy}, &\text{ if } (x,y)\neq (1,1),\\
	1, &\text{ if } (x,y)= (1,1).\\
\end{array}
\right.
\]
Assuming $x>0$, then it follows from  (\ref{prevision}) that $\mu=\frac{z}{x}$. Thus, as every $(x,y,z,\mu)\in \Pi'$ is coherent, it holds that $\Pi'\subseteq \Pi$. Of course, if $x>0$ and $(x,y,z,\mu)\notin \Pi'$, then $(x,y,z,\mu)$ is not coherent, i.e.  $(x,y,z,\mu)\notin \Pi$. \\
Let us consider now the case $x=0$.  From Remark \ref{coherenceS}, coherence  requires  that $z=0$. We show that the assessment $\M=(0,y,0,\mu)$  is coherent if and only if  $y\in [0,1]$ and $\mu\geq 0$, and hence $\Pi''\subseteq \Pi$. As $x=0$, by (\ref{EQ:iteS'}) it holds that 
\begin{equation}\label{EQ:iteS'x=0}
\begin{array}{l}
(B|K)|_{S}(A|H)
 =
\left\{
\begin{array}{ll}
1, &\text{ if } AHBK \vee AH\no{K} \text{ is true},\\
0, &\text{ if }  AH\no{B}K \text{ is true},\\
\mu,& \text{ if } \no{A}H \vee \no{H}\,\no{K} \vee \no{H}\,\no{B}K \text{ is true},\\
1+\mu,& \text{ if } \no{H}BK  \text{ is true,}\\
\end{array}
\right.
\end{array}
\end{equation}
that is, when $x=0$,  
\[
(B|K)|_{S}(A|H) =[AHBK+AH\no{K}+(1+\mu)\no{H}BK]|(AH\vee \no{H}BK).
\]
Then, we have $\F=\{A|H,B|K,(AH\vee \no{H})\wedge(BK\vee \no{K})|(H\vee K),[AHBK+AH\no{K}+(1+\mu)\no{H}BK]|(AH\vee \no{H}BK)\}$, with 
$\H_4=H\vee K \vee (H\vee K)\vee (AH\vee \no{H}BK)=H\vee K$.
Therefore, the constituents and the points $Q_h$'s associated with $(\F,\M)$  are
\[
\begin{array}{l}
C_1=AHBK, C_2=AH\no{B}K, C_3=AH\no{K}, C_4=\no{A}HBK, C_5=\no{A}H\no{B}K,\\C_6=\no{A}H\no{K},  C_7=\no{H}BK, C_8=\no{H}\no{B}K, C_0=\no{H}\,\no{K},
\end{array}
\]
and
\[
\begin{array}{l}
Q_1=(1,1,1,1), Q_2=(1,0,0,0), Q_3=(1,y,1,1), Q_4=(0,1,0,\mu), Q_5=(0,0,0,\mu), \\ Q_6=(0,y,0,\mu)
Q_7=(0,1,1,1+\mu), Q_8=(0,0,0,\mu), \M=Q_0=(0,y,0,\mu).
\end{array}
\]
The system $(\Sigma)$ in (\ref{SYST-SIGMA}) associated with the pair $(\F,\M)$  is\\
\begin{equation}
\left\{
    \begin{array}{llll}
     \lambda_1+\lambda_2+\lambda_3=0, \\
     \lambda_1+y\lambda_3+\lambda_4+y\lambda_6+\lambda_7=y,\\
     \lambda_1+\lambda_3+\lambda_7=0,\\
     \lambda_1+\lambda_3+\mu \lambda_4+\mu \lambda_5+\mu \lambda_6+(1+\mu) \lambda_7+\mu \lambda_8=\mu,\\
    \lambda_1+\lambda_2+\lambda_3+\lambda_4+\lambda_5+\lambda_6+\lambda_7+\lambda_8=1, \\ \lambda_i\geq0 \ \forall i=1, \dots, 8.\\
    \end{array}
    \label{sS}
\right.    
\end{equation}
We observe that $\M$ belongs to the segment with end points $Q_4$, $Q_5$; indeed $(0,y,0,\mu)=yQ_4+(1-y)Q_5=y(0,1,0,\mu)+(1-y)(0,0,0,\mu)$.\\
Then, the vector $\Lambda=(0,0,0,y,1-y,0,0,0)$ is a solution of (\ref{sS}), with
\[
\phi_{1}(\Lambda)=\sum_{h:C_h\subseteq H}\lambda_h=\lambda_1+\lambda_2+\lambda_3+\lambda_4+\lambda_5+\lambda_6=1>0,
\]
\[
\phi_{2}(\Lambda)=\sum_{h:C_h\subseteq K}\lambda_h=\lambda_1+\lambda_2+\lambda_4+\lambda_5+\lambda_7+\lambda_8=1>0,
\]
 \[
 \phi_{3}(\Lambda)=\sum_{h:C_h\subseteq H\vee K}\lambda_h=\lambda_1+\lambda_2+\lambda_3+\lambda_4+\lambda_5+\lambda_6+\lambda_7+\lambda_8=1>0,
 \]
\[
\phi_{4}(\Lambda)=\sum_{h:C_h\subseteq (AH \vee \no{H}BK)}\lambda_h=\lambda_1+\lambda_2+\lambda_3+\lambda_7=0.
\]
Let $\mathcal{S}'=\{(0,0,0,y,1-y,0,0,0)\}$  denote a subset of the set $\mathcal{S}$ of all solutions of (\ref{sS}). We have that $M_1'=M_2'=M_3'=1$ and $M_4'=0$   (as defined in (\ref{EQ:I0'})). 
It follows that $I'_0=\{4\}$. We set  $\F'_0=\{[AHBK+AH\no{K}+(1+\mu)\no{H}BK]|(AH\vee \no{H}BK)\}$ and $\M'_0=\mu$. 
By exploiting the betting scheme, we show that the assessment $\mu$ on the conditional random quantity  $[AHBK+AH\no{K}+(1+\mu)\no{H}BK]|(AH\vee \no{H}BK)$ is coherent for every $\mu\geq 0$.
By recalling (\ref{EQ:RG}), the random gain for the assessment $\mu$ is 
\[
G=s(AH+\no{H}BK)(AHBK+AH\no{K}+(1+\mu)\no{H}BK-\mu).
\]
Without loss of generality we can assume $s=1$. 
The constituents in $AH\vee \no{H}BK$ are $C_1=AHBK, C_2=AH\no{B}K, C_3=AH\no{K}, C_4=\no{H}BK$ and the corresponding values for the random gain $G$ are $g_1=1-\mu$, $g_2=-\mu$, $g_3=1-\mu=g_1, g_4=1$. Hence the set of values of $G$ restricted to $AH\vee \no{H}BK$ is $\G_{AH\vee \no{H}BK}=\{g_1, g_2, g_4\}=\{1-\mu, -\mu, 1\}$. 
We distinguish two cases: $(i)$ $\mu<0$; $(ii)$ $\mu\geq 0$.\\
Case $(i)$. We observe that, $1-\mu>0$,  $-\mu>0$, and hence $\min \G_{AH\vee \no{H}BK}>0$. Then the assessment $\M_0$ on $\F_0$ is incoherent.\\
Case $(ii)$. We observe that, 
$\min \G_{AH\vee \no{H}BK}=-\mu\leq 0$ and 
$\max \G_{AH\vee \no{H}BK}=1>0$. Then, as
\[
\min \G_{AH\vee \no{H}BK}\cdot \max \G_{AH\vee \no{H}BK} \leq 0,
\]
the assessment $\M_0$ on $\F_0$ is coherent.
Then, as System (\ref{sS}) is solvable and $\M_0$ on $\F_0$ is coherent, 
by Theorem  \ref{CNES-PREV-I_0-INT'} 
it follows that 
every assessment $(0, y, 0, \mu)$ on $\F$ is coherent if and only if  $(0,y,0,\mu)\in \Pi''=\{(0,y,0,\mu): y\in[0,1], \mu \geq 0\}.$ Thus, $\Pi''\subseteq \Pi$. 
Of course, if  $(0,y,z,\mu)\notin\Pi''$ the assessment $(0,y,z,\mu)$ is not coherent and hence $(0,y,z,\mu)\notin\Pi$.
Therefore $\Pi=\Pi' \cup \Pi''$.
 \end{proof} 
\subsection{Proof of Theorem \ref{THM:LUS}}\label{SEC:APPS2}
{\bf Theorem}\ref{THM:LUS}:\emph{
	Let $A,B,H,K$ be any logically independent events. Given a coherent assessment $(x,y)$ on $\{A|H,B|K\}$, with $x\neq0$, for the iterated conditional $(B|K)|_{S}(A|H)$ the  extension $\mu_S=\prev((B|K)|_{S}(A|H))$ is coherent if and only if 
	$\mu_S \in [\mu_S', \mu_S'']$, where
\begin{equation*}
	\mu_S'=\max\{\frac{x+y-1}{x}, 0\} \text{ and }\mu_S''=\left\{
	\begin{array}{ll}
	\frac{x+y-2xy}{x(1-xy)},& \text{ if } (x,y)\neq (1,1);\\
	1,& \text{ if } (x,y)=(1,1).\\
	\end{array}
	\right.
\end{equation*} }
\begin{proof}
 We simply write $\mu$ instead of $\mu_S$. As, $x>0$, from Theorem \ref{THM:PIS} it follows that the set of all coherent assessments $(x,y,z,\mu)$ on 
$\F=\{A|H,B|K$, $(A|H)\wedge_S (B|K), (B|K)|_S(A|H)\}$
is the set $\Pi'$ given in Theorem \ref{THM:PIS}. Then,  $\mu$ is a coherent extension of $(x,y)$ if and only if $\mu \in [\mu', \mu'']$, where $\mu'$ and $\mu''$ are 
\[
\mu'=\frac{z'}{x}=\max\{\frac{x+y-1}{x}, 0\} \text{ and }\mu''=\frac{z''}{x}=\left\{
	\begin{array}{ll}
	\frac{x+y-2xy}{x(1-xy)},& \text{ if } (x,y)\neq (1,1);\\
	1,& \text{ if } (x,y)=(1,1).\\
	\end{array}
	\right.
\] 
\end{proof}

\section*{Acknowledgements}
We thank the anonymous reviewers for their  useful comments and suggestions.
Both authors acknowledge support by INdAM-GNAMPA research group. Giuseppe Sanfilippo acknowledges support by the FFR project of University of Palermo, Italy. 


\end{document}